\documentclass[3p,9pt,twocolumn]{elsarticle}

\usepackage[bitstream-charter]{mathdesign}
\usepackage[T1]{fontenc}

\usepackage{amsmath,amsthm,amsfonts,amscd}
\usepackage{rotating}
\usepackage{amsfonts}
\usepackage{amscd}
\usepackage[normalem]{ulem}
\usepackage{graphicx}
\usepackage[hang,small]{caption}
\usepackage[labelformat=simple]{subcaption}

\usepackage{wrapfig}
\usepackage{multirow}
\usepackage{xparse}
\usepackage{environ}
\usepackage{stackengine,wasysym}
\usepackage{mathdots}

\usepackage{tensor}
\usepackage{varwidth,calc}

\usepackage{pgfplotstable}
\usepackage{pgfplots}
\pgfplotsset{compat=1.14}
\pgfplotsset{
	colormap={my basis colormap}{
		rgb255=(0, 114, 189);
		rgb255=(54, 106, 148);
		rgb255=(108, 98, 107);
		rgb255=(163, 91, 66);
		rgb255=(217, 83, 25);
	}
}
\pgfplotsset{
	colormap={my parula}{
		rgb255=(53.0655, 42.406, 134.9460);
		rgb255=(20.2336, 132.6061, 211.9511);
		rgb255=(55.5463, 184.8859, 157.9118);
		rgb255=(208.7187, 186.8470,  89.1966);
		rgb255=(248.9565, 250.6905, 13.7190);
	}
}
\pdfpageattr {/Group << /S /Transparency /I true /CS /DeviceRGB>>}
\usepackage{tikz}
\usepgfmodule{oo}
\usetikzlibrary{calc}
\usetikzlibrary{shapes}
\usetikzlibrary{plotmarks}
\usetikzlibrary{backgrounds}
\usetikzlibrary{decorations.pathmorphing,decorations.pathreplacing,decorations.markings}
\usetikzlibrary{arrows.meta,arrows,fit,matrix,positioning,shapes.geometric,positioning}
\usetikzlibrary{external}
\usepackage{appendix}
\usepackage{chngcntr}
\usepackage{etoolbox}
\usepackage{apptools}
\AtAppendix{\counterwithin{theorem}{section}}
\AtBeginEnvironment{subappendices}{%
	\chapter*{Appendix}
	\addcontentsline{toc}{chapter}{Appendices}
}

\usepackage{breakcites}
\usepackage[inline]{enumitem}

\usepackage{scalerel,stackengine,wasysym}

\usepackage{accents}

\usepackage{algorithm}
\usepackage[noend]{algpseudocode}

\DeclareMathAlphabet{\mathpzc}{OT1}{pzc}{m}{it}

\usepackage{color}
\definecolor{lightgray}{gray}{0.80}
\definecolor{lightgray}{gray}{0.80}

\usepackage[framemethod=TikZ]{mdframed}

\usepackage[titles]{tocloft}
\cftsetindents{figure}{0em}{3.5em}
\cftsetindents{table}{0em}{3.5em}

\usepackage{xcolor,colortbl}
\usepackage{makecell}

\usepackage{titlesec}
\titleformat{\chapter}[block]
{\normalfont\huge\bfseries}{\thechapter.}{1em}{\huge}
\titlespacing*{\chapter}{0pt}{-19pt}{0pt}

\usepackage[T3,T1]{fontenc}

\newcommand{\Bezier}{B\'ezier~}

\theoremstyle{plain}
\newtheorem{theorem}{Theorem}[section]

\newtheorem{corollary}[theorem]{Corollary}

\theoremstyle{remark}
\newtheorem{remark}[theorem]{Remark}

\newtheorem{examplex}[theorem]{Example}
\newenvironment{example}
{\pushQED{\qed}\examplex}
{\popQED\endexamplex}

\theoremstyle{definition}
\newtheorem{definition}[theorem]{Definition}

\usepackage{mathtools}

\newcommand{\mbf}[1]{{\boldsymbol{#1}}}
\newcommand{\RR}{{\mathbb{R}}}

\newcommand{\NN}{{\mathbb{N}}}

\stackMath
\newcommand\reallywidehat[1]{%
	\savestack{\tmpbox}{\stretchto{%
			\scaleto{%
				\scalerel*[\widthof{\ensuremath{#1}}]{\kern-.6pt\bigwedge\kern-.6pt}%
				{\rule[-\textheight/2]{1ex}{\textheight}}
			}{\textheight}%
		}{0.5ex}}%
	\stackon[1pt]{#1}{\tmpbox}%
}
\stackMath
\newcommand\reallywidecheck[1]{%
	\savestack{\tmpbox}{\stretchto{%
			\scaleto{%
				\scalerel*[\widthof{\ensuremath{#1}}]{\kern-.6pt\bigvee\kern-.6pt}%
				{\rule[-\textheight/2]{1ex}{\textheight}}
			}{\textheight}%
		}{0.5ex}}%
	\stackon[1pt]{#1}{\tmpbox}%
}

\ExplSyntaxOn
\NewEnviron{bmatrixT}
{
	\marine_transpose:V \BODY
}

\int_new:N \l_marine_transpose_row_int
\int_new:N \l_marine_transpose_col_int
\seq_new:N \l_marine_transpose_rows_seq
\seq_new:N \l_marine_transpose_arow_seq
\prop_new:N \l_marine_transpose_matrix_prop
\tl_new:N \l_marine_transpose_last_tl
\tl_new:N \l_marine_transpose_body_tl

\cs_new_protected:Nn \marine_transpose:n
{
	\seq_set_split:Nnn \l_marine_transpose_rows_seq { \\ } { #1 }
	\int_zero:N \l_marine_transpose_row_int
	\prop_clear:N \l_marine_transpose_matrix_prop
	\seq_map_inline:Nn \l_marine_transpose_rows_seq
	{
		\int_incr:N \l_marine_transpose_row_int
		\int_zero:N \l_marine_transpose_col_int
		\seq_set_split:Nnn \l_marine_transpose_arow_seq { & } { ##1 }
		\seq_map_inline:Nn \l_marine_transpose_arow_seq
		{
			\int_incr:N \l_marine_transpose_col_int
			\prop_put:Nxn \l_marine_transpose_matrix_prop
			{
				\int_to_arabic:n { \l_marine_transpose_row_int }
				,
				\int_to_arabic:n { \l_marine_transpose_col_int }
			}
			{ ####1 }
		}
	}
	\tl_clear:N \l_marine_transpose_body_tl
	\int_step_inline:nnnn { 1 } { 1 } { \l_marine_transpose_col_int }
	{
		\int_step_inline:nnnn { 1 } { 1 } { \l_marine_transpose_row_int }
		{
			\tl_put_right:Nx \l_marine_transpose_body_tl
			{
				\prop_item:Nn \l_marine_transpose_matrix_prop { ####1,##1 }
				\int_compare:nF { ####1 = \l_marine_transpose_row_int } { & }
			}
		}
		\tl_put_right:Nn \l_marine_transpose_body_tl { \\ }
	}
	\begin{bmatrix*}[r]
		\l_marine_transpose_body_tl
	\end{bmatrix*}
}
\cs_generate_variant:Nn \marine_transpose:n { V }
\cs_generate_variant:Nn \prop_put:Nnn { Nx }
\ExplSyntaxOff

\makeatletter
\catcode`! 3
\def\Transpose #1{\romannumeral0\expandafter
	\Mar@Transpose@a\romannumeral`^^@\Mar@DoOneRow #1\\!\\}

\def\Mar@DoOneRow #1\\{\Mar@DoOneRow@a {}#1&^^@&}%

\def\Mar@DoOneRow@a #1#2&{%
	\if^^@\detokenize{#2}\expandafter\@gobble\fi
	\Mar@DoOneRow@a {#1#2\\}%
}%

\def\Mar@Transpose@a #1#2\\{\ifx!#2\expandafter\Mar@FinishTranspose\fi
	\expandafter\Mar@Transpose@b\romannumeral`^^@\Mar@DoOneRow@a {}#2&^^@&#1}

\def\Mar@Transpose@b #1#2^^@\\{\Mar@Join {}#2^^@!#1}

\def\Mar@Join #1#2\\#3!#4\\%
{\if^^@\detokenize{#3}\expandafter\Mar@EndJoin\fi
	\Mar@Join {#1#2&#4\\}#3!}%

\def\Mar@EndJoin\Mar@Join #1^^@!^^@\\{\Mar@Transpose@a {#1^^@\\}}

\def\Mar@FinishTranspose
#1&^^@&#2\\^^@\\{ #2}

\catcode`! 12
\makeatother

\newcommand{\standsimp}{\triangle}

\newcommand{\bcoord}{x}
\newcommand{\coord}{t}

\newcommand{\nel}{{\ncells{1}}}

\newcommand{\seg}[2]{{#2^{(#1)}}}
\newcommand{\ratseg}[2]{{#2^{w,(#1)}}}
\newcommand{\nseg}{{m}}

\newcommand{\interval}{I}
\newcommand{\domain}{\Omega}

\newcommand{\ncells}[1]{\mathfrak{t}_{#1}}

\newcommand{\degree}{p}

\newcommand{\bez}{b}
\newcommand{\ratbez}{\bez^w}
\newcommand{\bsp}{B}

\newcommand{\rwts}{w}
\newcommand{\brwts}{\mbf{\rwts}}
\newcommand{\brwtss}{\mbf{W}}

\newcommand{\mdbcoord}{\tau}
\newcommand{\knot}{\xi}
\newcommand{\bknot}{\mbf{\knot}}
\newcommand{\bknots}{\mbf{\Xi}}
\newcommand{\mmult}{m}
\newcommand{\nsum}{{\mu}}

\newcommand{\splSpace}{\mathcal{R}}

\newcommand{\ndof}{{n}}

\newcommand{\periodic}[1]{{#1}^{per}}

\newcommand{\refMat}{\mbf{R}}
\newcommand{\refMatloc}{\mbf{S}}
\newcommand{\extMat}{\mbf{H}}
\newcommand{\extMatel}{H}
\newcommand{\extMatint}{\mbf{H}^c}
\newcommand{\extMatintel}{H^c}
\newcommand{\extMatinv}{\mbf{G}}

\newcommand{\extMatintinv}{\mbf{G}^c}
\newcommand{\extMatintinvel}{G^c}
\newcommand{\extMatP}{\mbf{E}}
\newcommand{\extMatPel}{E}
\newcommand{\extMatPinv}{\mbf{D}}
\newcommand{\extMatPinvel}{D}
\newcommand{\reduced}[1]{\bar{#1}}
\newcommand{\refined}[1]{\tilde{#1}}

\newcommand{\polar}[1]{{#1}^{pol}}
\newcommand{\domainP}{\polar{\domain}}
\newcommand{\bspP}{\polar{\bsp}}

\newcommand{\polarMap}{F}
\newcommand{\bpolarMap}{\mbf{\polarMap}}
\newcommand{\bpolarMapC}{\bpolarMap}
\newcommand{\pcoorda}{s}
\newcommand{\pcoordb}{t}

\newcommand{\radius}{\rho}
\newcommand{\polari}{l}

\newcommand{\spd}{{d}}
\newcommand{\ptop}[1]{\hat{#1}}

\newcommand{\identMat}{\mbf{I}}
\newcommand{\exchangeMat}{\mbf{J}}

\definecolor{myBlue}{rgb} {0,0.4470,0.7410}
\definecolor{myRed}{rgb} {0.8500,0.3250,0.0980}
\definecolor{myGray1}{rgb} {0,0,0}
\definecolor{myGray2}{rgb} {0.6,0.6,0.6}
\definecolor{myGray3}{rgb} {0.45,0.45,0.45}
\definecolor{myGray4}{rgb} {0.8,0.8,0.8}
\definecolor{myGray5}{rgb} {0.55,0.55,0.55}

\definecolor{myCPlot1}{rgb} {0,    0.4470,    0.7410}
\definecolor{myCPlot2}{rgb} {0.8500,   0.3250,    0.0980}
\definecolor{myCPlot3}{rgb} {0.9290,    0.6940,    0.1250}
\definecolor{myCPlot4}{rgb} {0.4940,    0.1840,    0.5560}
\definecolor{myCPlot5}{rgb} {0.4660,    0.6740,    0.1880}
\definecolor{myCPlot6}{rgb} {0.3010,    0.7450,    0.9330}
\definecolor{myCPlot7}{rgb} {0.6350,    0.0780,    0.1840}

\definecolor{myCPColor}{rgb} {0.9255, 0.6941, 0.355}

\tikzstyle{scpColor}=[circle, fill=myCPColor]

\tikzset{
	show curve controls/.style={
		decoration={
			show path construction,
			curveto code={
				\draw[myGray1,dashed,eThickness]
				(\tikzinputsegmentfirst)
				-- (\tikzinputsegmentsupporta)
				-- (\tikzinputsegmentsupportb)
				-- (\tikzinputsegmentlast)
				;
				\fill[scpColor] (\tikzinputsegmentfirst) circle(3pt);
				\fill[scpColor] (\tikzinputsegmentsupporta) circle(3pt);
				\fill[scpColor] (\tikzinputsegmentsupportb) circle(3pt);
				\fill[scpColor] (\tikzinputsegmentlast) circle(3pt);
				\draw[#1,line width=1pt]
				(\tikzinputsegmentfirst)
				.. controls (\tikzinputsegmentsupporta)
				and (\tikzinputsegmentsupportb) ..
				(\tikzinputsegmentlast);
			}
		},decorate
	}
}

\tikzset{
	on each segment/.style={
		decorate,
		decoration={
			show path construction,
			moveto code={},
			lineto code={
				\path [#1]
				(\tikzinputsegmentfirst) -- (\tikzinputsegmentlast);
			},
			curveto code={
				\path [#1] (\tikzinputsegmentfirst)
				.. controls
				(\tikzinputsegmentsupporta) and (\tikzinputsegmentsupportb)
				..
				(\tikzinputsegmentlast);
			},
			closepath code={
				\path [#1]
				(\tikzinputsegmentfirst) -- (\tikzinputsegmentlast);
			},
		},
	},
	mid arrow/.style={postaction={decorate,decoration={
				markings,
				mark=at position .6 with {\arrow[#1]{Latex}}
	}}},
}

\tikzset{
	bThickness/.style={line width=#1\pgflinewidth},
	bThickness/.default={2},
}

\tikzset{
	eThickness/.style={line width=#1\pgflinewidth},
	eThickness/.default={0.5},
}

\pgfooclass{myAxis}{ 
	\method myAxis() { 
	}
	\method plot(#1,#2,#3,#4,#5,#6,#7) { 
		\draw[-latex,thick,myBlue] (#1,#2) --+ (#4:#3) node [right,sloped,black] {$#6$};
		\draw[-latex,thick,myRed] (#1,#2) --+ ({#4+#5}:#3) node [above,sloped,black] {$#7$};
	}
}
\pgfoonew \myAxis=new myAxis()

\tikzset{cross/.style={cross out, draw, 
		minimum size=2*(#1-\pgflinewidth), 
		inner sep=0pt, outer sep=0pt}}
	
\pgfooclass{myNewAxis}{ 
	\method myNewAxis() { 
	}
	\method plotX(#1,#2,#3,#4) { 
		\draw[-latex,thick,myBlue] (#1) --+ (#3:#2) node [right,sloped,black] {$#4$};
	}
	\method plotY(#1,#2,#3,#4,#5) {
		\draw[-latex,thick,myRed] (#1) --+ ({#3+#4}:#2) node [above,sloped,black] {$#5$};
	}
}
\pgfoonew \myNewAxis=new myNewAxis()

\pgfooclass{myBezierElements}{ 
	\method myBezierElements() { 
	}
	\method plot(#1,#2,#3,#4,#5,#6,#7,#8,#9) { 
		\begin{scope}[rotate=#6,xslant=cot(180-#5),xscale=1,yscale=sin(#5)]
			\draw[eThickness, step={#4/#1}, #7, #8, #9] (#2,#3) grid ({#2+#4},{#3+#4}); 
		\end{scope}
	}
	\method plotPoints(#1,#2,#3,#4,#5,#6,#7,#8,#9) { 
		\begin{scope}[rotate=#6,xslant=cot(180-#5),xscale=1,yscale=sin(#5)]
			\draw[step=#4, #7, black] (#2,#3) grid (#2+#4,#3+#4); 
			\def\dx{#4*1.0/3}
			\draw[step=\dx, #8, #9, black] (#2,#3) grid (#2+#4,#3+#4); 
			\foreach \x in {1,...,#1}{
				\foreach \y in {1,...,#1}{
					\node [draw, black, fill=white, shape=rectangle, minimum width=3pt, minimum height=3pt, anchor=center,scale=0.4] at (#2+\x*\dx-\dx,#3+\y*\dx-\dx) {};
				}
			}
		\end{scope}
	}
	\method plotGrid(#1,#2,#3,#4,#5,#6,#7,#8,#9) { 
		\begin{scope}[rotate=#6,xslant=cot(180-#5),xscale=1,yscale=sin(#5)]
			\draw[step=#4, #7, black] (#2,#3) grid (#2+#4,#3+#4); 
			\def\dx{#4/3}
			\draw[step=\dx, #8, #9, black] (#2,#3) grid (#2+#4,#3+#4); 
		\end{scope}
	}
	\method plotPointsWLabels(#1,#2,#3,#4,#5,#6) { 
		\begin{scope}[rotate=#6,xslant=cot(180-#5),xscale=1,yscale=sin(#5)]
			\draw[step=#4, thick, black] (#2,#3) grid (#2+#4,#3+#4); 
			\def\dx{#4/3}
			\draw[step=\dx, ultra thin, black] (#2,#3) grid (#2+#4,#3+#4); 
			\foreach \x in {1,...,#1}{
				\foreach \y in {1,...,#1}{
					\node [draw, myBlue, shape=rectangle, minimum width=3pt, minimum height=3pt, anchor=center] at (#2+\x*\dx-\dx,#3+\y*\dx-\dx) {
						\pgfmathparse{\x-1}
						\pgfmathprintnumber[    
						fixed,
						fixed zerofill,
						precision=0
						]{\pgfmathresult}%
						\pgfmathparse{\y-1}
						\pgfmathprintnumber[    
						fixed,
						fixed zerofill,
						precision=0
						]{\pgfmathresult}%
					};
				}
			}
		\end{scope}
	}
	\method plotVertices(#1,#2,#3,#4,#5,#6,#7,#8) { 
		\begin{scope}[rotate=#6,xslant=cot(180-#5),xscale=1,yscale=sin(#5)]
			\pgfmathsetmacro{\ncp}{#1+1}
			\foreach \x in {1,...,\ncp}{
				\foreach \y in {1,...,\ncp}{
					\node[draw,black,circle,fill=#7,scale=#8] at ({#2+(\x-1)*#4/#1},{#3+(\y-1)*#4/#1}) {};
				}
			}	
		\end{scope}
	}
	\method shade(#1,#2,#3,#4,#5,#6,#7) { 
		\begin{scope}[rotate=#6,xslant=cot(180-#5),xscale=1,yscale=sin(#5)]
			\fill[#7] (#2,#3) rectangle +(#4,#4);
		\end{scope}
	}
	\method plotC1Vertices(#1,#2,#3,#4,#5,#6,#7,#8) {
		\begin{scope}[rotate=#6,xslant=cot(180-#5),xscale=1,yscale=sin(#5)]
			\pgfmathsetmacro{\nel}{#1}
			\foreach \i in {1,...,\nel}{
				\foreach \j in {1,...,\nel}{
					\pgfmathsetmacro{\xl}{#2+(\i-1)*#4/#1}
					\pgfmathsetmacro{\xr}{#2+(\i)*#4/#1}
					\pgfmathsetmacro{\yt}{#3+(\j)*#4/#1}
					\pgfmathsetmacro{\yb}{#3+(\j-1)*#4/#1}
					\foreach \ii in {1,...,2}{
						\foreach \jj in {1,...,2}{
							\pgfmathsetmacro{\x}{(\xl*(3-\ii)+\xr*(\ii))/3}
							\pgfmathsetmacro{\y}{(\yb*(3-\jj)+\yt*(\jj))/3}
							\node[draw,black,regular polygon,regular polygon sides=4,fill=#7,scale=#8] at (\x,\y) {};
						}
					}
				}
			}
		\end{scope}
	}
}
\pgfoonew \myBezierElements=new myBezierElements()

\newcommand\pgfmathsinandcos[3]{%
	\pgfmathsetmacro#1{sin(#3)}%
	\pgfmathsetmacro#2{cos(#3)}%
}
\newcommand\LongitudePlane[3][current plane]{%
	\pgfmathsinandcos\sinEl\cosEl{#2} 
	\pgfmathsinandcos\sint\cost{#3} 
	\tikzset{#1/.style={cm={\cost,\sint*\sinEl,0,\cosEl,(0,0)}}}
}
\newcommand\LatitudePlane[3][current plane]{%
	\pgfmathsinandcos\sinEl\cosEl{#2} 
	\pgfmathsinandcos\sint\cost{#3} 
	\pgfmathsetmacro\yshift{\cosEl*\sint}
	\tikzset{#1/.style={cm={\cost,0,0,\cost*\sinEl,(0,\yshift)}}} %
}
\newcommand\DrawLongitudeCircle[2][2]{
	\LongitudePlane{\angEl}{#2}
	\tikzset{current plane/.prefix style={scale=#1}}
	\pgfmathsetmacro\angVis{atan(sin(#2)*cos(\angEl)/sin(\angEl))} %
	\draw[current plane, ultra thin] (\angVis:1) arc (\angVis:\angVis+180:1);
	\draw[current plane,dashed, ultra thin] (\angVis-180:1) arc (\angVis-180:\angVis:1);
}
\newcommand\DrawLongitudeCircleThick[2][2]{
	\LongitudePlane{\angEl}{#2}
	\tikzset{current plane/.prefix style={scale=#1}}
	\pgfmathsetmacro\angVis{atan(sin(#2)*cos(\angEl)/sin(\angEl))} %
	\draw[current plane, ultra thick] (\angVis-\angEl:1) arc (\angVis-\angEl:\angVis+180-\angEl:1);
	\draw[current plane, ultra thin] (\angVis:1) arc (\angVis+180-\angEl:\angVis+180:1);
	\draw[current plane,dashed, ultra thin] (\angVis-180:1) arc (\angVis-180:\angVis:1);
	\pgfmathsetmacro\H{#1*cos(\angEl)} 
	\coordinate[mark coordinate] (N) at (0,\H);
	\coordinate[mark coordinate] (N) at (0,-\H);
}
\newcommand\DrawLatitudeCircle[2][2]{
	\LatitudePlane{\angEl}{#2}
	\tikzset{current plane/.prefix style={scale=#1}}
	\pgfmathsetmacro\sinVis{sin(#2)/cos(#2)*sin(\angEl)/cos(\angEl)}
	\pgfmathsetmacro\angVis{asin(min(1,max(\sinVis,-1)))}
	\draw[current plane, ultra thin] (\angVis:1) arc (\angVis:-\angVis-180:1);
	\draw[current plane,dashed, ultra thin] (180-\angVis:1) arc (180-\angVis:\angVis:1);
}

\pgfooclass{myQuadElement}{ 
	\method myQuadElement() { 
	}
	\method plot(#1,#2,#3,#4,#5,#6,#7,#8) { 
		\begin{scope}
			\draw[eThickness,black] (#1,#2) -- (#3,#4) -- (#5,#6) -- (#7,#8) -- cycle;
			\node[coordinate] (v1) at (#1,#2) {};
			\node[coordinate] (v2) at (#3,#4) {};
			\node[coordinate] (v3) at (#5,#6) {};
			\node[coordinate] (v4) at (#7,#8) {};
		\end{scope}
	}
	\method plotColor(#1,#2,#3,#4,#5,#6,#7,#8,#9) { 
		\begin{scope}
			\draw[eThickness,dashed,#9] (#1,#2) -- (#3,#4) -- (#5,#6) -- (#7,#8) -- cycle;
			\node[coordinate] (v1) at (#1,#2) {};
			\node[coordinate] (v2) at (#3,#4) {};
			\node[coordinate] (v3) at (#5,#6) {};
			\node[coordinate] (v4) at (#7,#8) {};
		\end{scope}
	}
	\method plotInteriorDomainPoints(#1,#2,#3,#4,#5,#6) {
		\begin{scope}
			\pgfmathtruncatemacro{\p}{#5-1}
			\foreach \x in {2,...,\p}{
				\foreach \y in {2,...,\p}{
					\pgfmathsetmacro{\a}{(\x-1)/(#5-1)}
					\pgfmathsetmacro{\b}{(\y-1)/(#5-1)}
					\pgfmathsetmacro{\aone}{(1-\a)*(1-\b)}
					\pgfmathsetmacro{\atwo}{(\a)*(1-\b)}
					\pgfmathsetmacro{\athree}{(\a)*(\b)}
					\pgfmathsetmacro{\afour}{(1-\a)*(\b)}
					\node [draw, black, shape=rectangle, minimum width=3pt, minimum height=3pt, anchor=center,scale=0.5,fill=#6] (b\x\y) at (barycentric cs:#1=\aone,#2=\atwo,#3=\athree,#4=\afour) {};
				}
			}
		\end{scope}
	}
	\method plotEdgeDomainPoints(#1,#2,#3,#4,#5) {
		\begin{scope}
			\pgfmathtruncatemacro{\p}{#3-1}
			\foreach \x in {2,...,\p}{
				\pgfmathsetmacro{\a}{(\x-1)/(#3-1)}
				\pgfmathsetmacro{\aone}{(1-\a)}
				\pgfmathsetmacro{\atwo}{\a}
				\node [draw, black, shape=rectangle, minimum width=3pt, minimum height=3pt, anchor=center,scale=0.5,fill=#5] (e#4\x) at (barycentric cs:#1=\aone,#2=\atwo) {};
			}
		\end{scope}
	}
	\method plotCorner(#1,#2,#3) {
		\begin{scope}
			\node [draw, black, shape=rectangle, minimum width=3pt, minimum height=3pt, anchor=center,scale=0.5,fill=#3] (c#2) at (#1) {};
		\end{scope}
	}
}
\pgfoonew \myQuad=new myQuadElement()

\pgfplotstableset{
	create on use/hthree/.style={
		create col/copy column from table={/Users/deepesh/Documents/research/reports/utdiss2-05/from_dropbox/tikz/extraordinary_points/data/meshSize3.txt}{H}
	}
}
\pgfplotstableset{
	create on use/hfive/.style={
		create col/copy column from table={/Users/deepesh/Documents/research/reports/utdiss2-05/from_dropbox/tikz/extraordinary_points/data/meshSize5.txt}{H}
	}
}
\pgfplotstableset{
	create on use/hsix/.style={
		create col/copy column from table={/Users/deepesh/Documents/research/reports/utdiss2-05/from_dropbox/tikz/extraordinary_points/data/meshSize6.txt}{H}
	}
}
\pgfplotstableset{
	create on use/hseven/.style={
		create col/copy column from table={/Users/deepesh/Documents/research/reports/utdiss2-05/from_dropbox/tikz/extraordinary_points/data/meshSize7.txt}{H}
	}
}
\pgfplotstableset{
	create on use/hmep/.style={
		create col/copy column from table={/Users/deepesh/Documents/research/reports/utdiss2-05/from_dropbox/tikz/extraordinary_points/data/meshSizemep.txt}{H}
	}
}

\makeatletter
\def\ps@pprintTitle{%
	\let\@oddhead\@empty
	\let\@evenhead\@empty
	\def\@oddfoot{}%
	\let\@evenfoot\@oddfoot}
\makeatother

\graphicspath{{./images/}}
\usepackage{hyperref}

\newcommand{\DTA}{DTA-compatible}
\newcommand{\DT}[1]{{\color{black}#1}}
\newcommand{\HS}[1]{{\color{black}#1}}
\newcommand{\revision}[1]{{\color{black}#1}}

\begin{document}
	
\title{A general class of $C^1$ smooth rational splines:\\
Application to construction of exact ellipses and ellipsoids}

	\author[roma]{Hendrik Speleers}
	\ead{speleers@mat.uniroma2.it}
	
	\author[delft]{Deepesh Toshniwal\corref{cor1}}
	\ead{d.toshniwal@tudelft.nl}
	
	\address[roma]{Department of Mathematics, University of Rome Tor Vergata, Italy}
	\address[delft]{Delft Institute of Applied Mathematics, Delft University of Technology, The Netherlands}
	\cortext[cor1]{Corresponding author}
	
	\begin{abstract}
		In this paper, we describe a general class of $C^1$ smooth rational splines that enables, in particular, exact descriptions of ellipses and ellipsoids --- 
		\revision{some of the most important primitives for CAD and CAE.}
		The univariate rational splines are assembled by transforming multiple sets of NURBS basis functions via so-called design-through-analysis compatible extraction matrices; different sets of NURBS are allowed to have different polynomial degrees and weight functions.
		Tensor products of the univariate splines yield multivariate splines.
		In the bivariate setting, we describe how 
		\HS{similar design-through-analysis compatible transformations}
		of the tensor-product splines enable the construction of smooth surfaces containing one or two polar singularities.
		The material is self-contained, and is presented such that all tools can be easily implemented by CAD or CAE practitioners within existing software that support NURBS.
		To this end, we explicitly present the matrices (a) that describe our splines in terms of NURBS, and (b) that help refine the splines by performing (local) degree elevation and knot insertion.
		Finally, all $C^1$ spline constructions yield spline basis functions that are locally supported and form a convex partition of unity.
	\end{abstract}
	
	\begin{keyword}
		Piecewise-NURBS representations \sep Smooth parameterizations \sep Exact ellipses and ellipsoids
	\end{keyword}
	\maketitle
	
	\section{Introduction}

Multivariate splines are used extensively for computer-aided design (CAD) and, more recently, for computer-aided engineering (CAE).
Smoothness of such splines is a particularly valuable trait.
When the aim is to create a (freeform) geometric model for a smooth object, it helps if the splines used for the task are smooth themselves.
For instance, this circumvents situations where small displacements to control points may produce `non-smooth' features such as $C^0$ kinks, loss of curvature continuity, etc.
Similarly, when the aim is to numerically approximate the solution to high-order partial differential equations (PDEs) using isogeometric analysis (IGA)  --- a generalization of classical finite elements \DT{\cite{hughes2005isogeometric}} --- high smoothness of the approximating spaces can be beneficial.
For instance, it can allow us to directly discretize the PDEs without any auxiliary variables, thus yielding simpler and more efficient implementations.

In this paper, we discuss a general class of $C^1$ smooth rational splines that allow for the construction of $C^1$ smooth curves and surfaces.
These are an extension of classical $C^1$ non-uniform rational B-splines (NURBS) as they enjoy the flexibility of choosing locally unrelated weight functions as well as the option of local degree elevation --- they can \HS{be roughly} regarded as piecewise-NURBS. At the same time, they maintain intuitive control-point-based design.
Moreover, they enable simple (low-degree) and smooth descriptions of some of the most important primitives 
for CAD and CAE (but also for computer vision, graphics and robotics): closed, real, non-degenerate quadrics --- that is, ellipses in two dimensions and ellipsoids in three dimensions.

The ideas we present here build upon those from \cite{toshniwal2017multi}, in multiple directions, and their presentation is motivated by our primary objectives: \emph{self-contained, explicit, NURBS-compatible descriptions that can be easily and efficiently implemented within existing CAD software}.
The most important novel contributions are the following.
\begin{itemize}
	\item We describe the usage of classical univariate NURBS to assemble \emph{$C^1$ rational multi-degree spline basis functions} using an extraction matrix. \revision{
		The general framework was
		explained in \cite{toshniwal2017multi}, but we provide here a simplified exposition of the construction and a formal proof of the properties; see Remark~\ref{rmk:novelty}.} 
	We mainly stick to parametric smoothness, but a construction centred around the notion of geometric smoothness can be formulated as well; see Remark~\ref{rem:geometric-continuity}.
	\item We describe efficient refinement of the $C^1$ splines leveraging classical NURBS refinement. \revision{The novelty here relies in an explicit and simple construction of the refinement matrices.}
	\item We describe how tensor-product bivariate $C^1$ rational splines can be used to build $C^1$ smooth geometries that may contain one or two polar singularities; the $C^1$ smooth splines describing the geometries are called \emph{polar splines}. \revision{As above, the idea is based on building an extraction matrix.}
	\item We describe efficient refinement of polar splines. \revision{
		In particular, we provide an explicit and simple construction of the refinement matrices.}
	\item We provide explicit descriptions of ellipses and ellipsoids built using low-degree $C^1$ splines, and we detail their extraction in terms of NURBS so that they can be readily implemented and used in CAD or CAE software.
	Table~\ref{tab:conic-sections} summarizes the descriptions included in this paper.
\end{itemize}

\begin{table*}[t!]
	\centering
	\begin{tabular}{|m{3.5cm}|>{\raggedleft\arraybackslash}m{3.5cm}|>{\raggedleft\arraybackslash}m{1.7cm}|>{\raggedleft\arraybackslash}m{1.4cm}|>{\raggedleft\arraybackslash}m{2cm}|>{\raggedleft\arraybackslash}m{1.4cm}|}
		\hline\hline
		\rowcolor{gray!10}
		\thead{\textbf{Quadric}\\~} & \thead[r]{\textbf{Polynomial}\\\textbf{degree}} & \thead[r]{\textbf{\# rational}\\\textbf{pieces}} & \thead[r]{\textbf{\# DOFs}\\~} & \thead[r]{\textbf{\# DOFs for}\\\textbf{$C^{-1}$ NURBS}} &
		\thead[r]{\textbf{Section}\\~} \\
		\hline\hline
		\cellcolor[gray]{0.9} & uniform; $2$ & 4 & \multirow{3}{*}{4} & 12 & \ref{ssec:ellipse2}\\
		\cline{2-3}\cline{5-6}
		\cellcolor[gray]{0.9} & uniform; $3$ & 2 & & 8 & \ref{ssec:ellipse3}\\
		\cline{2-3}\cline{5-6}
		\multirow{-3}{*}{\cellcolor[gray]{0.9}\thead{Ellipse \\(special case: circle)}} & non-uniform; $(3,2,2)$ & 3 & & 10 & \ref{ssec:ellipse322}\\
		\hline
		\hline
		\cellcolor[gray]{0.9}
		& uniform; $(2,2)$ & 8 & \multirow{3}{*}{6} & 72 & \ref{ssec:ellipsoid22}\\
		\cline{2-3}\cline{5-6}
		\cellcolor[gray]{0.9} & uniform; $(2,3)$ & 4 & & 48 & \ref{ssec:ellipsoid23}\\
		\cline{2-3}\cline{5-6}
		\multirow{-3}{*}{\cellcolor[gray]{0.9}\thead{Ellipsoid\\(special case: sphere)}} & uniform; $(3,3)$ & 2 & & 32 & \ref{ssec:ellipsoid33}\\
		\hline\hline
	\end{tabular}
	\caption{\DT{An overview of the explicit $C^1$ descriptions of quadrics presented in this paper.
	The table also compares the number of degrees of freedom (DOFs) needed by our $C^1$ representation compared to those needed by an equivalent $C^{-1}$ NURBS representation.}\label{tab:conic-sections}}
\end{table*}

\subsection{Extraction matrices}
At the core of our approach is the notion of the so-called \emph{design-through-analysis (DTA) compatible extraction matrix}.\footnote{Our notion of \emph{\DTA{} extraction matrix} has been called \emph{IGA-suitable extraction matrix} in \cite{toshniwal2017multi}. The name reflects the fact that both design and analysis may profit \DT{from} the extraction operation.}
Roughly speaking, such matrix helps us assemble `simple splines' into `more general splines.'
Examples are the \Bezier extraction matrix introduced to assemble Bernstein polynomials into B-/T-splines \cite{borden2011isogeometric,scott2011extraction}; the multi-degree extraction matrix for assembling elements of extended Tchebycheff spaces into generalized Tchebycheffian B-splines \cite{speleers2019computation,toshniwal2020multi,hiemstra2020generalized}; and the unstructured spline extraction matrices for assembling tensor-product splines into splines on unstructured quadrilateral \revision{meshes \cite{toshniwal2017multi,toshniwal2017smooth,toshniwal2020discretedifferential}.}

Here, we apply the concept of extraction in the following context.
We start from multiple sets of (univariate or bivariate) NURBS basis functions defined on adjacent domains, and collect all of these functions in the set $\{b_{j}:j=1,\dots,m\}$.
Then, we assemble them into more general $C^1$ rational (polar) splines using a matrix $\mbf{C}$ (with entries $C_{ij}$), called the \HS{\emph{extraction matrix}.}
Denote this new set of splines by $\{N_{i}:i=1,\dots,\ndof\}$, where $\ndof<m$. These are defined as follows,
\begin{equation}\label{eq:basis_expansion}
   N_i = \sum_{j=1}^{m}C_{ij}b_{j},\quad i = 1, \dots, \ndof.
\end{equation}
We are particularly interested in matrices $\mbf{C}$ such that the functions $N_i$ 
\begin{itemize}
   \item satisfy certain smoothness constraints that may or may not be satisfied by the $b_j$, and
   \item possess the properties of non-negativity, locality, linear independence and partition of unity that the $b_j$ already possess.
\end{itemize}
Such extraction matrices are called \emph{\DTA}.
\begin{definition}[\DTA{} extraction] An extraction matrix $\mbf{C}$ is called \DTA{} if
   \begin{enumerate}[label=(\alph*)]
      \item $\mbf{C}$ is a full-rank matrix,
      \item each column of $\mbf{C}$ sums to 1,
      \item each entry in $\mbf{C}$ is non-negative, and
      \item $\mbf{C}$ imparts locality to the functions $N_i$ through sparsity.
   \end{enumerate}
   \label{def:dasExtraction}
\end{definition}
It is easy to see that the action of a \DTA{} extraction matrix on a convex partition of unity, local basis gives rise to another local basis that also forms a convex partition of unity.
Indeed, by summing over $i$ in Equation \eqref{eq:basis_expansion}, we have
\begin{equation*} 
	\sum_{i=1}^{\ndof} N_i = \sum_{i=1}^{\ndof}\sum_{j=1}^{m}{C}_{ij}b_{j}
	= \sum_{j=1}^{m}b_{j}\sum_{i=1}^{\ndof}{C}_{ij} = 1,
\end{equation*}
as the $b_{j}$ form a partition of unity.
Since $\mbf{C}$ has non-negative entries and is a full-rank matrix, non-negativity and linear independence of $N_i$ follow from the non-negativity and linear independence of~$b_{j}$.

\subsection{Related literature}
As mentioned in the \HS{previous section}, the construction of smooth univariate splines by joining simpler pieces has been recently explored in \cite{toshniwal2017multi,speleers2019computation,toshniwal2020multi} for  polynomial multi-degree splines, and in \cite{hiemstra2020generalized} for generalized Tchebycheffian splines.
These approaches have conceptual similarities with the notion of beta-splines \cite{barsky1988beta}.
The main differences are that the former approaches do not rely on symbolic computations while the latter \HS{does}, and the former approaches consider parametric continuity while the latter studies geometric continuity.
\revision{The use of smooth univariate rational splines for construction of circles has been previously explored in \cite{bangert1997circle,lu2009circular,lu2011cylindrical}.
It is known that a circle cannot be represented by a single (symmetric) periodic $C^1$ quadratic NURBS curve \cite[Section~7.5]{piegl2012nurbs} nor a $C^2$ cubic NURBS curve \cite[Section~13.7]{farin2002cagd}. However, it is possible to find $C^p$ smooth descriptions using NURBS of degree $2(p+1)$, which is shown to be the minimal degree in \cite{bangert1997circle}. On the other hand, \cite{lu2009circular,lu2011cylindrical} presented a $C^1$ piecewise quadratic NURBS description of the circle and used it for IGA. Our rational multi-degree splines form a flexible extension of the latter framework, and allow for a variety of exact descriptions of circles using low \mbox{(multi-)}degrees, as indicated in Table~\ref{tab:conic-sections}.}

In two dimensions, closed quadrics or, more generally, smooth closed surfaces of genus zero can be built using tensor-product splines by introducing polar singularities.
For such polar surfaces, subdivision schemes producing $C^1$ surfaces \cite{karvciauskas2007bicubic,myles2008} and $C^2$ surfaces \cite{karvciauskas2006c2polarjet,myles2009c2polar} have been previously worked out.
The corresponding limit surfaces consist of an infinite sequence of surface rings where the faces shrink to a point in the limit.
A more CAD-friendly finite construction was developed in \cite{myles2011c2}; this approach constructs `shape' basis functions for $C^2$ polar splines with bi-degree $(6,3)$.
These basis functions correspond to unique Fourier frequencies in the polar expansion of a quadratic surface.
The `shape' basis does not enjoy non-negativity and does not form a partition of unity, and extensions of it to higher smoothness leads to degrees of freedom that control non-intuitive shape parameters. \revision{Similar recent constructions for obtaining $C^1$ polar spline caps can be found in \cite{karvciauskas2020polar}.}
Curvature continuous polar NURBS surfaces were discussed in \cite{shi2013polar}, and \cite{shi2011gn} presented a construction of polar caps using periodic B-spline surfaces with $G^n$ continuity for arbitrary $n$.
On the CAE side, a standard circular serendipity-type element for IGA was proposed in \cite{lu2009circular}, and $C^k$ smooth basis functions over singular parametrizations of triangular domains were constructed in \cite{takacs2015construction}.
A design-through-analysis friendly construction of $C^k$ smooth polar surfaces was recently proposed in \cite{toshniwal2017multi}, and the current work builds further upon this construction.

\revision{A completely different approach for dealing with curves and surfaces is the use of implicit representations \cite{Gomes2009implicit}. Such representations enjoy nice geometric properties (especially for simple shapes). For instance, they allow for a straightforward point membership classification. On the other hand, explicit smooth B-spline representations are more convenient for direct geometric modeling and (local) modification.}

\subsection{Outline}

In Section~\ref{sec:rational-curves}, we present the construction and refinement of $C^1$ smooth rational multi-degree spline curves via explicitly defined extraction and refinement matrices.
The construction and refinement of $C^1$ smooth polar surfaces using tensor products of the univariate splines is detailed in Section~\ref{sec:rational-surfaces}.
The explicit descriptions of ellipses and ellipsoids using the univariate and bivariate $C^1$ splines are reported in Sections~\ref{sec:ellipse} and~\ref{sec:ellipsoid}, respectively.
Finally, we conclude the paper in Section~\ref{sec:conclusions}.
It should be mentioned that the text is written such that the theoretical sections --- \HS{Sections~\ref{sec:nurbs}--\ref{sec:refinement-curves} and \ref{sec:polar_map}--\ref{sec:refinement-surfaces}} --- can be skipped by readers interested only in implementing explicit descriptions of smooth quadrics.

	\section{Piecewise-rational curves} \label{sec:rational-curves}

In this section, we focus on a multi-degree extension of univariate NURBS splines. 
The multi-degree spline space is defined as a collection of classical NURBS spaces (with possibly \HS{different polynomial degrees and weight functions}) glued together $C^1$ smoothly. For such space we present a construction of a set of basis functions, with similar properties to classical NURBS. After discussing some preliminary material on NURBS in Section~\ref{sec:nurbs}, we elaborate how these basis functions can be computed through a \DTA{} extraction matrix in Section~\ref{sec:mdb-splines}.
A more general but also more complex algorithmic construction has been detailed in \cite[Section~2]{toshniwal2017multi} and further explored in \cite{speleers2019computation,toshniwal2020multi} for polynomial multi-degree splines. Then, in Section~\ref{sec:refinement-curves}, we give an explicit procedure how to compute a refined representation of a given curve. Finally, in Section~\ref{sec:ellipse}, we illustrate how this tool can be used to describe arbitrary ellipses in a $C^1$ smooth fashion using low-degree piecewise-rational curve representations suited for integrated design and analysis.

\subsection{Preliminaries on NURBS} \label{sec:nurbs}

We start by defining notation for NURBS basis functions, and introduce some classical relations that can be found, e.g., in \cite{rogers2001,piegl2012nurbs}.

Given a basic interval $\interval := [\bcoord_1, \bcoord_2] \subset \RR$, let us denote with $\bknot$ an \emph{open knot vector} of degree $\degree\in\NN$ and length $\ndof+\degree+1\in\NN$, i.e.,
\begin{equation} \label{eq:knotVector}
\begin{gathered}
  \bknot := [\knot_1,\knot_2,\dots,\knot_{\ndof+\degree+1}],\quad 
  \knot_{i+1} \geq \knot_{i},\\
  \knot_1 = \dots = \knot_{\degree+1} = \bcoord_1 < \knot_{\degree+2},\\
  \knot_{\ndof+1} = \dots = \knot_{\ndof+\degree+1} = \bcoord_2 > \knot_{\ndof}.
\end{gathered}
\end{equation}
The number of times a knot value $\knot_i$ is duplicated in the knot vector is called the knot's \emph{multiplicity}. The multiplicity of $\knot_i$ is denoted with $\mmult_i$, and we assume that $1\leq \mmult_i \leq \degree-1$.
The corresponding set of B-splines $\{\bez_{j,\degree}:j=1,\dots,\ndof\}$ are defined using the recursive relation,
\begin{equation*}
  \bez_{j,p}(\bcoord) := \frac{\bcoord-\knot_j}{\knot_{j+\degree}-\knot_j}\bez_{j,\degree-1}(\bcoord) + \frac{\knot_{j+\degree+1}-\bcoord}{\knot_{j+\degree+1}-\knot_{j+1}}\bez_{j+1,\degree-1}(\bcoord),
\end{equation*}
starting from
\begin{equation*}
  \bez_{j,0}(\bcoord) := 
  \begin{cases}
    1, & \text{if }\knot_j \leq \bcoord < \knot_{j+1},\\
    0, & \text{otherwise},
  \end{cases}
\end{equation*}
and under the convention that fractions with zero denominator have value zero. With the above definition, all the B-splines take the value zero at the end point $\bcoord_2$. 
Therefore, in order to avoid asymmetry over the interval $\interval$, it is common to assume the B-splines to be left continuous at $\bcoord_2$. We will follow suit.

Let us denote with $\brwts$ a \emph{weight vector} of length $\ndof$, i.e.,
\begin{equation}\label{eq:weights}
  \brwts:=[\rwts_1,\rwts_2,\dots,\rwts_\ndof], \quad \rwts_i>0.
\end{equation}
The corresponding set of NURBS $\{\ratbez_{j,\degree}:j=1,\dots,\ndof\}$ are defined by
\begin{equation*}
\ratbez_{j,\degree}(\bcoord) := \frac{w_j \bez_{j,\degree}(\bcoord)}{\sum_{i=1}^\ndof w_i \bez_{i,\degree}(\bcoord)}.
\end{equation*}
Each $\ratbez_{j,p}$ is non-negative on $\interval$ and is locally supported on $[\knot_j,\knot_{j+p+1}]$. Moreover, the functions $\ratbez_{j,\degree}$ are linearly independent and form a partition of unity. They satisfy the following end-point conditions:
\begin{align*} 
   &\ratbez_{1,\degree}(\bcoord_1)=1, \quad \ratbez_{j,\degree}(\bcoord_1)=0, \quad j=2,\ldots,\ndof, \\
   &\ratbez_{\ndof,\degree}(\bcoord_2)=1, \quad \ratbez_{j,\degree}(\bcoord_2)=0, \quad j=1,\ldots,\ndof-1.
\end{align*}

The NURBS space corresponding to $\bknot$ and $\brwts$ is denoted with $\splSpace[\bknot,\brwts]$ and is defined as the span of $\{\ratbez_{j,\degree}:j=1,\dots,\ndof\}$. This is a space of piecewise-rational functions of degree $\degree$ with smoothness $C^{\degree-\mmult_i}$ at knot $\knot_i$ and its dimension is $\ndof$. The assumption on the multiplicity will ensure us global $C^1$ smoothness.
Note that when $\rwts_1=\dots=\rwts_\ndof$, the \DT{members of this space are} piecewise-polynomial.

\begin{remark}
The structure of $\bknot$ in Equation \eqref{eq:knotVector} is such that $\degree$, $\mmult_i$, $\ndof$ and $\interval$ are embedded in it. Therefore, we will assume that once a knot vector $\bknot$ is known, so are the degree, smoothness, and dimension of the corresponding NURBS space $\splSpace[\bknot,\brwts]$.
\end{remark}

We identify a function $f\in \splSpace[\bknot,\brwts]$ with the vector of its coefficients $\left[f_1,\dots,f_\ndof\right]$, 
\begin{equation*}
f(\bcoord) = \sum_{j=1}^{\ndof} f_j \ratbez_{j,\degree}(\bcoord).
\end{equation*}
Only the first (last) $k+1$ basis functions contribute towards the $k$-th order derivative at the left (right) end point of $\interval$. In particular, we have
\begin{equation}\label{eq:C1-end}
\begin{aligned}
  f(\bcoord_1)&=f_1, \quad 
  \frac{d f}{d\bcoord}(\bcoord_1)= \frac{\degree}{\knot_{\degree+2}-\bcoord_1}\frac{\rwts_2}{\rwts_1}(f_2-f_1),\\
  f(\bcoord_2)&=f_\ndof, \quad  
  \frac{d f}{d\bcoord}(\bcoord_2)= \frac{\degree}{\bcoord_2-\knot_{\ndof}}\frac{\rwts_{\ndof-1}}{\rwts_\ndof}(f_\ndof-f_{\ndof-1}).
\end{aligned}
\end{equation}
A NURBS curve embedded in $\RR^\spd$, $\spd\geq2$, can be constructed as
\begin{equation*}
  \mbf{f}(\bcoord) = \sum_{j=1}^\ndof \mbf{f}_{j} \ratbez_{j,\degree}(\bcoord),
\end{equation*}
where $\mbf{f}_j \in \RR^{\spd}$ are the control points assigned to each basis function. All coordinate functions of this curve belong to $\splSpace[\bknot,\brwts]$ and therefore all the above relations hold for them.

\subsection{Rational multi-degree B-splines} \label{sec:mdb-splines}

\begin{figure}[t!]
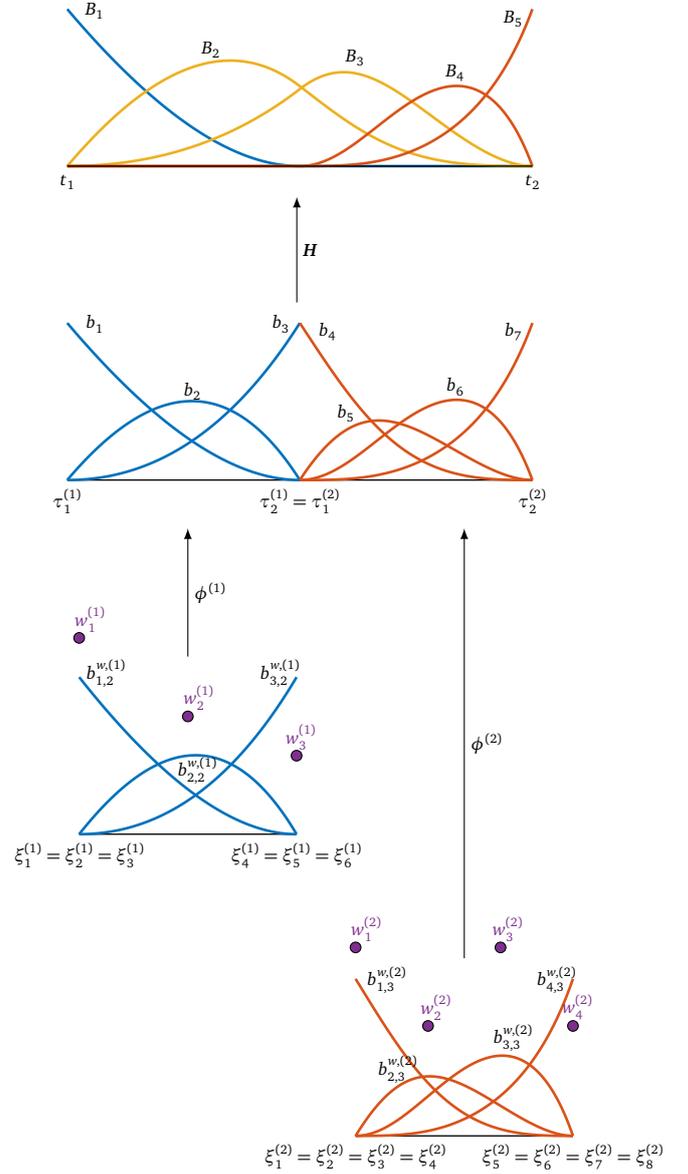

   \vspace*{-1cm}
	\centering
	\begin{tikzpicture}[anchor=south west]
		\node[] (n1) at (0,0) {%
	\tikzsetnextfilename{./tikz/images/mds_example_nurbs_1}%
	\input{./tikz/scripts/mds_example_nurbs_1}%
};
		\node[] (n2) at (3.3,-4.0) {%
	\tikzsetnextfilename{./tikz/images/mds_example_nurbs_2}%
	\input{./tikz/scripts/mds_example_nurbs_2}%
};
		\node[] (mds1) at (0.4,4.7) {%
	\tikzsetnextfilename{./tikz/images/mds_example_mds_c0}%
	\input{./tikz/scripts/mds_example_mds_c0}%
};
		\node[] (mds2) at (0.4,9) {%
	\tikzsetnextfilename{./tikz/images/mds_example_mds_c1}%
	\input{./tikz/scripts/mds_example_mds_c1}%
};
		\draw[-latex] ($(n1.north)+(0,-1)$) -- node[midway,right,scale=0.7] {$\phi^{(1)}$} ++(0,1.7);
		\draw[-latex] ($(n2.north)+(0,-1)$) -- node[midway,right,scale=0.7] {$\phi^{(2)}$} ++(0,5.7);
		\draw[-latex] ($(mds1.north)+(0,-1)$) -- node[midway,right,scale=0.7] {$\mbf{H}$} ++(0,1.4);
	\end{tikzpicture}
	\caption{A visual illustration of the notation and the construction of rational multi-degree B-splines as described in Section~\ref{sec:mdb-splines}.
	Here, the quadratic (blue) and cubic (red) NURBS shown at the bottom are used to build the $C^1$ multi-degree B-splines shown at the top.}
	\label{fig:mds_example}
\end{figure}

Consider $\nseg$ open knot vectors $\seg{i}{\bknot}$ of degree $\seg{i}{\degree}$, $i=1,\dots,\nseg$, defined as in Equation \eqref{eq:knotVector}. We denote the left and right end points of the interval $\seg{i}{\interval}$ associated to $\seg{i}{\bknot}$ with $\seg{i}{\bcoord_1}$ and $\seg{i}{\bcoord_2}$, respectively. 
The collection $\bknots := (\seg{1}{\bknot},\dots,\seg{\nseg}{\bknot})$ is called an \emph{$\nseg$-segment knot vector configuration}. 
The multi-degree spline spaces will be constructed by considering spline spaces over the knot vectors $\seg{i}{\bknot}$, which are glued together with certain smoothness requirements at the end points $\seg{i}{\bcoord_2}$ and $\seg{i+1}{\bcoord_1}$ for $i \in \{1,2,\dots,\nseg-1\}$. The equivalence class at the points $\seg{i}{\bcoord_2}$ and $\seg{i+1}{\bcoord_1}$ is called the \emph{$i$-th segment join}. 
We define the mapping $\seg{i}{\phi}$ for each segment $i=1,\ldots,\nseg$,
\begin{equation}\label{eq:phi}
   \seg{i}{\phi}(\bcoord) := \bcoord-\seg{i}{\bcoord_1} + \seg{1}{\mdbcoord_1} + \sum_{\ell=1}^{i-1}(\seg{\ell}{\bcoord_2}-\seg{\ell}{\bcoord_1}),
\end{equation}
for an arbitrarily chosen origin $\seg{1}{\mdbcoord_1} \in \RR$.
Then, $\seg{i}{\domain} := [\seg{i}{\mdbcoord_1},\seg{i}{\mdbcoord_2}] := \seg{i}{\phi}( [\seg{i}{\bcoord_1},\seg{i}{\bcoord_2}]) \subset \RR$, and we construct the composed interval
\begin{equation*}
   \domain := [\coord_1,\coord_2] := \seg{1}{\domain}\cup\dots\cup\seg{\nseg}{\domain}.
\end{equation*}
Note that $\seg{i}{\mdbcoord_2} = \seg{i+1}{\mdbcoord_1}$, $i \in\{ 1,\ldots,\nseg-1\}$ and $\coord_1=\seg{1}{\mdbcoord_1}$.
Moreover, let $\brwtss := (\seg{1}{\brwts},\ldots,\seg{\nseg}{\brwts})$ be a sequence of weight vectors defined as in Equation \eqref{eq:weights}.
\revision{We refer the reader to Figure~\ref{fig:mds_example} for a visual illustration of the notation of the above concepts, in case $\nseg=2$ and $\seg{1}{\degree}=2$, $\seg{2}{\degree}=3$.}

The space of rational multi-degree splines is defined as
\begin{equation*}
\begin{split}
	\splSpace[\bknots,\brwtss] := \bigl\{\,f \in C^1(\domain): \;
	&f\circ \seg{i}{\phi} \in \splSpace[\seg{i}{\bknot},\seg{i}{\brwts}],\\
	&\quad 1 \leq i \leq \nseg\bigr\},
\end{split}
\end{equation*}
and the periodic space of rational multi-degree splines as
\begin{equation*}
\begin{split}
	\periodic{\splSpace}[\bknots,\brwtss] := \bigl\{\,f \in \splSpace[\bknots,\brwtss]: \;
	&f(\coord_1)=f(\coord_2),\\
	&\frac{d f}{d\coord}(\coord_1)=\frac{d f}{d\coord}(\coord_2)\,\bigr\}.
\end{split}
\end{equation*}
The elements of $\splSpace[\bknots,\brwtss]$ and $\periodic{\splSpace}[\bknots,\brwtss]$ are piecewise-NURBS functions such that the pieces meet with $C^1$ continuity at each segment join. 
It is clear that classical NURBS spaces are a special case of the rational multi-degree spline spaces.

In the following, we build a suitable basis for the spaces $\splSpace[\bknots,\brwtss]$ and $\periodic{\splSpace}[\bknots,\brwtss]$.
On the $i$-th knot vector $\seg{i}{\bknot}$, we have $\seg{i}{\ndof}$ NURBS $\ratseg{i}{\bez_{j,\seg{i}{\degree}}}$ of degree $\seg{i}{\degree}$ that span the spline space $\splSpace[\seg{i}{\bknot},\seg{i}{\brwts}]$. 
In the first step, we map these basis functions from $\seg{i}{\interval}$ to $\seg{i}{\domain}$ using $\seg{i}{\phi}$ in Equation \eqref{eq:phi}, and extend them on the entire interval $\domain$ by defining them to be zero outside $\seg{i}{\domain}$.
More precisely, specifying the cumulative local dimensions $\nsum_i$ for $i=0,\ldots,\nseg$, 
\begin{equation*}
  \nsum_{0}:=0, \quad 
  \nsum_i:=\sum_{\ell=1}^{i}\seg{\ell}{\ndof}=\nsum_{i-1}+\seg{i}{\ndof}, \quad i>0,
\end{equation*}
we define for $i=1,\ldots,\nseg$ and $j=1,\ldots,\seg{i}{\ndof}$,
\begin{equation*}
  \bez_{\nsum_{i-1} + j}(\coord) := 
  \begin{cases}
  \ratseg{i}{\bez_{j,\seg{i}{\degree}}}(\bcoord),& 
  \text{if } [\seg{i}{\mdbcoord_1},\seg{i}{\mdbcoord_2}) \ni \coord = \seg{i}{\phi}(\bcoord),\\[0.2cm]
  \ratseg{i}{\bez_{j,\seg{i}{\degree}}}(\seg{i}{\bcoord_2}),& 
  \text{if } i=\nseg \text{ and } \coord={\coord_2},\\[0.2cm]
  0, & \text{otherwise}.
  \end{cases}
\end{equation*}
For the sake of simplicity, we dropped the reference to the (local) degree and weight in the notation.
From the properties of NURBS, it is clear that the functions $\bez_{1},\ldots,\bez_{\nsum_{\nseg}}$ are linearly independent and
form a non-negative partition of unity \DT{on $\domain$}.
We arrange these basis functions in a column vector $\mbf{\bez}$ of length $\nsum_{\nseg}$.
\revision{We refer the reader again to Figure~\ref{fig:mds_example} for a visual illustration of the notation of the above concepts.}

Now, we construct extraction matrices $\extMat$ and $\periodic{\extMat}$ such that the functions in $\{\bsp_i:i=1,\dots,\ndof\}$ and $\{\periodic{\bsp}_i:i=1,\dots,\periodic{\ndof}\}$, defined by
\begin{equation} \label{eq:mdbExtraction}
\mbf{\bsp} := \extMat\mbf{\bez},\quad
\periodic{\mbf{\bsp}} := \periodic{\extMat}\mbf{\bez},
\end{equation}
span $\splSpace[\bknots,\brwtss]$ and $\periodic{\splSpace}[\bknots,\brwtss]$, respectively.
The key here, and the reason our approach can be efficiently implemented by design, is that these extraction matrices can be explicitly specified.
To this end, we define counters $\eta_i$ for $i = 0, \dots, \nseg$,
\begin{equation*}
   \eta_{0}:=0, \quad 
   \eta_i:=\sum_{\ell=1}^{i}(\seg{\ell}{\ndof}-2)=
   \nsum_i-2i,
   \quad i>0,
\end{equation*}
and parameters $\seg{i}{\alpha}$ and $\seg{i}{\beta}$ for $i = 1, \dots, \nseg-1$,
\begin{align*}
\begin{gathered}
	\seg{i}{\alpha} := \frac{\seg{i}{\degree}}{\seg{i}{\bcoord_2}-\seg{i}{\knot_{\seg{i}{\ndof}}}}\frac{\seg{i}{\rwts_{\seg{i}{\ndof}-1}}}{\seg{i}{\rwts_\seg{i}{\ndof}}} > 0,\\
	\seg{i}{\beta} := \frac{\seg{i+1}{\degree}}{\seg{i+1}{\knot_{\seg{i+1}{\degree}+2}}-\seg{i+1}{\bcoord_1}}\frac{\seg{i+1}{\rwts_2}}{\seg{i+1}{\rwts_1}} > 0.
\end{gathered}
\end{align*}
\HS{In the periodic setting, $\seg{\nseg}{\alpha}$ and $\seg{\nseg}{\beta}$ are computed using the above equations by identifying the index $i+1$ with $1$.
Recall Equation \eqref{eq:C1-end} to see the motivation behind the definition of the above parameters.}
Then, we define a common sparse matrix $\extMatint$ of size $\eta_\nseg \times (\nsum_\nseg-2)$, whose non-zero entries $\extMatintel_{ij}$ are identified as follows: for $i=1,\ldots,\nseg$ and $j=1,\ldots,\seg{i}{\ndof}-2$,
\begin{equation}\label{eq:mdbC1Extraction:a}
   \extMatintel_{\eta_{i-1}+j, \nsum_{i-1}+j} := 1,
\end{equation}
and for $i=1,\ldots,\nseg-1$,
\begin{equation} \label{eq:mdbC1Extraction:b}
\begin{gathered}
	\extMatintel_{\eta_{i}, \nsum_{i}-1} := \extMatintel_{\eta_{i}, \nsum_{i}} := \dfrac{\seg{i}{\alpha}}{\seg{i}{\alpha} + \seg{i}{\beta}},\\
	\extMatintel_{\eta_{i}+1, \nsum_{i}-1} := \extMatintel_{\eta_{i}+1, \nsum_{i}} := \dfrac{\seg{i}{\beta}}{\seg{i}{\alpha} + \seg{i}{\beta}}.
\end{gathered}
\end{equation}
The desired extraction matrices in Equation \eqref{eq:mdbExtraction} are then specified as \DT{follows:
\begin{equation} \label{eq:mdbC1Extraction}
\begin{gathered}
	\extMat := 
	\left[
	\begin{array}{c|c|c}
	1 & \mbf{0} & 0\\
	\mbf{0} & \extMatint & \mbf{0}\\
	0 & \mbf{0} & 1
	\end{array}
	\right],\\
	\periodic{\extMat} := 
	\left[
	\begin{array}{c|c|c}
	\frac{\seg{\nseg}{\beta}}{\seg{\nseg}{\alpha} + \seg{\nseg}{\beta}}  & \multirow{3}{*}{$\extMatint$} &  \frac{\seg{\nseg}{\beta}}{\seg{\nseg}{\alpha} + \seg{\nseg}{\beta}}\\
	\mbf{0} & & \mbf{0}\\
	\frac{\seg{\nseg}{\alpha}}{\seg{\nseg}{\alpha} + \seg{\nseg}{\beta}} & & \frac{\seg{\nseg}{\alpha}}{\seg{\nseg}{\alpha} + \seg{\nseg}{\beta}}
	\end{array}
	\right].
\end{gathered}
\end{equation}
The} number of rows in the two matrices are denoted with $\ndof := \eta_\nseg+2$ and $\periodic{\ndof} := \eta_\nseg$, respectively.
The sparse and simple structure of both matrices means that it is easy to verify that both have full rank.
Indeed, this conclusion can be directly deduced from the full rank of $\extMatint$, which in turn is implied by Equation \eqref{eq:mdbC1Extraction:a}.
Moreover, their entries are non-negative, and the column sum is equal to one. Hence, we conclude that these matrices are \DTA.

How these matrices help us build $C^1$ splines can be understood by taking into account Equation \eqref{eq:C1-end} at any end point $\seg{i}{\mdbcoord_2}$.
Let us discuss the non-periodic setting and fix $i\in\{ 1, \dots, \nseg-1\}$; the argument can be directly applied to the periodic setting \HS{as well}.
Only four functions ($\bez_{\nsum_{i}-1}$, $\bez_{\nsum_{i}}$, $\bez_{\nsum_{i}+1}$, $\bez_{\nsum_{i}+2}$) are not $C^{1}$ at $\seg{i}{\mdbcoord_2}$.
More precisely, only these four functions have non-vanishing values and first derivatives here.
In view of \eqref{eq:C1-end}, a spline $f$ given by
\begin{equation*}
   f(\coord)=\sum_{j=1}^{4} f_j \bez_{\nsum_{i}-2+j}(\coord)
\end{equation*}
will be $C^1$ \HS{at $\seg{i}{\mdbcoord_2}$ if}
\begin{equation}\label{eq:C1-condition}
f_2 = f_3, \quad
\seg{i}{\alpha} (f_2-f_1) =
\seg{i}{\beta} (f_4-f_3).
\end{equation}
We can verify that the entries of $\extMat$ satisfy exactly such relations.
Indeed, for some $j$, the matrix $\extMat$ defines two new functions $\bsp_j$ and $\bsp_{j+1}$ such that
\begin{equation*}
   \begin{bmatrix}
      \bsp_j\\
      \bsp_{j+1}
   \end{bmatrix}
   = 
   \underbrace{
   \begin{bmatrix*}[c]
   1 & \frac{\seg{i}{\alpha}}{\seg{i}{\alpha}+\seg{i}{\beta}} & \frac{\seg{i}{\alpha}}{\seg{i}{\alpha}+\seg{i}{\beta}} & 0 \\
   0 & \frac{\seg{i}{\beta}}{\seg{i}{\alpha}+\seg{i}{\beta}} & \frac{\seg{i}{\beta}}{\seg{i}{\alpha}+\seg{i}{\beta}} & 1
   \end{bmatrix*}}_{=:\, \seg{i}{\reduced{\extMat}}}
   \begin{bmatrix}
      \bez_{\nsum_{i}-1}\\
      \bez_{\nsum_{i}}\\
      \bez_{\nsum_{i}+1}\\
      \bez_{\nsum_{i}+2}
   \end{bmatrix}
   + \begin{bmatrix}
   \seg{i}{\reduced{\bsp}_j}\\
   \seg{i}{\reduced{\bsp}_{j+1}}
   \end{bmatrix},
\end{equation*}
where $\seg{i}{\reduced{\bsp}_j}$ and $\seg{i}{\reduced{\bsp}_{j+1}}$ are at least $C^{1}$ at $\seg{i}{\mdbcoord_2}$.
When setting $[f_1, f_2, f_3, f_4]$ equal to the first or the second row of $\seg{i}{\reduced{\extMat}}$, we see that Equation \eqref{eq:C1-condition} is satisfied.
The following result follows from the above discussion.

\begin{theorem}\label{thm:univariate_dta_compatibility}
  The $C^1$ smooth piecewise-rational functions in the sets  $\{\bsp_i:i=1,\dots,\ndof\}$ and $\{\periodic{\bsp}_i:i=1,\dots,\periodic{\ndof}\}$ are linearly independent, locally supported, and form a convex partition of unity on~$\domain$.
\end{theorem}

\begin{remark} \label{rmk:novelty}
  \revision{In \cite[Section~2.3.5]{toshniwal2017multi} it was observed that the $C^1$ smooth piecewise-rational basis functions enjoy the properties described in Theorem~\ref{thm:univariate_dta_compatibility}. However, a formal proof was missing. It was also pointed out that the property of non-negativity is in general not present in case of $C^2$ or higher smoothness. On the other hand, this is possible when restricting to polynomial pieces \cite{toshniwal2020multi}.}
\end{remark}

\begin{remark} \label{rmk:simple-betas}
  With the aim of designing quadric curves, also called conics, it is natural to choose local NURBS spaces of the same degree $p$ and defined on the same \DT{uniform} knot vector $\bknot$.
  Moreover, it is common to set $\seg{i}{\rwts_1}=\seg{i}{\rwts_{\seg{i}{\ndof}}}=1$. Under these circumstances, the ratios in Equation~\eqref{eq:mdbC1Extraction:b} read as
\begin{equation*}
\begin{gathered}
	 \frac{\seg{i}{\alpha}}{\seg{i}{\alpha}+\seg{i}{\beta}} = \DT{\frac{\seg{i}{\rwts_{\seg{i}{\ndof}-1}}}{\seg{i}{\rwts_{\seg{i}{\ndof}-1}}+\seg{i+1}{\rwts_2}},}\\
	\frac{\seg{i}{\beta}}{\seg{i}{\alpha}+\seg{i}{\beta}} = \frac{\seg{i+1}{\rwts_2}}{\seg{i}{\rwts_{\seg{i}{\ndof}-1}}+\seg{i+1}{\rwts_2}}.
\end{gathered}
\end{equation*}
Finally, if there is additional symmetry in the choice of weights, so $\seg{i}{\rwts_{\seg{i}{\ndof}-1}}=\seg{i+1}{\rwts_2}$, we simply get
\begin{equation*}
 \frac{\seg{i}{\alpha}}{\seg{i}{\alpha}+\seg{i}{\beta}} = \frac{1}{2}, \quad
 \frac{\seg{i}{\beta}}{\seg{i}{\alpha}+\seg{i}{\beta}} = \frac{1}{2}.
\end{equation*}
\end{remark}

Once we have computed a \DTA{} extraction matrix $\extMat$ (or $\periodic{\extMat}$), given $\ndof$ control points $\mbf{f}_i \in \RR^{\spd}$, $\spd\geq2$, we can construct a piecewise-rational curve $\mbf{f}$ embedded in $\RR^{\spd}$,
\begin{equation*}
 \mbf{f}(\coord) = \sum_{i=1}^{\ndof} \mbf{f}_i \bsp_i(\coord).
\end{equation*}
For a fixed curve, the transpose of $\extMat$ (or $\periodic{\extMat}$) defines the relationship between control points of the $\bez_j$ (discontinuous at the segment joins) and control points of the smooth $\bsp_i$. More precisely, if
\begin{equation*}
  \sum_{j=1}^{\nsum_{\nseg}}\mbf{g}_{j}\bez_j(\coord) = \mbf{f}(\coord) = \sum_{i=1}^{\ndof} \mbf{f}_i \bsp_i(\coord),
\end{equation*}
then
\begin{equation} \label{eq:mdbControlPointExtraction}
  \mbf{g}_{j} = \sum_{i=1}^{\ndof} \extMatel_{ij}\mbf{f}_i.
\end{equation}

\begin{remark}\label{rem:geometric-continuity}
  When dealing with curves, the proposed piecewise-NURBS framework can also be formulated in the context of geometric continuity \cite{dyn1989}. In such case, the $C^1$ smoothness condition at the segment join in \eqref{eq:C1-condition} is replaced by the $G^1$ smoothness condition
\begin{equation*}
  f_2 = f_3, \quad
  \seg{i}{\alpha} (f_2-f_1) = \seg{i}{\gamma}
  \seg{i}{\beta} (f_4-f_3),
\end{equation*}
for a given geometric shape parameter $\seg{i}{\gamma}>0$, resulting in the matrix
\begin{equation*} 
  \seg{i}{\reduced{\extMat}}=\begin{bmatrix*}[c]
   1 & \frac{\seg{i}{\alpha}}{\seg{i}{\alpha}+\seg{i}{\gamma}\seg{i}{\beta}} & \frac{\seg{i}{\alpha}}{\seg{i}{\alpha}+\seg{i}{\gamma}\seg{i}{\beta}} & 0 \\
   0 & \frac{\seg{i}{\gamma}\seg{i}{\beta}}{\seg{i}{\alpha}+\seg{i}{\gamma}\seg{i}{\beta}} & \frac{\seg{i}{\gamma}\seg{i}{\beta}}{\seg{i}{\alpha}+\seg{i}{\gamma}\seg{i}{\beta}} & 1
  \end{bmatrix*}.
\end{equation*}
It is clear that this matrix is still \DTA.
\end{remark}

\subsection{Refinement of piecewise-rational curves} \label{sec:refinement-curves}
The rational spline spaces defined in the previous section can be refined in a multitude of ways.
We could reduce the smoothness at segment joins, raise the polynomial degrees of local NURBS spaces, and/or insert new knots in local NURBS spaces \cite[Section~2.4.3]{toshniwal2017multi}. A combination of these possibilities could be judiciously employed to achieve spline spaces that provide higher resolution or approximation power exactly where needed. In this section, we present an explicit construction of refined representations of a given piecewise-rational curve.

Before delving into the details of the refinement procedure, we first define two matrices $\extMatinv$ and $\periodic{\extMatinv}$ that can be regarded as right inverses of the extraction matrices $\extMat$ and $\periodic{\extMat}$, respectively.
Looking at the structure of the matrix $\extMatint$ specified in Equations \eqref{eq:mdbC1Extraction:a}--\eqref{eq:mdbC1Extraction:b}, we can define a sparse matrix $\extMatintinv$ of size $(\nsum_\nseg-2) \times \eta_\nseg$, whose non-zero entries $\extMatintinvel_{ij}$ are identified as follows: for $i=1,\ldots,\nseg$ and $j=1,\ldots,\seg{i}{\ndof}-2$,
\begin{equation*}
   \extMatintinvel_{\nsum_{i-1}+j,\eta_{i-1}+j} := 1.
\end{equation*}
From its construction it is clear that the product $\extMatint\extMatintinv$ is equal to the identity matrix.
Similarly, keeping in mind Equation \eqref{eq:mdbC1Extraction}, the matrices
\begin{equation} \label{eq:mdbC1Extraction-inv}
   \extMatinv := 
   \left[
   \begin{array}{ccc}
   1 & \mbf{0} & 0\\ \hline
   \mbf{0} & \extMatintinv & \mbf{0}\\ \hline
   0 & \mbf{0} & 1
   \end{array}
   \right], \quad
   \periodic{\extMatinv} := 
   \left[
   \begin{array}{ccc}
   \mbf{0}\\ \hline
   \extMatintinv\\ \hline
   \mbf{0}
   \end{array}
   \right]
\end{equation}
give rise to products $\extMat\extMatinv$ and $\periodic{\extMat}\periodic{\extMatinv}$ that are equal to identity matrices.

Now, let $\splSpace$ be a given spline space and let us denote the target refined space with $\refined{\splSpace}$.
\DT{For simplicity of notation, we drop the superscript $per$ in case of periodicity.}
Then, we consider the two unique representations of a curve $\mbf{f}$ with coordinate functions in $\splSpace \subset \refined{\splSpace}$,
\begin{equation} \label{eq:two-representations}
  \sum_{i=1}^{\refined{\ndof}} \refined{\mbf{f}}_i \refined{\bsp}_i(\coord) = \mbf{f}(\coord) = \sum_{i=1}^{\ndof} \mbf{f}_i \bsp_i(\coord).
\end{equation}
Let us collect the control points in the row vectors $\refined{\mbf{F}}:=[\refined{\mbf{f}}_1,\dots,\refined{\mbf{f}}_{\refined{\ndof}}]$ and $\mbf{F}:=[\mbf{f}_1,\dots,\mbf{f}_{\ndof}]$. We now seek the refinement matrix $\refMat$ of size $\ndof \times \refined{\ndof}$ that helps us compute $\refined{\mbf{F}}$ from $\mbf{F}$, i.e.,
\begin{equation*}
  \refined{\mbf{F}}
  = 
  \mbf{F} \refMat.
\end{equation*}
Assume that $\extMat$ and $\refined{\extMat}$ are the extraction matrices corresponding to the spaces $\splSpace$ and $\refined{\splSpace}$, respectively.
Incorporating these matrices in the representations in Equation \eqref{eq:two-representations} results in
\begin{equation*}
\refined{\mbf{F}}\refined{\extMat}\refined{\mbf{\bez}}
= 
\mbf{F}{\extMat}{\mbf{\bez}},
\end{equation*}
where the column vectors ${\mbf{\bez}}$ and $\refined{\mbf{\bez}}$ collect the local NURBS basis functions.
Let the matrix $\refMatloc$ be such that ${\mbf{\bez}} = \refMatloc\refined{\mbf{\bez}}$; this matrix can be computed with standard NURBS refinement techniques.
Then, we have
\begin{equation*}
\refined{\mbf{F}}\refined{\extMat}\refined{\mbf{\bez}}
= 
\mbf{F}{\extMat}\refMatloc\refined{\mbf{\bez}}.
\end{equation*}
This implies that we can compute $\refMat$ by solving the following (overdetermined) linear system with a unique solution,
\begin{equation*} 
  \refMat\refined{\extMat} = \extMat\refMatloc.
\end{equation*}
After multiplication of both sides of this system with the matrix $\refined{\extMatinv}$ (corresponding to $\refined{\extMat}$) as defined in Equation \eqref{eq:mdbC1Extraction-inv}, we arrive at
\begin{equation} \label{eq:ref-system}
  \refMat = \extMat\refMatloc\refined{\extMatinv}.
\end{equation}
Note that the application of $\refined{\extMatinv}$ in Equation \eqref{eq:ref-system} means that a subset of $\refined{\ndof}$ columns of $\extMat\refMatloc$ are selected to form $\refMat$.

\begin{remark} \label{rmk:choice-G}
  The definition of $\refined{\extMatinv}$ is done for the sake of simplicity of computation of $\refMat$ in Equation \eqref{eq:ref-system}, but is not unique. Any matrix that is a right inverse of $\refined{\extMat}$ would be a valid choice as well, such as the standard Moore--Penrose right inverse $\refined{\extMat}^T(\refined{\extMat}\refined{\extMat}^T)^{-1}$.
\end{remark}

\subsection{Circles and ellipses} \label{sec:ellipse}

We now present the general construction of ellipses (and as a special case also circles) using the $C^1$ rational splines introduced thus far.
We present three approaches for doing so using splines of low(est) degree, i.e., $C^1$ splines of quadratic degree, cubic degree and mixed quadratic/cubic multi-degree.
All approaches will construct four $C^1$ piecewise-NURBS functions $\bsp_i$ and associated control points $\mbf{f}_i$, $i = 1, \dots, 4$, such that the curve~$\mbf{f}$,
\begin{equation*}
	\mbf{f}(\coord) := (f_x(\coord), f_y(\coord)) := \sum_{i = 1}^4 \mbf{f}_i \bsp_i(\coord),
\end{equation*}
describes the exact ellipse centred at $(0,0)$ and with axis lengths $(a_x, a_y)$,
\begin{equation}\label{eq:ellipse}
\left(\frac{f_x}{a_x}\right)^2
+ \left(\frac{f_y}{a_y}\right)^2
= 1.
\end{equation}
Since the splines $\bsp_i$ form a partition of unity, these ellipses can be affinely transformed by directly applying the transformation to the control points $\mbf{f}_i$.
Subdivided or higher-degree representations can be easily obtained by refining the 
\DT{representations provided here} (see Section~\ref{sec:refinement-curves}).

To visually illustrate the smoothness of $\mbf{f}$, we will also show the curve $\tilde{\mbf{f}}$ obtained by perturbing one of the control points. Since all $\bsp_i$ are smooth, the perturbed curve will also be smooth. For uniformity throughout the examples, we will choose the control points of the perturbed curve as
\begin{equation} \label{eq:conic2-perturb}
\tilde{\mbf{f}}_i := 
\begin{cases}
   \mbf{f}_i + (0,a_y), & i = 1,\\
   \mbf{f}_i, & i = 2, 3, 4.
\end{cases}
\end{equation}

\begin{sidewaysfigure*}[p]
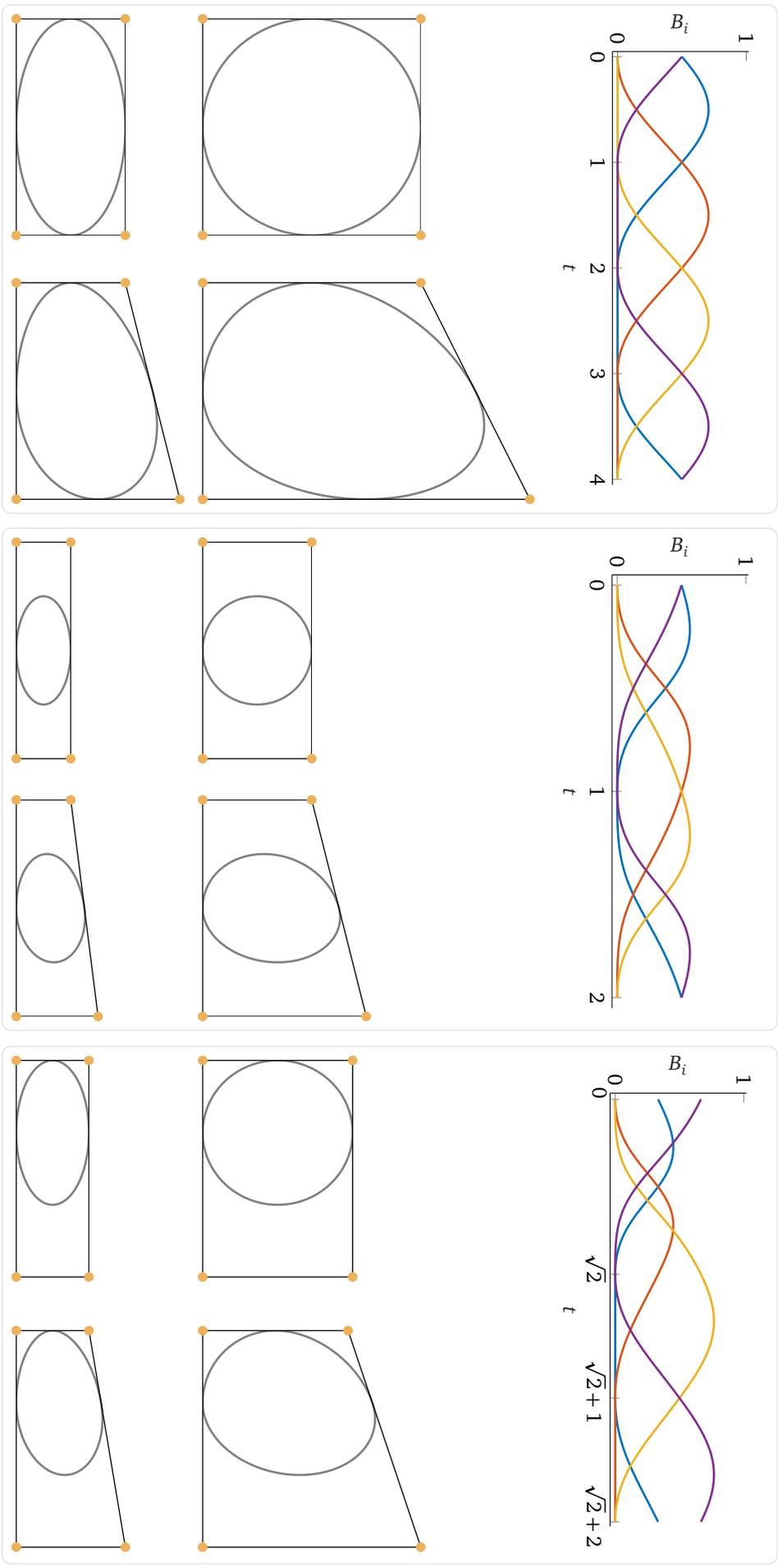

   \centering
   \begin{tikzpicture}[anchor=center]
   \matrix[matrix of nodes,column sep=0.4cm] (B) {
      \node[] (b1) {%
	\tikzsetnextfilename{./tikz/images/circle1_b}%
	\input{./tikz/scripts/circle1_b}%
}; & \node[] (b2) {%
	\tikzsetnextfilename{./tikz/images/circle2_b}%
	\input{./tikz/scripts/circle2_b}%
}; & \node[] (b3) {%
	\tikzsetnextfilename{./tikz/images/circle3_b}%
	\input{./tikz/scripts/circle3_b}%
};\\		
   };
   \matrix[matrix of nodes,column sep=0.4cm,anchor=north west] (ce1) at (b1.south west) {
      \node[] (c1) {%
	\tikzsetnextfilename{./tikz/images/circle1}%
	\input{./tikz/scripts/circle1}%
}; &
      \node[] (c1p) {%
	\tikzsetnextfilename{./tikz/images/circle1_perturb}%
	\input{./tikz/scripts/circle1_perturb}%
};\\
      \node[] (e1) {%
	\tikzsetnextfilename{./tikz/images/ellipse1}%
	\input{./tikz/scripts/ellipse1}%
}; &
      \node[] (e1p) {%
	\tikzsetnextfilename{./tikz/images/ellipse1_perturb}%
	\input{./tikz/scripts/ellipse1_perturb}%
};\\
   };
   \coordinate (a) at (ce1.south -| b2.west);
   \matrix[matrix of nodes,row sep=0.06\textwidth,column sep=0.3cm,anchor=south west] (ce2) at (a) {
      \node[] (c2) {%
	\tikzsetnextfilename{./tikz/images/circle2}%
%
%
\begin{tikzpicture}

\begin{axis}[%
ticks=none,
width=0.16\textwidth,
height=0.08\textwidth,
at={(0\textwidth,0\textwidth)},
scale only axis,
xmin=-2,
xmax=2,
ymin=-1,
ymax=1,
axis background/.style={fill=white},
axis line style={draw=none},
tick style={draw=none}
]
\addplot [color=gray, line width=1.0pt]
  table[row sep=crcr]{%
0	1\\
0.0201999591920017	0.999795960008162\\
0.0407993338884264	0.999167360532889\\
0.0617965597791463	0.998088766192398\\
0.0831889081455806	0.996533795493934\\
0.104972375690608	0.994475138121547\\
0.127141568981064	0.991884580703336\\
0.149689583812371	0.988733042078639\\
0.172607879924953	0.984990619136961\\
0.195886151638364	0.980626644343459\\
0.219512195121951	0.975609756097561\\
0.243471773190749	0.969907983088784\\
0.267748478701825	0.963488843813387\\
0.292323597828896	0.956319462393383\\
0.317175974710221	0.948366701791359\\
0.342281879194631	0.939597315436242\\
0.36761487964989	0.929978118161926\\
0.39314572304263	0.919476177208136\\
0.418842224744608	0.908059023836549\\
0.444669170759896	0.895694885871136\\
0.470588235294118	0.882352941176471\\
0.49655791679138	0.868003591739001\\
0.522533495736906	0.852618757612667\\
0.548467017652524	0.836172189532363\\
0.574307304785894	0.818639798488665\\
0.6	0.8\\
0.625487646293888	0.780234070221066\\
0.650709805216243	0.759326510399472\\
0.675603217158177	0.737265415549598\\
0.700102006120367	0.714042842570554\\
0.724137931034483	0.689655172413793\\
0.747640685075148	0.664103460328556\\
0.770538243626062	0.637393767705382\\
0.792757260666906	0.609537468626748\\
0.81422351233672	0.58055152394775\\
0.834862385321101	0.550458715596331\\
0.85459940652819	0.519287833827893\\
0.87336080929187	0.487073810415886\\
0.8910741301059	0.453857791225416\\
0.907668828691339	0.419687142312094\\
0.923076923076923	0.384615384615384\\
0.937233630375823	0.348702053467648\\
0.950078003120125	0.31201248049922\\
0.961553550411926	0.27461749705767\\
0.971608832807571	0.236593059936909\\
0.98019801980198	0.198019801980198\\
0.987281399046105	0.158982511923688\\
0.992825827022718	0.119569549621363\\
0.996805111821086	0.0798722044728435\\
0.999200319872051	0.039984006397441\\
1	0\\
0.999200319872051	-0.039984006397441\\
0.996805111821086	-0.0798722044728435\\
0.992825827022718	-0.119569549621363\\
0.987281399046105	-0.158982511923688\\
0.98019801980198	-0.198019801980198\\
0.971608832807571	-0.236593059936909\\
0.961553550411926	-0.274617497057669\\
0.950078003120125	-0.31201248049922\\
0.937233630375824	-0.348702053467648\\
0.923076923076923	-0.384615384615384\\
0.907668828691339	-0.419687142312094\\
0.8910741301059	-0.453857791225416\\
0.873360809291869	-0.487073810415886\\
0.85459940652819	-0.519287833827893\\
0.834862385321101	-0.550458715596331\\
0.81422351233672	-0.58055152394775\\
0.792757260666906	-0.609537468626748\\
0.770538243626062	-0.637393767705383\\
0.747640685075149	-0.664103460328556\\
0.724137931034483	-0.689655172413793\\
0.700102006120367	-0.714042842570554\\
0.675603217158177	-0.737265415549598\\
0.650709805216243	-0.759326510399472\\
0.625487646293888	-0.780234070221066\\
0.6	-0.8\\
0.574307304785894	-0.818639798488665\\
0.548467017652524	-0.836172189532363\\
0.522533495736906	-0.852618757612667\\
0.49655791679138	-0.868003591739\\
0.470588235294117	-0.882352941176471\\
0.444669170759896	-0.895694885871136\\
0.418842224744609	-0.908059023836549\\
0.39314572304263	-0.919476177208136\\
0.367614879649891	-0.929978118161926\\
0.342281879194631	-0.939597315436242\\
0.317175974710221	-0.948366701791359\\
0.292323597828896	-0.956319462393383\\
0.267748478701825	-0.963488843813387\\
0.243471773190749	-0.969907983088784\\
0.219512195121951	-0.975609756097561\\
0.195886151638364	-0.980626644343458\\
0.172607879924953	-0.984990619136961\\
0.149689583812371	-0.988733042078639\\
0.127141568981064	-0.991884580703336\\
0.104972375690608	-0.994475138121547\\
0.0831889081455807	-0.996533795493934\\
0.0617965597791463	-0.998088766192398\\
0.0407993338884264	-0.999167360532889\\
0.0201999591920017	-0.999795960008162\\
0	-1\\
0	-1\\
-0.0201999591920017	-0.999795960008162\\
-0.0407993338884263	-0.999167360532889\\
-0.0617965597791464	-0.998088766192398\\
-0.0831889081455808	-0.996533795493934\\
-0.104972375690608	-0.994475138121547\\
-0.127141568981064	-0.991884580703336\\
-0.149689583812371	-0.988733042078639\\
-0.172607879924953	-0.984990619136961\\
-0.195886151638364	-0.980626644343458\\
-0.219512195121951	-0.975609756097561\\
-0.243471773190749	-0.969907983088784\\
-0.267748478701826	-0.963488843813388\\
-0.292323597828896	-0.956319462393383\\
-0.317175974710222	-0.948366701791359\\
-0.342281879194631	-0.939597315436242\\
-0.36761487964989	-0.929978118161926\\
-0.39314572304263	-0.919476177208136\\
-0.418842224744608	-0.908059023836549\\
-0.444669170759896	-0.895694885871136\\
-0.470588235294117	-0.882352941176471\\
-0.49655791679138	-0.868003591739\\
-0.522533495736906	-0.852618757612667\\
-0.548467017652524	-0.836172189532363\\
-0.574307304785894	-0.818639798488665\\
-0.6	-0.8\\
-0.625487646293888	-0.780234070221066\\
-0.650709805216243	-0.759326510399472\\
-0.675603217158177	-0.737265415549598\\
-0.700102006120367	-0.714042842570554\\
-0.724137931034483	-0.689655172413793\\
-0.747640685075149	-0.664103460328556\\
-0.770538243626062	-0.637393767705382\\
-0.792757260666906	-0.609537468626748\\
-0.81422351233672	-0.58055152394775\\
-0.834862385321101	-0.55045871559633\\
-0.85459940652819	-0.519287833827894\\
-0.87336080929187	-0.487073810415886\\
-0.8910741301059	-0.453857791225416\\
-0.90766882869134	-0.419687142312094\\
-0.923076923076923	-0.384615384615385\\
-0.937233630375823	-0.348702053467648\\
-0.950078003120125	-0.31201248049922\\
-0.961553550411926	-0.27461749705767\\
-0.971608832807571	-0.236593059936909\\
-0.98019801980198	-0.198019801980198\\
-0.987281399046105	-0.158982511923688\\
-0.992825827022718	-0.119569549621363\\
-0.996805111821086	-0.0798722044728435\\
-0.999200319872051	-0.039984006397441\\
-1	0\\
-0.999200319872051	0.039984006397441\\
-0.996805111821086	0.0798722044728435\\
-0.992825827022718	0.119569549621363\\
-0.987281399046105	0.158982511923688\\
-0.98019801980198	0.198019801980198\\
-0.971608832807571	0.236593059936909\\
-0.961553550411926	0.274617497057669\\
-0.950078003120124	0.31201248049922\\
-0.937233630375824	0.348702053467648\\
-0.923076923076923	0.384615384615385\\
-0.907668828691339	0.419687142312094\\
-0.8910741301059	0.453857791225416\\
-0.87336080929187	0.487073810415886\\
-0.85459940652819	0.519287833827894\\
-0.834862385321101	0.55045871559633\\
-0.81422351233672	0.580551523947751\\
-0.792757260666906	0.609537468626748\\
-0.770538243626062	0.637393767705383\\
-0.747640685075149	0.664103460328556\\
-0.724137931034483	0.689655172413793\\
-0.700102006120367	0.714042842570554\\
-0.675603217158177	0.737265415549598\\
-0.650709805216243	0.759326510399472\\
-0.625487646293888	0.780234070221066\\
-0.6	0.8\\
-0.574307304785894	0.818639798488665\\
-0.548467017652524	0.836172189532363\\
-0.522533495736906	0.852618757612667\\
-0.49655791679138	0.868003591739\\
-0.470588235294118	0.882352941176471\\
-0.444669170759896	0.895694885871136\\
-0.418842224744609	0.908059023836549\\
-0.39314572304263	0.919476177208136\\
-0.367614879649891	0.929978118161925\\
-0.342281879194631	0.939597315436242\\
-0.317175974710222	0.948366701791359\\
-0.292323597828896	0.956319462393383\\
-0.267748478701826	0.963488843813387\\
-0.243471773190748	0.969907983088784\\
-0.219512195121951	0.975609756097561\\
-0.195886151638364	0.980626644343459\\
-0.172607879924953	0.984990619136961\\
-0.14968958381237	0.988733042078639\\
-0.127141568981064	0.991884580703336\\
-0.104972375690608	0.994475138121547\\
-0.0831889081455808	0.996533795493934\\
-0.0617965597791464	0.998088766192398\\
-0.0407993338884264	0.999167360532889\\
-0.0201999591920017	0.999795960008162\\
0	1\\
};

\addplot [color=black, line width=0.5pt]
  table[row sep=crcr]{%
2	1\\
2	-1\\
-2	-1\\
-2	1\\
2	1\\
};

\addplot [only marks,mark=*,myCPColor,mark options={fill=myCPColor}]
  table[row sep=crcr]{%
2	1\\
2	-1\\
-2	-1\\
-2	1\\
};

\end{axis}
\end{tikzpicture}%
%
}; &
      \node[] (c2p) {%
	\tikzsetnextfilename{./tikz/images/circle2_perturb}%
%
%
\begin{tikzpicture}

\begin{axis}[%
ticks=none,
width=0.16\textwidth,
height=0.12\textwidth,
at={(0\textwidth,0\textwidth)},
scale only axis,
xmin=-2,
xmax=2,
ymin=-1,
ymax=2,
axis background/.style={fill=white},
axis line style={draw=none},
tick style={draw=none}
]
\addplot [color=gray, line width=1.0pt]
  table[row sep=crcr]{%
0	1.5\\
0.0201999591920017	1.50474443991022\\
0.0407993338884264	1.50895503746878\\
0.0617965597791463	1.51259662348694\\
0.0831889081455806	1.51563258232236\\
0.104972375690608	1.51802486187845\\
0.127141568981064	1.51973399458972\\
0.149689583812371	1.52071913083467\\
0.172607879924953	1.52093808630394\\
0.195886151638364	1.52034740492705\\
0.219512195121951	1.51890243902439\\
0.243471773190749	1.51655744839592\\
0.267748478701825	1.51326572008114\\
0.292323597828896	1.50897971051951\\
0.317175974710221	1.5036512118019\\
0.342281879194631	1.49723154362416\\
0.36761487964989	1.48967177242888\\
0.39314572304263	1.48092295904152\\
0.418842224744608	1.47093643586833\\
0.444669170759896	1.4596641144178\\
0.470588235294118	1.44705882352941\\
0.49655791679138	1.43307467824005\\
0.522533495736906	1.41766747868453\\
0.548467017652524	1.40079513781356\\
0.574307304785894	1.38241813602015\\
0.6	1.3625\\
0.625487646293888	1.3410078023407\\
0.650709805216243	1.3179126774513\\
0.675603217158177	1.29319034852547\\
0.700102006120367	1.26682165929956\\
0.724137931034483	1.23879310344828\\
0.747640685075148	1.20909734358616\\
0.770538243626062	1.17773371104816\\
0.792757260666906	1.14470867694514\\
0.81422351233672	1.11003628447025\\
0.834862385321101	1.07373853211009\\
0.85459940652819	1.03584569732938\\
0.87336080929187	0.996396590483327\\
0.8910741301059	0.955438729198185\\
0.907668828691339	0.913028424265547\\
0.923076923076923	0.869230769230769\\
0.937233630375823	0.824119527314994\\
0.950078003120125	0.777776911076443\\
0.961553550411926	0.730293252255787\\
0.971608832807571	0.681766561514196\\
0.98019801980198	0.63230198019802\\
0.987281399046105	0.582011128775835\\
0.992825827022718	0.531011359107214\\
0.996805111821086	0.479424920127796\\
0.999200319872051	0.427378048780488\\
1	0.375\\
0.999200319872051	0.322422031187525\\
0.996805111821086	0.269776357827476\\
0.992825827022718	0.217195097648465\\
0.987281399046105	0.164809220985691\\
0.98019801980198	0.112747524752475\\
0.971608832807571	0.0611356466876969\\
0.961553550411926	0.0100951353471953\\
0.950078003120125	-0.0402574102964118\\
0.937233630375824	-0.0898111197210383\\
0.923076923076923	-0.138461538461538\\
0.907668828691339	-0.186111217092713\\
0.8910741301059	-0.23267019667171\\
0.873360809291869	-0.27805638816036\\
0.85459940652819	-0.322195845697329\\
0.834862385321101	-0.365022935779817\\
0.81422351233672	-0.406480406386067\\
0.792757260666906	-0.446519361778415\\
0.770538243626062	-0.485099150141643\\
0.747640685075149	-0.522187172317371\\
0.724137931034483	-0.557758620689655\\
0.700102006120367	-0.591796157769466\\
0.675603217158177	-0.624289544235925\\
0.650709805216243	-0.655235226147243\\
0.625487646293888	-0.68463589076723\\
0.6	-0.7125\\
0.574307304785894	-0.738841309823678\\
0.548467017652524	-0.763678383400434\\
0.522533495736906	-0.787034104750305\\
0.49655791679138	-0.808935199042203\\
0.470588235294117	-0.829411764705883\\
0.444669170759896	-0.848496821727824\\
0.418842224744609	-0.866225879682179\\
0.39314572304263	-0.882636528280858\\
0.367614879649891	-0.897768052516411\\
0.342281879194631	-0.911661073825503\\
0.317175974710221	-0.924357218124341\\
0.292323597828896	-0.93589881106229\\
0.267748478701825	-0.946328600405679\\
0.243471773190749	-0.955689505098234\\
0.219512195121951	-0.964024390243903\\
0.195886151638364	-0.97137586701746\\
0.172607879924953	-0.977786116322702\\
0.149689583812371	-0.983296734881582\\
0.127141568981064	-0.987948602344455\\
0.104972375690608	-0.991781767955801\\
0.0831889081455807	-0.994835355285962\\
0.0617965597791463	-0.997147483542153\\
0.0407993338884264	-0.998755203996669\\
0.0201999591920017	-0.999694450112222\\
0	-1\\
0	-1\\
-0.0201999591920017	-0.999795449908182\\
-0.0407993338884263	-0.999163197335554\\
-0.0617965597791464	-0.998074431938841\\
-0.0831889081455808	-0.996499133448874\\
-0.104972375690608	-0.994406077348066\\
-0.127141568981064	-0.991762849413886\\
-0.149689583812371	-0.988535870315015\\
-0.172607879924953	-0.9846904315197\\
-0.195886151638364	-0.980190743841186\\
-0.219512195121951	-0.975\\
-0.243471773190749	-0.969080452623726\\
-0.267748478701826	-0.962393509127789\\
-0.292323597828896	-0.954899844921168\\
-0.317175974710222	-0.946559536354057\\
-0.342281879194631	-0.937332214765101\\
-0.36761487964989	-0.927177242888403\\
-0.39314572304263	-0.916053914739482\\
-0.418842224744608	-0.903921679909194\\
-0.444669170759896	-0.890740392950015\\
-0.470588235294117	-0.876470588235294\\
-0.49655791679138	-0.861073780305298\\
-0.522533495736906	-0.844512789281364\\
-0.548467017652524	-0.826752090430474\\
-0.574307304785894	-0.807758186397985\\
-0.6	-0.7875\\
-0.625487646293888	-0.765949284785436\\
-0.650709805216243	-0.743081049851436\\
-0.675603217158177	-0.71887399463807\\
-0.700102006120367	-0.693310948656919\\
-0.724137931034483	-0.666379310344827\\
-0.747640685075149	-0.638071478504019\\
-0.770538243626062	-0.608385269121813\\
-0.792757260666906	-0.577324309788455\\
-0.81422351233672	-0.544898403483309\\
-0.834862385321101	-0.511123853211009\\
-0.85459940652819	-0.476023738872404\\
-0.87336080929187	-0.439628137879355\\
-0.8910741301059	-0.401974281391831\\
-0.90766882869134	-0.363106638687524\\
-0.923076923076923	-0.323076923076923\\
-0.937233630375823	-0.281944013948082\\
-0.950078003120125	-0.239773790951638\\
-0.961553550411926	-0.196638877991369\\
-0.971608832807571	-0.152618296529969\\
-0.98019801980198	-0.107797029702971\\
-0.987281399046105	-0.0622655007949127\\
-0.992825827022718	-0.0161189717018734\\
-0.996805111821086	0.0305431309904153\\
-0.999200319872051	0.0776179528188725\\
-1	0.125\\
-0.999200319872051	0.172581967213115\\
-0.996805111821086	0.220255591054313\\
-0.992825827022718	0.267912514946194\\
-0.987281399046105	0.315445151033386\\
-0.98019801980198	0.362747524752476\\
-0.971608832807571	0.409716088328076\\
-0.961553550411926	0.456250490388387\\
-0.950078003120124	0.502254290171607\\
-0.937233630375824	0.547635606354126\\
-0.923076923076923	0.592307692307693\\
-0.907668828691339	0.636189431514689\\
-0.8910741301059	0.679205748865356\\
-0.87336080929187	0.721287935556388\\
-0.85459940652819	0.762373887240357\\
-0.834862385321101	0.802408256880734\\
-0.81422351233672	0.84134252539913\\
-0.792757260666906	0.879134994621728\\
-0.770538243626062	0.915750708215298\\
-0.747640685075149	0.951161307235232\\
-0.724137931034483	0.985344827586207\\
-0.700102006120367	1.01828544712683\\
-0.675603217158177	1.04997319034853\\
-0.650709805216243	1.08040359854738\\
-0.625487646293888	1.10957737321196\\
-0.6	1.1375\\
-0.574307304785894	1.16418136020151\\
-0.548467017652524	1.18963533601734\\
-0.522533495736906	1.21387941534714\\
-0.49655791679138	1.23693430110745\\
-0.470588235294118	1.25882352941176\\
-0.444669170759896	1.27957310026004\\
-0.418842224744609	1.29921112372304\\
-0.39314572304263	1.31776748397882\\
-0.367614879649891	1.33527352297593\\
-0.342281879194631	1.35176174496644\\
-0.317175974710222	1.3672655426765\\
-0.292323597828896	1.38181894546394\\
-0.267748478701826	1.39545638945233\\
-0.243471773190748	1.40821250932604\\
-0.219512195121951	1.42012195121951\\
-0.195886151638364	1.4312192059316\\
-0.172607879924953	1.44153846153846\\
-0.14968958381237	1.45111347436192\\
-0.127141568981064	1.45997745716862\\
-0.104972375690608	1.46816298342541\\
-0.0831889081455808	1.47570190641248\\
-0.0617965597791464	1.48262529199405\\
-0.0407993338884264	1.48896336386345\\
-0.0201999591920017	1.49474546011018\\
0	1.5\\
};

\addplot [color=black, line width=0.5pt]
  table[row sep=crcr]{%
2	2\\
2	-1\\
-2	-1\\
-2	1\\
2	2\\
};

\addplot [only marks,mark=*,myCPColor,mark options={fill=myCPColor}]
  table[row sep=crcr]{%
2	2\\
2	-1\\
-2	-1\\
-2	1\\
};

\end{axis}
\end{tikzpicture}%
%
};\\
      \node[] (e2) {%
	\tikzsetnextfilename{./tikz/images/ellipse2}%
%
%
\begin{tikzpicture}

\begin{axis}[%
ticks=none,
width=0.16\textwidth,
height=0.04\textwidth,
at={(0\textwidth,0\textwidth)},
scale only axis,
xmin=-2,
xmax=2,
ymin=-0.5,
ymax=0.5,
axis background/.style={fill=white},
axis line style={draw=none},
tick style={draw=none}
]
\addplot [color=gray, line width=1.0pt]
  table[row sep=crcr]{%
0	0.5\\
0.0201999591920017	0.499897980004081\\
0.0407993338884264	0.499583680266445\\
0.0617965597791463	0.499044383096199\\
0.0831889081455806	0.498266897746967\\
0.104972375690608	0.497237569060773\\
0.127141568981064	0.495942290351668\\
0.149689583812371	0.494366521039319\\
0.172607879924953	0.49249530956848\\
0.195886151638364	0.490313322171729\\
0.219512195121951	0.487804878048781\\
0.243471773190749	0.484953991544392\\
0.267748478701825	0.481744421906694\\
0.292323597828896	0.478159731196692\\
0.317175974710221	0.47418335089568\\
0.342281879194631	0.469798657718121\\
0.36761487964989	0.464989059080963\\
0.39314572304263	0.459738088604068\\
0.418842224744608	0.454029511918275\\
0.444669170759896	0.447847442935568\\
0.470588235294118	0.441176470588235\\
0.49655791679138	0.4340017958695\\
0.522533495736906	0.426309378806334\\
0.548467017652524	0.418086094766182\\
0.574307304785894	0.409319899244333\\
0.6	0.4\\
0.625487646293888	0.390117035110533\\
0.650709805216243	0.379663255199736\\
0.675603217158177	0.368632707774799\\
0.700102006120367	0.357021421285277\\
0.724137931034483	0.344827586206896\\
0.747640685075148	0.332051730164278\\
0.770538243626062	0.318696883852691\\
0.792757260666906	0.304768734313374\\
0.81422351233672	0.290275761973875\\
0.834862385321101	0.275229357798165\\
0.85459940652819	0.259643916913947\\
0.87336080929187	0.243536905207943\\
0.8910741301059	0.226928895612708\\
0.907668828691339	0.209843571156047\\
0.923076923076923	0.192307692307692\\
0.937233630375823	0.174351026733824\\
0.950078003120125	0.15600624024961\\
0.961553550411926	0.137308748528835\\
0.971608832807571	0.118296529968454\\
0.98019801980198	0.099009900990099\\
0.987281399046105	0.0794912559618442\\
0.992825827022718	0.0597847748106816\\
0.996805111821086	0.0399361022364217\\
0.999200319872051	0.0199920031987205\\
1	0\\
0.999200319872051	-0.0199920031987205\\
0.996805111821086	-0.0399361022364217\\
0.992825827022718	-0.0597847748106816\\
0.987281399046105	-0.0794912559618442\\
0.98019801980198	-0.0990099009900991\\
0.971608832807571	-0.118296529968454\\
0.961553550411926	-0.137308748528835\\
0.950078003120125	-0.15600624024961\\
0.937233630375824	-0.174351026733824\\
0.923076923076923	-0.192307692307692\\
0.907668828691339	-0.209843571156047\\
0.8910741301059	-0.226928895612708\\
0.873360809291869	-0.243536905207943\\
0.85459940652819	-0.259643916913947\\
0.834862385321101	-0.275229357798165\\
0.81422351233672	-0.290275761973875\\
0.792757260666906	-0.304768734313374\\
0.770538243626062	-0.318696883852691\\
0.747640685075149	-0.332051730164278\\
0.724137931034483	-0.344827586206897\\
0.700102006120367	-0.357021421285277\\
0.675603217158177	-0.368632707774799\\
0.650709805216243	-0.379663255199736\\
0.625487646293888	-0.390117035110533\\
0.6	-0.4\\
0.574307304785894	-0.409319899244333\\
0.548467017652524	-0.418086094766182\\
0.522533495736906	-0.426309378806334\\
0.49655791679138	-0.4340017958695\\
0.470588235294117	-0.441176470588235\\
0.444669170759896	-0.447847442935568\\
0.418842224744609	-0.454029511918275\\
0.39314572304263	-0.459738088604068\\
0.367614879649891	-0.464989059080963\\
0.342281879194631	-0.469798657718121\\
0.317175974710221	-0.47418335089568\\
0.292323597828896	-0.478159731196692\\
0.267748478701825	-0.481744421906694\\
0.243471773190749	-0.484953991544392\\
0.219512195121951	-0.487804878048781\\
0.195886151638364	-0.490313322171729\\
0.172607879924953	-0.49249530956848\\
0.149689583812371	-0.494366521039319\\
0.127141568981064	-0.495942290351668\\
0.104972375690608	-0.497237569060773\\
0.0831889081455807	-0.498266897746967\\
0.0617965597791463	-0.499044383096199\\
0.0407993338884264	-0.499583680266445\\
0.0201999591920017	-0.499897980004081\\
0	-0.5\\
0	-0.5\\
-0.0201999591920017	-0.499897980004081\\
-0.0407993338884263	-0.499583680266445\\
-0.0617965597791464	-0.499044383096199\\
-0.0831889081455808	-0.498266897746967\\
-0.104972375690608	-0.497237569060773\\
-0.127141568981064	-0.495942290351668\\
-0.149689583812371	-0.494366521039319\\
-0.172607879924953	-0.49249530956848\\
-0.195886151638364	-0.490313322171729\\
-0.219512195121951	-0.48780487804878\\
-0.243471773190749	-0.484953991544392\\
-0.267748478701826	-0.481744421906694\\
-0.292323597828896	-0.478159731196692\\
-0.317175974710222	-0.47418335089568\\
-0.342281879194631	-0.469798657718121\\
-0.36761487964989	-0.464989059080963\\
-0.39314572304263	-0.459738088604068\\
-0.418842224744608	-0.454029511918275\\
-0.444669170759896	-0.447847442935568\\
-0.470588235294117	-0.441176470588235\\
-0.49655791679138	-0.4340017958695\\
-0.522533495736906	-0.426309378806334\\
-0.548467017652524	-0.418086094766181\\
-0.574307304785894	-0.409319899244333\\
-0.6	-0.4\\
-0.625487646293888	-0.390117035110533\\
-0.650709805216243	-0.379663255199736\\
-0.675603217158177	-0.368632707774799\\
-0.700102006120367	-0.357021421285277\\
-0.724137931034483	-0.344827586206896\\
-0.747640685075149	-0.332051730164278\\
-0.770538243626062	-0.318696883852691\\
-0.792757260666906	-0.304768734313374\\
-0.81422351233672	-0.290275761973875\\
-0.834862385321101	-0.275229357798165\\
-0.85459940652819	-0.259643916913947\\
-0.87336080929187	-0.243536905207943\\
-0.8910741301059	-0.226928895612708\\
-0.90766882869134	-0.209843571156047\\
-0.923076923076923	-0.192307692307692\\
-0.937233630375823	-0.174351026733824\\
-0.950078003120125	-0.15600624024961\\
-0.961553550411926	-0.137308748528835\\
-0.971608832807571	-0.118296529968454\\
-0.98019801980198	-0.0990099009900991\\
-0.987281399046105	-0.0794912559618442\\
-0.992825827022718	-0.0597847748106816\\
-0.996805111821086	-0.0399361022364217\\
-0.999200319872051	-0.0199920031987205\\
-1	0\\
-0.999200319872051	0.0199920031987205\\
-0.996805111821086	0.0399361022364217\\
-0.992825827022718	0.0597847748106816\\
-0.987281399046105	0.0794912559618442\\
-0.98019801980198	0.0990099009900991\\
-0.971608832807571	0.118296529968454\\
-0.961553550411926	0.137308748528835\\
-0.950078003120124	0.15600624024961\\
-0.937233630375824	0.174351026733824\\
-0.923076923076923	0.192307692307692\\
-0.907668828691339	0.209843571156047\\
-0.8910741301059	0.226928895612708\\
-0.87336080929187	0.243536905207943\\
-0.85459940652819	0.259643916913947\\
-0.834862385321101	0.275229357798165\\
-0.81422351233672	0.290275761973875\\
-0.792757260666906	0.304768734313374\\
-0.770538243626062	0.318696883852691\\
-0.747640685075149	0.332051730164278\\
-0.724137931034483	0.344827586206897\\
-0.700102006120367	0.357021421285277\\
-0.675603217158177	0.368632707774799\\
-0.650709805216243	0.379663255199736\\
-0.625487646293888	0.390117035110533\\
-0.6	0.4\\
-0.574307304785894	0.409319899244333\\
-0.548467017652524	0.418086094766182\\
-0.522533495736906	0.426309378806334\\
-0.49655791679138	0.4340017958695\\
-0.470588235294118	0.441176470588235\\
-0.444669170759896	0.447847442935568\\
-0.418842224744609	0.454029511918275\\
-0.39314572304263	0.459738088604068\\
-0.367614879649891	0.464989059080963\\
-0.342281879194631	0.469798657718121\\
-0.317175974710222	0.47418335089568\\
-0.292323597828896	0.478159731196692\\
-0.267748478701826	0.481744421906694\\
-0.243471773190748	0.484953991544392\\
-0.219512195121951	0.48780487804878\\
-0.195886151638364	0.490313322171729\\
-0.172607879924953	0.49249530956848\\
-0.14968958381237	0.494366521039319\\
-0.127141568981064	0.495942290351668\\
-0.104972375690608	0.497237569060773\\
-0.0831889081455808	0.498266897746967\\
-0.0617965597791464	0.499044383096199\\
-0.0407993338884264	0.499583680266445\\
-0.0201999591920017	0.499897980004081\\
0	0.5\\
};

\addplot [color=black, line width=0.5pt]
  table[row sep=crcr]{%
2	0.5\\
2	-0.5\\
-2	-0.5\\
-2	0.5\\
2	0.5\\
};

\addplot [only marks,mark=*,myCPColor,mark options={fill=myCPColor}]
  table[row sep=crcr]{%
2	0.5\\
2	-0.5\\
-2	-0.5\\
-2	0.5\\
};

\end{axis}
\end{tikzpicture}%
%
}; &
      \node[] (e2p) {%
	\tikzsetnextfilename{./tikz/images/ellipse2_perturb}%
%
%
\begin{tikzpicture}

\begin{axis}[%
ticks=none,
width=0.16\textwidth,
height=0.06\textwidth,
at={(0\textwidth,0\textwidth)},
scale only axis,
xmin=-2,
xmax=2,
ymin=-0.5,
ymax=1,
axis background/.style={fill=white},
axis line style={draw=none},
tick style={draw=none}
]
\addplot [color=gray, line width=1.0pt]
  table[row sep=crcr]{%
0	0.75\\
0.0201999591920017	0.752372219955111\\
0.0407993338884264	0.754477518734388\\
0.0617965597791463	0.75629831174347\\
0.0831889081455806	0.757816291161179\\
0.104972375690608	0.759012430939226\\
0.127141568981064	0.75986699729486\\
0.149689583812371	0.760359565417337\\
0.172607879924953	0.76046904315197\\
0.195886151638364	0.760173702463525\\
0.219512195121951	0.759451219512195\\
0.243471773190749	0.758278724197961\\
0.267748478701825	0.756632860040568\\
0.292323597828896	0.754489855259757\\
0.317175974710221	0.751825605900948\\
0.342281879194631	0.748615771812081\\
0.36761487964989	0.744835886214442\\
0.39314572304263	0.740461479520758\\
0.418842224744608	0.735468217934166\\
0.444669170759896	0.729832057208899\\
0.470588235294118	0.723529411764706\\
0.49655791679138	0.716537339120024\\
0.522533495736906	0.708833739342265\\
0.548467017652524	0.700397568906782\\
0.574307304785894	0.691209068010076\\
0.6	0.68125\\
0.625487646293888	0.670503901170351\\
0.650709805216243	0.658956338725652\\
0.675603217158177	0.646595174262735\\
0.700102006120367	0.633410829649779\\
0.724137931034483	0.619396551724138\\
0.747640685075148	0.604548671793079\\
0.770538243626062	0.588866855524079\\
0.792757260666906	0.572354338472571\\
0.81422351233672	0.555018142235123\\
0.834862385321101	0.536869266055046\\
0.85459940652819	0.517922848664689\\
0.87336080929187	0.498198295241664\\
0.8910741301059	0.477719364599092\\
0.907668828691339	0.456514212132774\\
0.923076923076923	0.434615384615385\\
0.937233630375823	0.412059763657497\\
0.950078003120125	0.388888455538222\\
0.961553550411926	0.365146626127893\\
0.971608832807571	0.340883280757098\\
0.98019801980198	0.31615099009901\\
0.987281399046105	0.291005564387917\\
0.992825827022718	0.265505679553607\\
0.996805111821086	0.239712460063898\\
0.999200319872051	0.213689024390244\\
1	0.1875\\
0.999200319872051	0.161211015593762\\
0.996805111821086	0.134888178913738\\
0.992825827022718	0.108597548824233\\
0.987281399046105	0.0824046104928457\\
0.98019801980198	0.0563737623762375\\
0.971608832807571	0.0305678233438485\\
0.961553550411926	0.00504756767359765\\
0.950078003120125	-0.0201287051482059\\
0.937233630375824	-0.0449055598605191\\
0.923076923076923	-0.0692307692307692\\
0.907668828691339	-0.0930556085463563\\
0.8910741301059	-0.116335098335855\\
0.873360809291869	-0.13902819408018\\
0.85459940652819	-0.161097922848665\\
0.834862385321101	-0.182511467889908\\
0.81422351233672	-0.203240203193033\\
0.792757260666906	-0.223259680889208\\
0.770538243626062	-0.242549575070822\\
0.747640685075149	-0.261093586158686\\
0.724137931034483	-0.278879310344828\\
0.700102006120367	-0.295898078884733\\
0.675603217158177	-0.312144772117962\\
0.650709805216243	-0.327617613073622\\
0.625487646293888	-0.342317945383615\\
0.6	-0.35625\\
0.574307304785894	-0.369420654911839\\
0.548467017652524	-0.381839191700217\\
0.522533495736906	-0.393517052375152\\
0.49655791679138	-0.404467599521101\\
0.470588235294117	-0.414705882352941\\
0.444669170759896	-0.424248410863912\\
0.418842224744609	-0.43311293984109\\
0.39314572304263	-0.441318264140429\\
0.367614879649891	-0.448884026258206\\
0.342281879194631	-0.455830536912752\\
0.317175974710221	-0.462178609062171\\
0.292323597828896	-0.467949405531145\\
0.267748478701825	-0.47316430020284\\
0.243471773190749	-0.477844752549117\\
0.219512195121951	-0.482012195121951\\
0.195886151638364	-0.48568793350873\\
0.172607879924953	-0.488893058161351\\
0.149689583812371	-0.491648367440791\\
0.127141568981064	-0.493974301172227\\
0.104972375690608	-0.495890883977901\\
0.0831889081455807	-0.497417677642981\\
0.0617965597791463	-0.498573741771077\\
0.0407993338884264	-0.499377601998335\\
0.0201999591920017	-0.499847225056111\\
0	-0.5\\
0	-0.5\\
-0.0201999591920017	-0.499897724954091\\
-0.0407993338884263	-0.499581598667777\\
-0.0617965597791464	-0.49903721596942\\
-0.0831889081455808	-0.498249566724437\\
-0.104972375690608	-0.497203038674033\\
-0.127141568981064	-0.495881424706943\\
-0.149689583812371	-0.494267935157507\\
-0.172607879924953	-0.49234521575985\\
-0.195886151638364	-0.490095371920593\\
-0.219512195121951	-0.4875\\
-0.243471773190749	-0.484540226311863\\
-0.267748478701826	-0.481196754563895\\
-0.292323597828896	-0.477449922460584\\
-0.317175974710222	-0.473279768177028\\
-0.342281879194631	-0.46866610738255\\
-0.36761487964989	-0.463588621444201\\
-0.39314572304263	-0.458026957369741\\
-0.418842224744608	-0.451960839954597\\
-0.444669170759896	-0.445370196475007\\
-0.470588235294117	-0.438235294117647\\
-0.49655791679138	-0.430536890152649\\
-0.522533495736906	-0.422256394640682\\
-0.548467017652524	-0.413376045215237\\
-0.574307304785894	-0.403879093198993\\
-0.6	-0.39375\\
-0.625487646293888	-0.382974642392718\\
-0.650709805216243	-0.371540524925718\\
-0.675603217158177	-0.359436997319035\\
-0.700102006120367	-0.34665547432846\\
-0.724137931034483	-0.333189655172414\\
-0.747640685075149	-0.31903573925201\\
-0.770538243626062	-0.304192634560906\\
-0.792757260666906	-0.288662154894227\\
-0.81422351233672	-0.272449201741654\\
-0.834862385321101	-0.255561926605504\\
-0.85459940652819	-0.238011869436202\\
-0.87336080929187	-0.219814068939678\\
-0.8910741301059	-0.200987140695916\\
-0.90766882869134	-0.181553319343762\\
-0.923076923076923	-0.161538461538462\\
-0.937233630375823	-0.140972006974041\\
-0.950078003120125	-0.119886895475819\\
-0.961553550411926	-0.0983194389956847\\
-0.971608832807571	-0.0763091482649843\\
-0.98019801980198	-0.0538985148514853\\
-0.987281399046105	-0.0311327503974563\\
-0.992825827022718	-0.00805948585093669\\
-0.996805111821086	0.0152715654952076\\
-0.999200319872051	0.0388089764094362\\
-1	0.0625\\
-0.999200319872051	0.0862909836065573\\
-0.996805111821086	0.110127795527157\\
-0.992825827022718	0.133956257473097\\
-0.987281399046105	0.157722575516693\\
-0.98019801980198	0.181373762376238\\
-0.971608832807571	0.204858044164038\\
-0.961553550411926	0.228125245194193\\
-0.950078003120124	0.251127145085804\\
-0.937233630375824	0.273817803177063\\
-0.923076923076923	0.296153846153846\\
-0.907668828691339	0.318094715757344\\
-0.8910741301059	0.339602874432678\\
-0.87336080929187	0.360643967778194\\
-0.85459940652819	0.381186943620178\\
-0.834862385321101	0.401204128440367\\
-0.81422351233672	0.420671262699565\\
-0.792757260666906	0.439567497310864\\
-0.770538243626062	0.457875354107649\\
-0.747640685075149	0.475580653617616\\
-0.724137931034483	0.492672413793103\\
-0.700102006120367	0.509142723563414\\
-0.675603217158177	0.524986595174263\\
-0.650709805216243	0.540201799273688\\
-0.625487646293888	0.554788686605982\\
-0.6	0.56875\\
-0.574307304785894	0.582090680100756\\
-0.548467017652524	0.594817668008671\\
-0.522533495736906	0.606939707673569\\
-0.49655791679138	0.618467150553726\\
-0.470588235294118	0.629411764705882\\
-0.444669170759896	0.63978655013002\\
-0.418842224744609	0.649605561861521\\
-0.39314572304263	0.658883741989412\\
-0.367614879649891	0.667636761487965\\
-0.342281879194631	0.675880872483221\\
-0.317175974710222	0.683632771338251\\
-0.292323597828896	0.690909472731972\\
-0.267748478701826	0.697728194726166\\
-0.243471773190748	0.704106254663019\\
-0.219512195121951	0.710060975609756\\
-0.195886151638364	0.715609602965798\\
-0.172607879924953	0.720769230769231\\
-0.14968958381237	0.725556737180961\\
-0.127141568981064	0.72998872858431\\
-0.104972375690608	0.734081491712707\\
-0.0831889081455808	0.737850953206239\\
-0.0617965597791464	0.741312645997027\\
-0.0407993338884264	0.744481681931723\\
-0.0201999591920017	0.747372730055091\\
0	0.75\\
};

\addplot [color=black, line width=0.5pt]
  table[row sep=crcr]{%
2	1\\
2	-0.5\\
-2	-0.5\\
-2	0.5\\
2	1\\
};

\addplot [only marks,mark=*,myCPColor,mark options={fill=myCPColor}]
  table[row sep=crcr]{%
2	1\\
2	-0.5\\
-2	-0.5\\
-2	0.5\\
};

\end{axis}
\end{tikzpicture}%
%
};\\
   };
   \coordinate (a) at (ce2.south -| b3.west);
   \matrix[matrix of nodes,row sep=0.04\textwidth,column sep=0.5cm,anchor=south west] (ce3) at (a) {
      \node[] (c3) {%
	\tikzsetnextfilename{./tikz/images/circle3}%
	\input{./tikz/scripts/circle3}%
}; &
      \node[] (c3p) {%
	\tikzsetnextfilename{./tikz/images/circle3_perturb}%
	\input{./tikz/scripts/circle3_perturb}%
};\\
      \node[] (e3) {%
	\tikzsetnextfilename{./tikz/images/ellipse3}%
	\input{./tikz/scripts/ellipse3}%
}; &
      \node[] (e3p) {%
	\tikzsetnextfilename{./tikz/images/ellipse3_perturb}%
	\input{./tikz/scripts/ellipse3_perturb}%
};\\
    };
   \coordinate (a) at (b1.east |- B.north);
   \node[fit=(a)(e1)(e1p), inner sep=1pt,draw,color=black!14,rounded corners] (c1) {};
   \coordinate (a) at (b2.east |- B.north);
   \node[fit=(a)(e2)(e2p), inner sep=1pt,draw,color=black!14,rounded corners] (c2) {};
   \coordinate (a) at (b3.east |- B.north);
   \node[fit=(a)(e3)(e3p), inner sep=1pt,draw,color=black!14,rounded corners] (c3) {};
   \node[below] at (c1.south) {\small(a) degree $2$};
   \node[below] at (c2.south) {\small(b) degree $3$};
   \node[below] at (c3.south) {\small(c) multi-degree $(3,2,2)$};
   \end{tikzpicture}
   \caption{Different $C^1$ smooth descriptions of unit circles and ellipses with axis lengths $(1, \frac{1}{2})$.
      The conics in figure boxes (a)--(c) are made up of $C^1$ splines of degree~$2$, degree~$3$ and multi-degree $(3,2,2)$, respectively.
      In each box, the middle and bottom rows of figures show the $C^1$ circles and ellipses built using \HS{the $C^1$ B-spline} functions shown at the top; the $t$-axis markers correspond to the breakpoint locations.
      Furthermore, in each of the bottom two rows of figures, the curve on the left shows the exact conic, and the curve on the right is obtained by raising one of the control points of the exact conic; see Examples~\ref{ex:ellipse_2},~\ref{ex:ellipse_3} and~\ref{ex:ellipse_322} for details.
      (Even though all circles and ellipses have the same respective dimensions, they are scaled differently only to accommodate their control nets into the figure.)
   }
   \label{fig:conic2}
\end{sidewaysfigure*}

\subsubsection{$C^1$ description of degree $2$}\label{ssec:ellipse2}
Here we present a $C^1$ quadratic description of the ellipse in Equation~\eqref{eq:ellipse} using $4$ rational pieces.
Consider the domain $\domain = [0, 4]$ in the periodic setting.
Choose $\nseg = 4$ and 
\begin{equation*}
\seg{i}{\bknot} = [0,0,0,1,1,1], \quad
\seg{i}{\brwts} = \Bigl[1,\tfrac{\sqrt{2}}{2},1\Bigl], \quad i = 1, \dots, 4.
\end{equation*}
We can define four $C^1$ quadratic piecewise-NURBS functions $\bsp_i$  on $\domain$ using the extraction matrix
\begin{equation}\label{eq:matH-ellipse_2}
\periodic{\extMat} = \begin{bmatrix*}[c]
\frac{1}{2} & 1 & \frac{1}{2} & \frac{1}{2} & 0 & 0 & 0 & 0 & 0 & 0 & 0 & \frac{1}{2}\\
0 & 0 & \frac{1}{2} & \frac{1}{2} & 1 & \frac{1}{2} & \frac{1}{2} & 0 & 0 & 0 & 0 & 0\\
0 & 0 & 0 & 0 & 0 & \frac{1}{2} & \frac{1}{2} & 1 & \frac{1}{2} & \frac{1}{2} & 0 & 0\\
\frac{1}{2} & 0 & 0 & 0 & 0 & 0 & 0 & 0 & \frac{1}{2} & \frac{1}{2} & 1 & \frac{1}{2}
\end{bmatrix*};
\end{equation}
see Equation \eqref{eq:mdbC1Extraction} taking into account Remark~\ref{rmk:simple-betas}. These spline basis functions are shown in Figure~\ref{fig:conic2} (a, top row).
Finally, we can build an ellipse centred at $(0,0)$ and with axis lengths $(a_x, a_y)$ by combining the splines $\bsp_{i}$ with the control points $\mbf{f}_i$ defined as
\begin{equation*}
\mbf{f}_1 = (a_x, a_y) = -\mbf{f}_3,\quad
\mbf{f}_2 = (a_x, -a_y) = -\mbf{f}_4.
\end{equation*}

\begin{remark}\label{rmk:verification}
To verify that the curve $\mbf{f}$ satisfies Equation \eqref{eq:ellipse}, we can proceed as follows.
The simplest approach is to numerically evaluate $\mbf{f}(t)$ at all $t$ and plug the result in that equation.
Alternatively, this verification can also be performed analytically by looking at the explicit \HS{expressions} of the rational pieces that form $\mbf{f}$.
For instance, consider the first quadratic rational piece, $\seg{1}{\mbf{g}}$, that is a part of $\mbf{f}$.
As discussed in Equation \eqref{eq:mdbControlPointExtraction}, we can get the control points of this piece, denoted with $\seg{1}{\mbf{g}_{j}}$, \HS{$j\in\{1,2,3\}$,} by applying the transpose of (a submatrix of) $\periodic{\extMat}$ from Equation \eqref{eq:matH-ellipse_2} to a vector containing the points $\mbf{f}_i$, i.e., 
\begin{equation*}
\seg{1}{\mbf{g}_{j}} = \sum_{i = 1}^4 \periodic{\extMatel}_{ij} \mbf{f}_i.
\end{equation*}
This yields the control points
\begin{equation*}
\seg{1}{\mbf{g}_1} = (0, a_y),\quad
\seg{1}{\mbf{g}_2} = (a_x, a_y),\quad
\seg{1}{\mbf{g}_3} = (a_x, 0).
\end{equation*}
Combining the above control points with the NURBS basis defined on the first segment,
\begin{equation*}
\ratbez_{1,2}(\coord) = \frac{(1-\coord)^2}{w(\coord)},\ 
\ratbez_{2,2}(\coord) = \frac{\sqrt{2}\coord(1-\coord)}{w(\coord)},\ 
\ratbez_{3,2}(\coord) = \frac{\coord^2}{w(\coord)},
\end{equation*}
where $w(\coord) := (1-\coord)^2 + \sqrt{2}\coord(1-\coord) + \coord^2$, some simple algebra shows that $\seg{1}{\mbf{g}}$ indeed satisfies Equation \eqref{eq:ellipse}.
Verifications for the other pieces of $\mbf{f}$ can be similarly done.
\end{remark}

\begin{example}\label{ex:ellipse_2}
Choosing $a_x = a_y = 1$, we obtain a circle of radius $1$, as shown in Figure~\ref{fig:conic2} (a, middle row).
This $C^1$ quadratic description is equivalent to the one used in \cite{lu2009circular}.
The choice $a_x = 2a_y = 1$ yields an ellipse with axis lengths $(1, \frac{1}{2})$, as shown in Figure~\ref{fig:conic2} (a, bottom row).
The perturbed versions of these conics, with the control points chosen as in Equation \eqref{eq:conic2-perturb}, are shown as well and they remain clearly smooth.
\end{example}

\subsubsection{$C^1$ description of degree $3$}\label{ssec:ellipse3}
Here we present a $C^1$ cubic representation of the ellipse in Equation~\eqref{eq:ellipse} that uses only 2 rational pieces.
Consider the domain $\domain = [0, 2]$ in the periodic setting.
Choose $\nseg = 2$ and
\begin{equation*}
\seg{i}{\bknot} = [0,0,0,0,1,1,1,1], \quad
\seg{i}{\brwts} = \Bigl[1, \tfrac{1}{3},\tfrac{1}{3},1\Bigr], \quad i = 1, 2.
\end{equation*}
We can define four $C^1$ cubic piecewise-NURBS functions $\bsp_i$ on $\domain$ using the extraction matrix
\begin{equation} \label{eq:matH-ellipse_3}
\periodic{\extMat} = \begin{bmatrix*}[c]
\frac{1}{2} & 1 & 0 & 0 & 0 & 0 & 0 & \frac{1}{2}\\
0 & 0 & 1 & \frac{1}{2} & \frac{1}{2} & 0 & 0 & 0\\
0 & 0 & 0 & \frac{1}{2} & \frac{1}{2} & 1 & 0 & 0 \\
\frac{1}{2} & 0 & 0 & 0 & 0 & 0 & 1 & \frac{1}{2}
\end{bmatrix*}.
\end{equation}
These basis functions are shown in Figure~\ref{fig:conic2} (b, top row).
Choosing the associated control points $\mbf{f}_i$ as
\begin{equation*}
\mbf{f}_1 = (2a_x,a_y) = -\mbf{f}_3,\quad
\mbf{f}_2 = (2a_x,-a_y) = -\mbf{f}_4,
\end{equation*}
we get a $C^1$ cubic description of an ellipse centred at $(0,0)$ with axis lengths $(a_x, a_y)$.
This can be verified in the vein of Remark~\ref{rmk:verification}.

\begin{example}\label{ex:ellipse_3}
Choosing $a_x = a_y = 1$, we obtain a circle of radius $1$, as shown in Figure~\ref{fig:conic2} (b, middle row). This $C^1$ cubic description is equivalent to the one used in \cite{toshniwal2017multi}.
The choice $a_x = 2a_y = 1$ yields an ellipse with axis lengths $(1, \frac{1}{2})$, as shown in Figure~\ref{fig:conic2} (b, bottom row).
The perturbed versions of these conics, with the control points chosen as in Equation \eqref{eq:conic2-perturb}, are shown as well and they remain clearly smooth.
It can be observed that, compared to the description from Section~\ref{ssec:ellipse2}, the control points here are at a greater distance from the curve.
This is completely analogous to the behavior of classical NURBS.
\end{example}

\subsubsection{$C^1$ description of multi-degree $(3,2,2)$}\label{ssec:ellipse322}
Many, many different $C^1$ low-degree representations of the circle can be cooked up.
Instead of attempting the impossible task of presenting them all, we present a single example that uses our splines in a more general setting than the above two descriptions.
Consider the domain $\domain = [0, \sqrt{2}+2]$ in the periodic setting.
Choose $\nseg = 3$ and
\begin{align*}
\seg{1}{\bknot} &= [0,0,0,0,\sqrt{2},\sqrt{2},\sqrt{2},\sqrt{2}], \quad
\seg{1}{\brwts} = \Bigl[1, \tfrac{1}{3},\tfrac{1}{3},1\Bigr],\\
\seg{2}{\bknot} &= \seg{3}{\bknot} = [0,0,0,1,1,1], \quad
\seg{2}{\brwts} = \seg{3}{\brwts} = \Bigl[1, \tfrac{\sqrt{2}}{2},1\Bigr].
\end{align*}
\begin{remark}
Opting for geometric continuity (see Remark~\ref{rem:geometric-continuity}) would have allowed the choice of $\domain$ with integral length, but the present setting is sufficient for illustrative purposes.
\end{remark}
Then, we can define four $C^1$ multi-degree piecewise-NURBS functions $\bsp_i$ on $\domain$ using the extraction matrix
\begin{equation*} 
\periodic{\extMat} = \begin{bmatrix*}[c]
\frac{1}{3} & 1 & 0 & 0 & 0 & 0 & 0 & 0 & 0 & \frac{1}{3}\\
0 & 0 & 1 & \frac{1}{3} & \frac{1}{3} & 0 & 0 & 0 & 0 & 0\\
0 & 0 & 0 & \frac{2}{3} & \frac{2}{3} & 1 & \frac{1}{2} & \frac{1}{2} & 0 & 0 \\
\frac{2}{3} & 0 & 0 & 0 & 0 & 0 & \frac{1}{2} & \frac{1}{2} & 1 & \frac{2}{3}
\end{bmatrix*}.
\end{equation*}
These basis functions are shown in Figure~\ref{fig:conic2} (c, top row).
Choosing the associated control points $\mbf{f}_i$ as
\begin{equation*}
\begin{gathered}
	\mbf{f}_1 = (2a_x,a_y),\quad
	\mbf{f}_2 = (2a_x,-a_y),\\
	\mbf{f}_3 = (-a_x,-a_y),\quad
	\mbf{f}_4 = (-a_x,a_y),
\end{gathered}
\end{equation*}
we get a $C^1$ multi-degree description of an ellipse centred at $(0,0)$ with axis lengths $(a_x, a_y)$.
This can be verified in the vein of Remark~\ref{rmk:verification}.

\begin{example}\label{ex:ellipse_322}
Choosing $a_x = a_y = 1$, we obtain a circle of radius $1$, as shown in Figure~\ref{fig:conic2} (c, middle row).
The choice $a_x = 2a_y = 1$ yields an ellipse with axis lengths $(1, \frac{1}{2})$, as shown in Figure~\ref{fig:conic2} (c, bottom row).
The perturbed versions of these conics, with the control points chosen as in Equation \eqref{eq:conic2-perturb}, are shown as well and they remain clearly smooth.
Once again, compared to the description from Section~\ref{ssec:ellipse2}, the control points here lie at a greater distance from the cubic portion of the curve. 
\end{example}

	\section{Piecewise-rational polar surfaces} \label{sec:rational-surfaces}

\begin{figure}[t!]
   \begin{subfigure}{0.49\textwidth}
      \centering
      \resizebox{0.75\textwidth}{!}{%
	\tikzsetnextfilename{./tikz/images/polar_map_horiz}%
	\begin{tikzpicture}[scale=0.4,line cap=round, line join=round]
	\def\L{6}
	\def\S{-12}

	\begin{scope}
		\draw [step=\L/4, black, ultra thin] (0,0) grid (\L, \L);
		\draw [red, ultra thick] (0,0) -- (\L, 0);
		\draw [black, ultra thick] (0,0) -- (0, \L);
		\draw [black, ultra thick] (\L,0) -- (\L, \L);
		
		\node (a) at (\L, \L/2){};
		
		\draw [black, ultra thin] (\L/2-\S, \L/2) circle (\L/2);
		\draw [black, ultra thin] (\L/2-\S, \L/2) circle (\L/4);
		\draw [black, ultra thin] (\L/2-\S, \L/2) circle (\L/8);
		\draw [black, ultra thin] (\L/2-\S, \L/2) circle (3*\L/8);
		\draw [black, ultra thin] (\L/2-\S, 0) -- (\L/2-\S, \L);
		\draw [black, ultra thin] (\L-\S, \L/2) -- (\L/2-\S, \L/2);
		\draw [black, ultra thick] (\L/2-\S, \L/2) -- (-\S, \L/2);
		\fill[red] (\L/2-\S, \L/2) circle (5pt);
		
		\node (b) at (-\S, \L/2){};
		
		\path [draw, -latex'] ($(a)!0.25!(b)$) -- ($(a)!0.75!(b)$); 
	
	\end{scope}
\end{tikzpicture}%
}
      \caption{Single polar singularity}
   \end{subfigure} \\[0.25cm]
   \begin{subfigure}{0.49\textwidth}
      \centering
      \resizebox{0.75\textwidth}{!}{%
	\tikzsetnextfilename{./tikz/images/bipolar_map}%
	\begin{tikzpicture}[scale=0.4,line cap=round, line join=round, 
	point/.style = {draw, circle, fill=black, inner sep=0.7pt},
	mark coordinate/.style={inner sep=0pt,outer sep=0pt,minimum size=3pt,
		fill=red,circle}%
]

	\def\L{6}
	\def\S{-15}
	
	\draw [step=\L*1.0/4, black, ultra thin] (\S,-\L/2) grid (\L+\S, \L-\L/2);
	\draw [red, ultra thick] (\S,0-\L/2) -- (\L+\S, 0-\L/2);
	\draw [black, ultra thick] (\S,0-\L/2) -- (0+\S, \L-\L/2);
	\draw [black, ultra thick] (\L+\S,0-\L/2) -- (\L+\S, \L-\L/2);
	\draw [red, ultra thick] (\S,\L-\L/2) -- (\L+\S, \L-\L/2);
	
	\node (a) at (\L+\S, 0){};

	\def\R{\L/2} 
	\def\angEl{35} 
	\filldraw[ball color=white, ultra thin] (0,0) circle (\R);
	\foreach \t in {-45,0,45} { \DrawLatitudeCircle[\R]{\t} }
	\foreach \t in {-10, -55, -100, -145} { \DrawLongitudeCircle[\R]{\t} }
	\foreach \t in {-100} { \DrawLongitudeCircleThick[\R]{\t} }
	
	\node (b) at (-\R, 0){};
	
	\path [draw, -latex'] ($(a)!0.25!(b)$) -- node [midway, above] {} ($(a)!0.75!(b)$);
	
\end{tikzpicture}%
}
      \caption{Double polar singularity}
   \end{subfigure}
   \caption{A single edge or a pair of opposite edges of a tensor-product spline patch can be collapsed for creating geometries with polar singularities. The collapsed edges here are shown in red, and the black edges are identified with each other to enforce periodicity.}
   \label{fig:polarMap}
\end{figure}
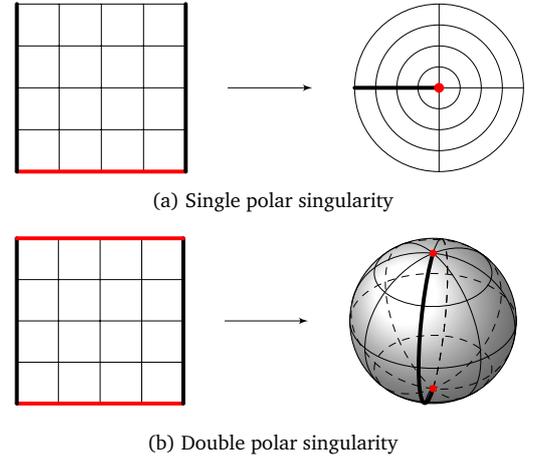

In this section, we describe how to construct $C^1$ smooth representations for \emph{polar surfaces} containing single or double polar singularities (e.g., hemispheres and spheres, respectively).
Such surfaces can be obtained by starting from a bivariate tensor-product (piecewise-NURBS) spline patch and collapsing one or two of its edges, respectively, as illustrated in Figure~\ref{fig:polarMap}.
Each of such edge collapses creates a \emph{polar point} and can be achieved by coalescing the control points related to basis functions with non-zero values on the edge.
In general, however, this control-point coalescing will introduce kinks at the poles and the surface representation will not be smooth.
To achieve overall smoothness, additional conditions need to be satisfied by the control points \cite[Section~3]{toshniwal2017multi}.
In Section~\ref{sec:polar_map}, we derive $C^1$ smoothness conditions at a polar point, and they enable us to build \emph{smooth polar splines} as linear combinations of bivariate tensor-product splines in Section~\ref{sec:polar-b-splines}. Then, in Section~\ref{sec:refinement-surfaces}, we give an explicit procedure how to compute a refined representation of a given polar surface.
Finally, in Section~\ref{sec:ellipsoid}, we present explicit descriptions of arbitrary ellipsoids using $C^1$ smooth low-degree polar spline representations suited for integrated design and analysis.

\subsection{Smoothness conditions at the polar points}\label{sec:polar_map}
A polar surface will be smooth at a polar point if it can be locally (re)parameterized in a smooth way.
Such \revision{parameterizations} can be specified in a constructive manner and we elaborate upon it in this section.
The resulting conditions will help us build smooth polar B-splines in the next section.

As shown in Figure~\ref{fig:polarMap}, we first describe the initial setup --- a tensor-product spline space on a rectangular domain.
We start from two univariate $C^1$ rational spline spaces $\splSpace^\pcoorda, \splSpace^\pcoordb$ defined on the univariate domains $\domain^\pcoorda := [\pcoorda_1, \pcoorda_2]$ and $\domain^\pcoordb := [\pcoordb_1, \pcoordb_2]$, respectively; the superscripts of $\pcoorda$ and $\pcoordb$ are meant to indicate the symbols used for the respective coordinates.
Using a Cartesian product, we build the rectangular domain $\domain := \domain^\pcoorda \times \domain^\pcoordb$, and on $\domain$ we define the tensor-product spline space $\splSpace := \splSpace^{\pcoorda} \otimes \splSpace^{\pcoordb}$.
Without loss of generality, we assume that $\pcoorda_1 = \pcoordb_1 = 0$.
This tensor-product spline space is spanned by tensor-product B-spline basis functions $\bsp_{ij}$, $i = 1, \dots, \ndof^\pcoorda$; $j = 1, \dots, \ndof^\pcoordb$.
Here, $\ndof^\pcoorda$ and $\ndof^\pcoordb$ denote the respective dimensions of the chosen univariate spline spaces; the basis functions spanning these spaces are denoted with $\bsp_i^\pcoorda$ and $\bsp_j^\pcoordb$.
Then, the tensor-product basis function $\bsp_{ij}$ is simply the product $\bsp_i^\pcoorda\bsp_j^\pcoordb$.
The functions $\bsp_{ij}$ are assumed to be periodic in $\pcoorda$ and non-periodic in $\pcoordb$.

Now, let us use the functions $\bsp_{ij}$ to map the domain $\domain$ to a polar surface using edge-collapse.
Then, the smoothness conditions at a collapsed edge will only involve those $\bsp_{ij}$ that have non-zero first derivatives there.
Observe that, if $\ndof^\pcoordb \geq 4$, then any $\bsp_{ij}$ with non-zero first derivatives at the bottom edge of $\domain$ will have zero first derivatives at its top edge, and vice versa.
The upshot is that, when we are collapsing both the bottom and top edges of $\domain$ into two polar points, as in Figure~\ref{fig:polarMap} (b), the smoothness conditions at those points are independent of each other and can be resolved separately for $\ndof^\pcoordb \geq 4$.

\subsubsection{Single polar point} \label{sec:pole-single}
In light of the above discussion, in the following we first focus on the case of a single collapsed edge, i.e., the one shown in Figure~\ref{fig:polarMap} (a).
We derive smoothness conditions that will help us build smooth polar spline functions (and, using them, smooth polar surfaces).
This is done by explicitly specifying the \revision{parameterization} {with respect to which} the spline functions are deemed smooth.
First, we construct a planar disk-like domain $\domainP$, called the \emph{polar parametric domain}, via a suitable polar map $\bpolarMap$; see Figure~\ref{fig:polarMap} (a).
Next, for an arbitrary $C^1$ spline $f \in \splSpace$, we define $\polar{f} : \domainP \rightarrow \RR$ to be $f \circ \bpolarMap^{-1}$. In general, $\polar{f}$ will be multivalued at the pole.
Finally, we derive the required smoothness conditions by asking for $\polar{f}$ to be $C^1$ smooth at the polar point.

We start by building $\bpolarMap$.
Assign the control point $\bpolarMapC_{ij}:=(\radius_j\cos(\theta_i),\radius_j\sin(\theta_i)) \in \RR^2$ to the basis function $\bsp_{ij}$, where
\HS{
\begin{equation}\label{eq:rho}
   \radius_j := \frac{j-1}{\ndof^\pcoordb-1} \in [0,1],
\end{equation}
and}
\begin{equation}\label{eq:theta}
   \theta_i := 2\pi + \frac{(1 - 2i)\pi}{\ndof^\pcoorda} \in [0,2\pi].
\end{equation}
The above choice of control-point values has been made in the interest of standardization and is not unique.
Using these control points, we can construct the disk-like domain $\domainP$ with the aid of the map 
$\bpolarMap$ from $\domain$ to $\domainP$,
\begin{equation*}
\domain \ni (s, t) \mapsto \bpolarMap(\pcoorda,\pcoordb) = (u,v) \in \domainP,
\end{equation*}
 defined as 
\begin{equation}\label{eq:degenMap}
\bpolarMap(\pcoorda,\pcoordb) := (\polarMap_u(\pcoorda,\pcoordb),\polarMap_v(\pcoorda,\pcoordb)) := \sum_{i=1}^{\ndof^\pcoorda} \sum_{j=1}^{\ndof^\pcoordb} \bpolarMapC_{ij} \bsp_{ij}(\pcoorda,\pcoordb).
\end{equation}
Note that the above construction will not necessarily yield an exactly circular domain $\domainP$; its shape will depend on the choice of $\splSpace$.
This domain will serve as the reference element for \DT{the smoothness of} polar configurations, i.e., we will define polar splines such that they are $C^1$ smooth functions over $\domainP$.
It is clear that for all $\pcoorda\in \domain^\pcoorda$,
\begin{equation}\label{eq:degenEdge}
\bpolarMap(\pcoorda,0) = (0,0),
\end{equation}
where $(0,0)\in \domainP$ is the polar point. Note that this implies
\begin{equation}\label{eq:degenEdgeDeriv}
  \frac{\partial \polarMap_u}{\partial \pcoorda}\bigg|_{\pcoordb = 0} \equiv 0 \equiv \frac{\partial \polarMap_v}{\partial \pcoorda}\bigg|_{\pcoordb = 0}.
\end{equation}

Let $\bspP_{ij}$ be the image of $\bsp_{ij}$ under the polar map 
$\bpolarMap:\domain \mapsto \domainP$ in Equation \eqref{eq:degenMap}
so that
\begin{equation*}
  \bspP_{ij}(u,v) = \bspP_{ij}(\bpolarMap(\pcoorda,\pcoordb)) = \bsp_{ij}(\pcoorda,\pcoordb).
\end{equation*}
Then, for given coefficients $f_{ij}$, a polar spline function $\polar{f}$ over $\domainP$ can be constructed as 
\begin{equation*}
  \polar{f}(u,v) = \sum_{i=1}^{\ndof^\pcoorda} \sum_{j=1}^{\ndof^\pcoordb} f_{ij} \bspP_{ij}(u,v).
\end{equation*}
We can pull $\polar{f}$ back to $\domain$ as follows,
\begin{equation}\label{eq:pulledBackPolarSplineVal}
\begin{split}
	f(\pcoorda,\pcoordb) := \polar{f}(\bpolarMap(\pcoorda,\pcoordb)) &= \sum_{i=1}^{\ndof^\pcoorda} \sum_{j=1}^{\ndof^\pcoordb} f_{ij} \bspP_{ij}\left(\bpolarMap(\pcoorda,\pcoordb)\right)\\
	&= \sum_{i=1}^{\ndof^\pcoorda} \sum_{j=1}^{\ndof^\pcoordb} f_{ij} \bsp_{ij}(\pcoorda,\pcoordb).
\end{split}
\end{equation}
Moreover, by using the chain rule we can also relate the partial derivatives of $f$ and $\polar{f}$:
\begin{equation*}
\begin{aligned}
\frac{\partial f}{\partial \pcoorda}(\pcoorda,\pcoordb) 
&= \frac{\partial \polar{f}}{\partial u}(u,v)\dfrac{\partial\polarMap_u}{\partial \pcoorda}(\pcoorda,\pcoordb) +  \frac{\partial \polar{f}}{\partial v}(u,v)\dfrac{\partial\polarMap_v}{\partial \pcoorda}(\pcoorda,\pcoordb),   \\
\frac{\partial f}{\partial \pcoordb}(\pcoorda,\pcoordb) 
&= \frac{\partial \polar{f}}{\partial u}(u,v)\dfrac{\partial\polarMap_u}{\partial \pcoordb}(\pcoorda,\pcoordb) +  \frac{\partial \polar{f}}{\partial v}(u,v)\dfrac{\partial\polarMap_v}{\partial \pcoordb}(\pcoorda,\pcoordb).
\end{aligned}
\end{equation*}
For $\polar{f}$ to be $C^1$ smooth at the polar point, there must exist real values $\alpha,\beta,\gamma$ such that
\begin{equation} \label{eq:C1-hermite}
\begin{gathered}
\lim\limits_{(u,v)\rightarrow(0,0)}
\left[\polar{f},\;\frac{\partial \polar{f}}{\partial u},\;\frac{\partial \polar{f}}{\partial v}\right](u,v) = [\alpha,\;\beta,\;\gamma].
\end{gathered}
\end{equation}
In view of \eqref{eq:degenEdge} and \eqref{eq:degenEdgeDeriv} this means for all $\pcoorda\in \domain^\pcoorda$,
\begin{equation*}
 f(s,0) = \alpha,
 \quad
\frac{\partial f}{\partial \pcoordb}(\pcoorda,0) 
= \beta\dfrac{\partial\polarMap_u}{\partial \pcoordb}(\pcoorda,0) +  \gamma\dfrac{\partial\polarMap_v}{\partial \pcoordb}(\pcoorda,0).
\end{equation*}
In particular, since only $\bsp_{ij}$, $j \leq 2$, have non-zero values and derivatives when $\pcoordb = 0$, the above condition translates to the following requirement for all $\pcoorda\in \domain^\pcoorda$,
\begin{equation}\label{eq:C1-constraints}
\begin{gathered}
	\sum_{i=1}^{\ndof^\pcoorda}f_{i1}\bsp_{i1}(\pcoorda,0) = \alpha,\\
	\sum_{i=1}^{\ndof^\pcoorda}\sum_{j=1}^{2}\left[f_{ij} - \bpolarMapC_{ij}\cdot(\beta,\gamma)\right]\frac{\partial \bsp_{ij}}{\partial \pcoordb}(\pcoorda,0) = 0.
\end{gathered}
\end{equation}

\subsubsection{Double polar point} \label{sec:pole-double}
Equation \eqref{eq:C1-constraints} shows the required smoothness conditions when the bottom edge of $\domain$ is being collapsed.
Next, if we also want to collapse the top edge of $\domain$, we can repeat the previous argument with minor changes.
We would, of course, need to choose a map $\ptop{\bpolarMap}$ that collapses the edge $\domain^\pcoorda \times \DT{\{\pcoordb_2\}}$ instead.
One way of achieving this could be by choosing the control points $\ptop{\bpolarMapC}_{ij}:=(\ptop{\radius}_j\cos(\ptop{\theta}_i),\ptop{\radius}_j\sin(\ptop{\theta}_i)) \in \RR^2$, where
\begin{equation*}
   \ptop{\radius}_j := 1 - \radius_j,
   \quad
   \ptop{\theta}_i := 2\pi - \theta_i.
\end{equation*}
Then, we can follow the same argument as in Section~\ref{sec:pole-single}. Asking for $C^1$ smoothness of $\polar{\ptop{f}}$ is equivalent to asking that there exist real values \HS{$\ptop{\alpha},\ptop{\beta},\ptop{\gamma}$} such that for all $\pcoorda \in \domain^\pcoorda$,
\begin{equation}\label{eq:C1-constraints-top}
\begin{gathered}
	\sum_{i=1}^{\ndof^\pcoorda}f_{i,\ndof^\pcoordb}\bsp_{i,\ndof^\pcoordb}(\pcoorda,\pcoordb_2) = \ptop{\alpha},\\
	\sum_{i=1}^{\ndof^\pcoorda}\sum_{j=\ndof^\pcoordb-1}^{\ndof^\pcoordb}\left[f_{ij} - \ptop{\bpolarMapC}_{ij}\cdot(\ptop{\beta}, \ptop{\gamma})\right]\frac{\partial \bsp_{ij}}{\partial \pcoordb}(\pcoorda,\pcoordb_2) = 0.
\end{gathered}
\end{equation}
Note once again that the smoothness at the polar point corresponding to $\pcoordb = 0$ is imposed with respect to the \revision{parameterization} $\bpolarMap(\domain)$, while that at the polar point corresponding to $\pcoordb = \pcoordb_2$ is imposed with respect to the \revision{parameterization} $\ptop{\bpolarMap}(\domain)$. The \DT{corresponding smoothness conditions in Equations \eqref{eq:C1-constraints} and \eqref{eq:C1-constraints-top} involve different coefficients} $f_{ij}$ for $\ndof^\pcoordb \geq 4$, and so can be resolved separately.

\begin{remark} \label{rmk:polar-local}
   The choices of $\theta_i, \radius_j, \ptop{\theta}_i, \ptop{\radius}_j$ are such that the maps $\bpolarMap$ and $\ptop{\bpolarMap}$ preserve the orientation of the parametric domain $\domain$.
\end{remark}


\subsection{Rational polar B-splines at the polar points} \label{sec:polar-b-splines}

We now elaborate how the derived $C^1$ smoothness constraints at a polar point will enable the computation of a \DTA{} extraction matrix.
This matrix represents a linear map to a set of polar spline basis functions that are $C^1$ smooth on the polar parametric domain.

\subsubsection{Single polar point} \label{sec:polar-spline-single}
As before, we start by considering the case of a single collapsed edge, i.e., the one shown in Figure~\ref{fig:polarMap} (a).
Let us arrange {the set of basis functions $\{\bspP_{ij}:i=1,\dots,\ndof^\pcoorda; j=1,\dots,\ndof^\pcoordb\}$} in a vector $\mbf{\bsp}$, where $\bspP_{ij}$ occupies the $(i + (j-1)\ndof^\pcoorda)$-th entry. Our goal is to construct $C^1$ smooth polar basis functions on the polar parametric domain $\polar{\domain}$ as suitable linear combinations of the functions $\bspP_{ij}$. In other words, we are looking for an extraction matrix $\extMatP$ such that the polar spline basis functions in $\{\polar{N}_{\polari}:\polari=1,\dots,\ndof\}$ defined by the following relation,
\begin{equation} \label{eq:polarExtraction}
  \mbf{N} := \extMatP\mbf{\bsp},
\end{equation}
are $C^1$ at the polar point.
For fixed $j$, the set $\{\bspP_{ij}:i=1,\dots,\ndof^\pcoorda\}$ is called the $(j-1)$-th polar ring of basis functions.
When $j>2$, all basis functions in the $(j-1)$-th ring already satisfy the $C^1$ continuity conditions at the polar point (their derivatives are identically zero there), so they can be included without modifications in the set of polar basis functions being created. The others will be substituted by three smooth polar basis functions.
This dictates that $\extMatP$ will be a matrix, with $\ndof := \ndof^\pcoorda (\ndof^\pcoordb - 2) + 3$ rows and $\ndof^\pcoorda \ndof^\pcoordb$ columns, taking the following sparse block-diagonal form:
\begin{equation} \label{eq:polarC1Extraction-single}
\extMatP := \begin{bmatrix*}[c]
  \reduced{\extMatP} & \\
   & \identMat
  \end{bmatrix*},
\end{equation}
where $\identMat$ is the identity matrix of size $\ndof^\pcoorda (\ndof^\pcoordb - 2) \times \ndof^\pcoorda (\ndof^\pcoordb - 2)$ and $\reduced{\extMatP}$ is a matrix of size $3 \times 2\ndof^\pcoorda$. The entry of $\reduced{\extMatP}$ corresponding to its $\polari$-th row and $\left(i + (j-1)\ndof^\pcoorda\right)$-th column is denoted with $\reduced{\extMatPel}_{\polari,(ij)}$. We can then rewrite Equation \eqref{eq:polarExtraction} as follows for $\polari = 1,2,3$,
\begin{equation*} 
\polar{N}_\polari(u,v) = \sum_{i=1}^{\ndof^\pcoorda}\sum_{j=1}^{2}\reduced{\extMatPel}_{\polari,(ij)}\bspP_{ij}(u,v).
\end{equation*}
We can pull these back to $\domain$ using Equation \eqref{eq:pulledBackPolarSplineVal} to obtain the equivalent representation for $\polari = 1,2,3$,~
\begin{equation} \label{eq:polarBasis_expansion_pulledBack}
N_\polari(\pcoorda,\pcoordb) = \sum_{i=1}^{\ndof^\pcoorda}\sum_{j=1}^{2}\reduced{\extMatPel}_{\polari,(ij)}\bsp_{ij}(\pcoorda,\pcoordb).
\end{equation}

We will enforce $C^1$ continuity at the polar point by requiring the basis functions $\polar{N}_\polari$ to satisfy a linearly independent Hermite data set at the polar point, in the spirit of Equation \eqref{eq:C1-hermite}. To this end, we will use three \emph{source basis functions} {$\{T_\polari:\polari=1,2,3\}$}, that provide us with the appropriate Hermite data. 
Given a non-degenerate triangle $\standsimp$ with vertices $\mbf{v}_1,\mbf{v}_2$ and $\mbf{v}_3$, let $(\lambda_1,\lambda_2,\lambda_3)$ be the unique barycentric coordinates of point $(u,v)$ with respect to $\standsimp$ such that
\begin{equation*}
   \lambda_1\mbf{v}_1+\lambda_2\mbf{v}_2+\lambda_3\mbf{v}_3 = (u,v), \quad \lambda_1+\lambda_2+\lambda_3 = 1.
\end{equation*}
Then, we define
\begin{equation*}
 T_\polari(u,v):=\lambda_\polari, \quad \polari=1,2,3.
\end{equation*}
These functions can be interpreted as triangular Bernstein polynomials of degree $1$. They are non-negative on the domain triangle $\standsimp$. Moreover, they are linearly independent, form a partition of unity, and span the space of bivariate polynomials of total degree less than or equal to $1$.
Then, we require that $N_\polari$ in Equation \eqref{eq:polarBasis_expansion_pulledBack} is a spline function $f$ such that it satisfies the continuity constraints in Equation \eqref{eq:C1-constraints}, with
\begin{equation*} 
 \alpha = T_\polari(0,0),
 \quad
\beta = \frac{\partial T_\polari}{\partial u}(0,0), \quad
\gamma = \frac{\partial T_\polari}{\partial v}(0,0),
\end{equation*}
for $\polari=1,2,3$.
In the interest of standardization, we choose the triangle $\standsimp$ as equilateral with vertices
\DT{
\begin{equation*}
	\mbf{v}_1 = (2\radius_2, 0),\ 
	\mbf{v}_2 = (-\radius_2, \sqrt{3}\radius_2),\ 
	\mbf{v}_3 = (-\radius_2, -\sqrt{3}\radius_2);
\end{equation*}
recall the definition of $\radius_2$ from Equation~\eqref{eq:rho}.}
After some calculations, we deduce that 
\begin{equation*}
\reduced{\extMatPel}_{1,(i1)}=
\reduced{\extMatPel}_{2,(i1)}=
\reduced{\extMatPel}_{3,(i1)}=\frac{1}{3},
\end{equation*}
and
\begin{equation*}
\begin{bmatrix*}[c]
\reduced{\extMatPel}_{1,(i2)}\\
\reduced{\extMatPel}_{2,(i2)}\\
\reduced{\extMatPel}_{3,(i2)}
\end{bmatrix*} = 
\begin{bmatrix*}[r]
\frac{1}{3} & 0 & \frac{1}{3}\\
-\frac{1}{6} & \frac{\sqrt{3}}{6} & \frac{1}{3}\\
-\frac{1}{6} & -\frac{\sqrt{3}}{6} & \frac{1}{3}
\end{bmatrix*}
\begin{bmatrix*}[c]
\cos(\theta_i)\\
\sin(\theta_i)\\
1
\end{bmatrix*}.
\end{equation*}
This relation says that 
$\left(\reduced{\extMatPel}_{1,(i2)},\reduced{\extMatPel}_{2,(i2)},\reduced{\extMatPel}_{3,(i2)}\right)$
are simply the barycentric coordinates of the control point $\bpolarMapC_{i2}:=(\radius_2\cos(\theta_i),\radius_2\sin(\theta_i))$ with respect to $\standsimp$.
It is easily checked that $\standsimp$  encloses the circle centred at $(0,0)$ with a radius of $\radius_2$, and hence $\left(\reduced{\extMatPel}_{1,(i2)},\reduced{\extMatPel}_{2,(i2)},\reduced{\extMatPel}_{3,(i2)}\right)$ are guaranteed to be non-negative.
In summary, $\reduced{\extMatP}$ is specified as
\begin{equation} \label{eq:polarC1Extraction}
  \reduced{\extMatP} := \begin{bmatrix*}[c]
  \frac{1}{3} & \cdots & \frac{1}{3} & \reduced{\extMatPel}_{1,(12)} & \cdots & \reduced{\extMatPel}_{1,(i2)} & \cdots & \reduced{\extMatPel}_{1,(\ndof^\pcoorda 2)}\\
  \frac{1}{3} & \cdots & \frac{1}{3} & \reduced{\extMatPel}_{2,(12)} & \cdots & \reduced{\extMatPel}_{2,(i2)} & \cdots & \reduced{\extMatPel}_{2,(\ndof^\pcoorda 2)}\\
  \frac{1}{3} & \cdots & \frac{1}{3} & \reduced{\extMatPel}_{3,(12)} & \cdots & \reduced{\extMatPel}_{3,(i2)} & \cdots & \reduced{\extMatPel}_{3,(\ndof^\pcoorda 2)}
  \end{bmatrix*}.
\end{equation}
This matrix has full rank and the column sum is equal to one, thus confirming that $\extMatP$ is \DTA. The following result follows from the above discussion.

\begin{theorem}\label{thm:polar_properties}
  The $C^1$ smooth polar spline functions in the set $\{\polar{N}_{\polari}:\polari=1,\dots,\ndof\}$ are linearly independent, locally supported, and form a convex partition of unity on $\domainP$.
\end{theorem} 

\begin{remark}
  As long as $\standsimp$ is chosen to be a triangle enclosing the first polar ring of control points $\bpolarMapC_{ij}$ for a given configuration, we are guaranteed non-negative extraction coefficients.
  It is only in the interest of standardization that we have chosen to fix $\standsimp$ as an equilateral triangle with a fixed pattern of vertices. 
\end{remark}


Given $\ndof$ control points $\mbf{f}_\polari \in \RR^{\spd}$, $\spd\geq3$, we can construct a  $C^1$ polar surface $\mbf{f}$ embedded in $\RR^{\spd}$,
\begin{equation*}
 \mbf{f}(u,v) = \sum_{\polari=1}^{\ndof} \mbf{f}_\polari \polar{N}_\polari(u,v),
\end{equation*}
or, equivalently, after pulling back to $\domain$,
\begin{equation} \label{eq:polarSurface}
 \mbf{f}(\pcoorda,\pcoordb) = \sum_{\polari=1}^{\ndof} \mbf{f}_\polari N_\polari(\pcoorda,\pcoordb).
\end{equation}
The behavior of $\mbf{f}$ at the polar point is going to be fully specified by the first three control points $\mbf{f}_1$, $\mbf{f}_2$ and $\mbf{f}_3$. These control points can be thought of as forming a \emph{control triangle}.
For a fixed surface, the transpose of the extraction matrix $\extMatP$ defines the relationship between control points of the $\bsp_{ij}$ and control points of the smooth $N_\polari$. More precisely, if
\begin{equation*}
  \sum_{i=1}^{\ndof^\pcoorda}\sum_{j=1}^{\ndof^\pcoordb}\mbf{g}_{ij}\bsp_{ij}(\pcoorda,\pcoordb) = \mbf{f}(\pcoorda,\pcoordb) = \sum_{\polari=1}^{\ndof}\mbf{f}_{\polari}N_{\polari}(\pcoorda,\pcoordb),
\end{equation*}
then
\begin{equation} \label{eq:polarControlPointExtraction}
  \mbf{g}_{ij} = \sum_{\polari=1}^\ndof \extMatPel_{\polari,(ij)}\mbf{f}_\polari.
\end{equation}
In particular, the transpose of the extraction matrix $\reduced{\extMatP}$ in Equation \eqref{eq:polarC1Extraction} essentially computes the control points of the zeroth and the first polar rings of basis functions $\bsp_{ij}$ as convex combinations of $\mbf{f}_1$, $\mbf{f}_2$ and $\mbf{f}_3$.
In other words, it forces the zeroth and the first polar rings of the control points $\mbf{g}_{ij}$ to be coplanar.
The plane in which they lie is the one passing through $\mbf{f}_1$, $\mbf{f}_2$ and $\mbf{f}_3$, and hence the control triangle must be tangent to the surface $\mbf{f}$ at the pole.

In light of \HS{Theorem~\ref{thm:polar_properties}} we can immediately conclude the following properties for the pulled-back functions $N_{\polari}$, $\polari=1,\dots,\ndof$.

\begin{corollary}\label{col:polar_properties_pullback}
  The spline functions in the set $\{N_{\polari}:\polari=1,\dots,\ndof\}$ are linearly independent, locally supported, and form a convex partition of unity on $\domain$. Moreover, any polar surface $\mbf{f}$ as in Equation \eqref{eq:polarSurface} with non-collinear control points $\mbf{f}_1$, $\mbf{f}_2$ and $\mbf{f}_3$ will have a well-defined tangent plane at the pole.
\end{corollary}

\begin{remark}
  According to Section~\ref{sec:mdb-splines}, the univariate sets of basis functions $\{\bsp_i^\pcoorda:i=1, \dots,\ndof^\pcoorda\}$ and $\{\bsp_j^\pcoordb:j=1,\dots,\ndof^\pcoordb\}$ are built from local NURBS basis functions through extraction matrices $\extMat^\pcoorda$ and $\extMat^\pcoordb$, respectively. Hence, by combining the matrices $\extMatP$, $\extMat^\pcoorda$ and $\extMat^\pcoordb$, the spline basis functions in the set $\{N_{\polari}:\polari=1,\dots,\ndof\}$ can be directly expressed in terms of local tensor-product NURBS basis functions. 
\end{remark}


\subsubsection{Double polar point} \label{sec:polar-spline-double}

When dealing with double polar surfaces, the spline construction can be obtained by collapsing a pair of two opposite edges as illustrated in Figure~\ref{fig:polarMap} (b). As explained in Section~\ref{sec:polar_map}, the smoothness treatment of the two poles can be done separately for $\ndof^\pcoordb\geq 4$.
In this case, each pole leads to a local extraction matrix by applying the same procedure as in Section~\ref{sec:polar-spline-single} and the combined global extraction matrix takes the following sparse block-diagonal form:
\begin{equation}  \label{eq:polarC1Extraction-double}
\extMatP := \begin{bmatrix*}[c]
  \seg{1}{\reduced{\extMatP}} & & \\
   & \identMat & \\
   & & \seg{2}{\reduced{\extMatP}}
  \end{bmatrix*},
\end{equation}
where $\identMat$ is the identity matrix of size $\ndof^\pcoorda (\ndof^\pcoordb - 4) \times \ndof^\pcoorda (\ndof^\pcoordb - 4)$ and $\seg{i}{\reduced{\extMatP}}$, $i=1,2$, are matrices of size $3 \times 2\ndof^\pcoorda$. 
By choosing the two polar \revision{parameterizations} $\bpolarMap(\domain)$ and $\ptop{\bpolarMap}(\domain)$ specified in Section~\ref{sec:polar_map}, it is easily verified that one can set
\begin{equation} \label{eq:polarC1Extraction-double-sub}
  \seg{1}{\reduced{\extMatP}} := \reduced{\extMatP}, \quad 
  \seg{2}{\reduced{\extMatP}} := \exchangeMat_3\reduced{\extMatP}\exchangeMat_{2\ndof^\pcoorda},
\end{equation}
where $\reduced{\extMatP}$ is the matrix defined in Equation \eqref{eq:polarC1Extraction} and $\exchangeMat_k$ is the exchange matrix of size $k \times k$, i.e., an anti-diagonal matrix of the form
\begin{equation*}
\exchangeMat_k := 
\begin{bmatrix*}[c]
 &  & & 1 \\
 & &\iddots& \\
 & 1 & & \\
1 & & &
\end{bmatrix*}.
\end{equation*}
The extraction matrix $\extMatP$ can then be used to compute the set of spline functions $\{N_{\polari}:\polari=1,\dots,\ndof\}$ in terms of the tensor-product functions $\{\bsp_{ij}:i=1,\dots,\ndof^\pcoorda; j=1,\dots,\ndof^\pcoordb\}$.
Similar to the single-pole result in Corollary~\ref{col:polar_properties_pullback}, these spline functions have the following properties.

\begin{corollary}
  The spline functions in the set $\{N_{\polari}:\polari=1,\dots,\ndof\}$ are linearly independent, locally supported, and form a convex partition of unity on $\domain$. Moreover, any polar surface $\mbf{f}$ as in Equation \eqref{eq:polarSurface} with non-collinear control points $\mbf{f}_1$, $\mbf{f}_2$ and $\mbf{f}_3$ and with non-collinear control points $\mbf{f}_{\ndof-2}$, $\mbf{f}_{\ndof-1}$ and $\mbf{f}_{\ndof}$ will have a well-defined tangent plane at both poles.
\end{corollary}

\subsection{Refinement of piecewise-rational polar surfaces}\label{sec:refinement-surfaces}

A polar spline surface can be refined in a manner similar to the one discussed in Section~\ref{sec:refinement-curves}.
We begin with the following observation. Consider the $3\times3$ matrix
\begin{equation*}
\reduced{\mbf{M}}:=\begin{bmatrix*}[r]
\frac{1}{3} & 0 & \frac{1}{3}\\
-\frac{1}{6} & \frac{\sqrt{3}}{6} & \frac{1}{3}\\
-\frac{1}{6} & -\frac{\sqrt{3}}{6} & \frac{1}{3}
\end{bmatrix*}
\begin{bmatrix*}[c]
0 & \cos(\theta_{\iota}) & \cos(\theta_{\kappa})\\
0 & \sin(\theta_{\iota}) & \sin(\theta_{\kappa})\\
1 & 1 & 1
\end{bmatrix*},
\end{equation*}
for some \HS{angles $\theta_{\iota}$ and $\theta_{\kappa}\not\in\{\theta_{\iota},\theta_{\iota}+\pi\}$ selected from the set in Equation \eqref{eq:theta}.} This a submatrix of $\reduced{\extMatP}$, defined in Equation \eqref{eq:polarC1Extraction}, consisting of three linearly independent columns. Its inverse is given by
\begin{align*}
\reduced{\mbf{M}}^{-1} &= \begin{bmatrix*}[c]
0 & \cos(\theta_{\iota}) & \cos(\theta_{\kappa})\\
0 & \sin(\theta_{\iota}) & \sin(\theta_{\kappa})\\
1 & 1 & 1
\end{bmatrix*}^{-1}
\begin{bmatrix*}[r]
\frac{1}{3} & 0 & \frac{1}{3}\\
-\frac{1}{6} & \frac{\sqrt{3}}{6} & \frac{1}{3}\\
-\frac{1}{6} & -\frac{\sqrt{3}}{6} & \frac{1}{3}
\end{bmatrix*}^{-1} 
\\
&=\frac{1}{\sin(\theta_{\kappa}-\theta_{\iota})}
\reduced{\mbf{L}}
\begin{bmatrix*}[r]
2 & -1 & -1\\
0 & \sqrt{3} & -\sqrt{3}\\
1 & 1 & 1
\end{bmatrix*},
\end{align*}
where
\begin{equation*}
\reduced{\mbf{L}} := \begin{bmatrix*}[c]
\sin(\theta_{\iota})-\sin(\theta_{\kappa}) & \cos(\theta_{\kappa})-\cos(\theta_{\iota}) & \sin(\theta_{\kappa}-\theta_{\iota})\\
\sin(\theta_{\kappa}) & -\cos(\theta_{\kappa}) & 0\\
-\sin(\theta_{\iota}) & \cos(\theta_{\iota}) & 0
\end{bmatrix*}.
\end{equation*}
Then, we define a sparse matrix $\reduced{\extMatPinv}$ of size $2\ndof^\pcoorda \times 3$, whose non-zero entries $\reduced{\extMatPinvel}_{ij}$ are identified as follows: for $j=1,2,3$,
\begin{equation*}
\reduced{\extMatPinvel}_{1j} := \reduced{M}^{-1}_{1j}, \quad
\reduced{\extMatPinvel}_{\ndof^\pcoorda+\iota,j} := \reduced{M}^{-1}_{2j}, \quad
\reduced{\extMatPinvel}_{\ndof^\pcoorda+\kappa,j} := \reduced{M}^{-1}_{2j}.
\end{equation*}
From its construction it is clear that the product $\reduced{\extMatP}\reduced{\extMatPinv}$ is equal to the identity matrix. Note that the product $(\exchangeMat_3\reduced{\extMatP}\exchangeMat_{2\ndof^\pcoorda})(\exchangeMat_{2\ndof^\pcoorda}\reduced{\extMatPinv}\exchangeMat_3)$ is also equal to the identity matrix. Hence, keeping in mind the definition of $\extMatP$ in Equations \eqref{eq:polarC1Extraction-double}--\eqref{eq:polarC1Extraction-double-sub}, the matrix
\begin{equation} \label{eq:polarC1Extraction-inv}
\extMatPinv := \begin{bmatrix*}[c]
\reduced{\extMatPinv} & & \\
& \identMat & \\
& & \exchangeMat_{2\ndof^\pcoorda}\reduced{\extMatPinv}\exchangeMat_3
\end{bmatrix*}
\end{equation}
gives rise to a product $\extMatP\extMatPinv$ that is equal to the identity matrix.
A similar matrix $\extMatPinv$ can be found (with only two diagonal blocks) for the matrix $\extMatP$ defined in Equation \eqref{eq:polarC1Extraction-single}.

Any refinement matrix of polar splines can then be built using the following procedure.
First, we compute the control points of the tensor-product basis functions $\bsp_{ij}$ as in Equation \eqref{eq:polarControlPointExtraction}.
Then, we refine the tensor-product control points using a tensor product of univariate refinement matrices (see Section~\ref{sec:refinement-curves}).
Denote this matrix with $\refMatloc$, and denote the polar spline extraction matrices before and after refinement with $\extMatP$ and $\refined{\extMatP}$, respectively.
Then, control points of the refined polar spline basis functions can be obtained by applying a matrix $\refMat$ to the original set of polar control points, where $\refMat$ is computed by solving the following (overdetermined) linear system with a unique solution,
\begin{equation*}
\refMat\refined{\extMatP} = \extMatP\refMatloc.
\end{equation*}
After multiplication of both sides of this system with the matrix $\refined{\extMatPinv}$ (corresponding to $\refined{\extMatP}$) as defined in Equation \eqref{eq:polarC1Extraction-inv}, we arrive at
\begin{equation} \label{eq:ref-system-polar}
\refMat = \extMatP\refMatloc\refined{\extMatPinv}.
\end{equation}

\begin{remark} \label{rmk:choice-D}
	Similar to Remark~\ref{rmk:choice-G}, there is some flexibility in the definition of $\refined{\extMatPinv}$ for the computation of $\refMat$ in Equation \eqref{eq:ref-system-polar}, as long as it is a right inverse of $\refined{\extMatP}$. The current choice is again done for the sake of simplicity. For the computation of $\reduced{\mbf{M}}^{-1}$, from a numerical point of view, it is advised to select the angles $\theta_{\iota}$ and $\theta_{\kappa}$ from the set in Equation \eqref{eq:theta} such that their difference is as close as possible to $\pm\pi/2$. For instance, taking $\iota=\kappa+\lfloor \ndof^\pcoorda/4 + 1/2\rfloor$ gives 
	\begin{equation*}
	\theta_{\kappa}-\theta_{\iota}=\frac{2\pi}{\ndof^\pcoorda}\left\lfloor \frac{\ndof^\pcoorda}{4} + \frac{1}{2}\right\rfloor,
	\end{equation*}
	and then $\det(\reduced{\mbf{M}})=\sin(\theta_{\kappa}-\theta_{\iota})\geq 0.86$ for $\ndof^\pcoorda\geq3$.
\end{remark}

\subsection{Ellipsoids and spheres} \label{sec:ellipsoid}

Let us now present explicit descriptions of ellipsoids (and as a special case also spheres) built using the $C^1$ polar spline framework discussed thus far.
Once again, in the interest of providing the simplest possible representations, we focus on smooth descriptions of low(est) degree only.
More precisely, we present three explicit descriptions of arbitrary ellipsoids using smooth polar splines.
The first one uses $8$ rational pieces of bi-degree $(2,2)$; the second one uses $4$ rational pieces of bi-degree $(2,3)$; and the last one uses $2$ rational pieces of bi-degree $(3,3)$.
All approaches will define six $C^1$ polar spline functions $N_\polari$ and associated control points $\mbf{f}_\polari$, $\polari = 1, \dots, 6$, such that the bivariate surface $\mbf{f}$,
\begin{equation*}
\mbf{f}(\pcoorda,\pcoordb) := (f_x(\pcoorda,\pcoordb), f_y(\pcoorda,\pcoordb), f_z(\pcoorda,\pcoordb)) := \sum_{\polari = 1}^6 \mbf{f}_\polari N_\polari(\pcoorda,\pcoordb),
\end{equation*}
describes the exact ellipsoid centred at $(0,0,0)$ and with axis lengths $(a_x, a_y, a_z)$, 
\begin{equation} \label{eq:ellipsoid}
   \left(\frac{f_x}{a_x}\right)^2
   + \left(\frac{f_y}{a_y}\right)^2
   + \left(\frac{f_z}{a_z}\right)^2
      = 1.
\end{equation}
To this end, we will make repeated use of the identity matrix $\identMat_k$ of size $k \times k$ and the exchange matrix $\exchangeMat_k$ of size $k \times k$.

To visually illustrate the smoothness of $\mbf{f}$, we will also show the surface $\tilde{\mbf{f}}$ obtained by perturbing one of the control points. This surface will also be smooth at the poles. For uniformity throughout the examples, we will choose the control points of the perturbed \DT{surface} as
\begin{equation}\label{eq:conic3-perturb}
\tilde{\mbf{f}}_\polari := 
\begin{cases}
   \mbf{f}_\polari + (0,0,4a_x), & \polari = 3,\\
   \mbf{f}_\polari, & \polari \neq3.
\end{cases}
\end{equation}

\tikzexternaldisable

\mdfsetup{%
   middlelinecolor=black!14,
   middlelinewidth=1pt,
   roundcorner=10pt}

\begin{sidewaysfigure*}[p]
   \centering
   \begin{subfigure}{0.32\linewidth}
      \begin{mdframed}
         \hspace*{1cm}\includegraphics[trim=630 370 70 240, clip, width=0.9\linewidth]{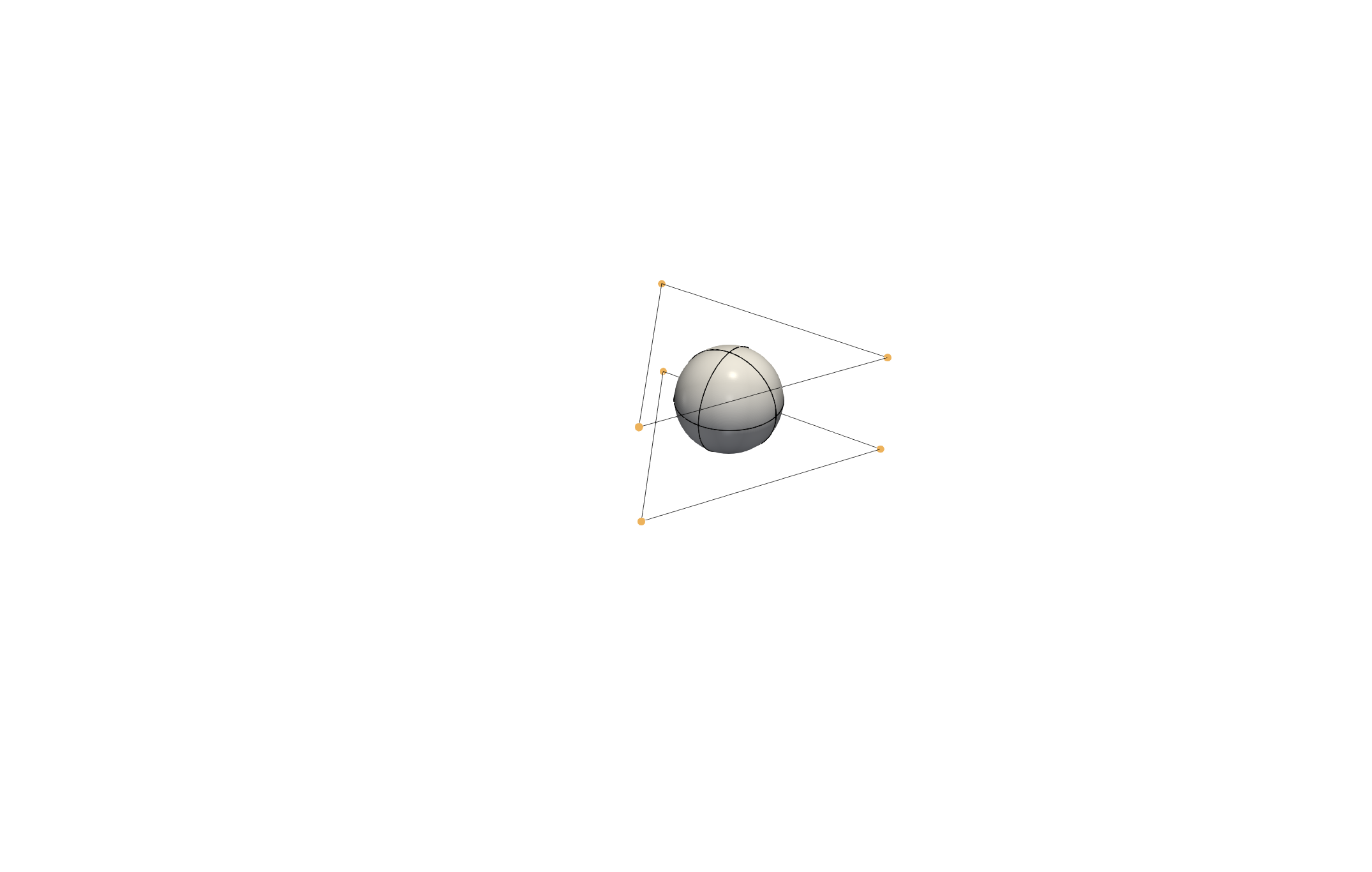}\\
         \includegraphics[trim=170 280 340 400, clip, width=0.9\linewidth]{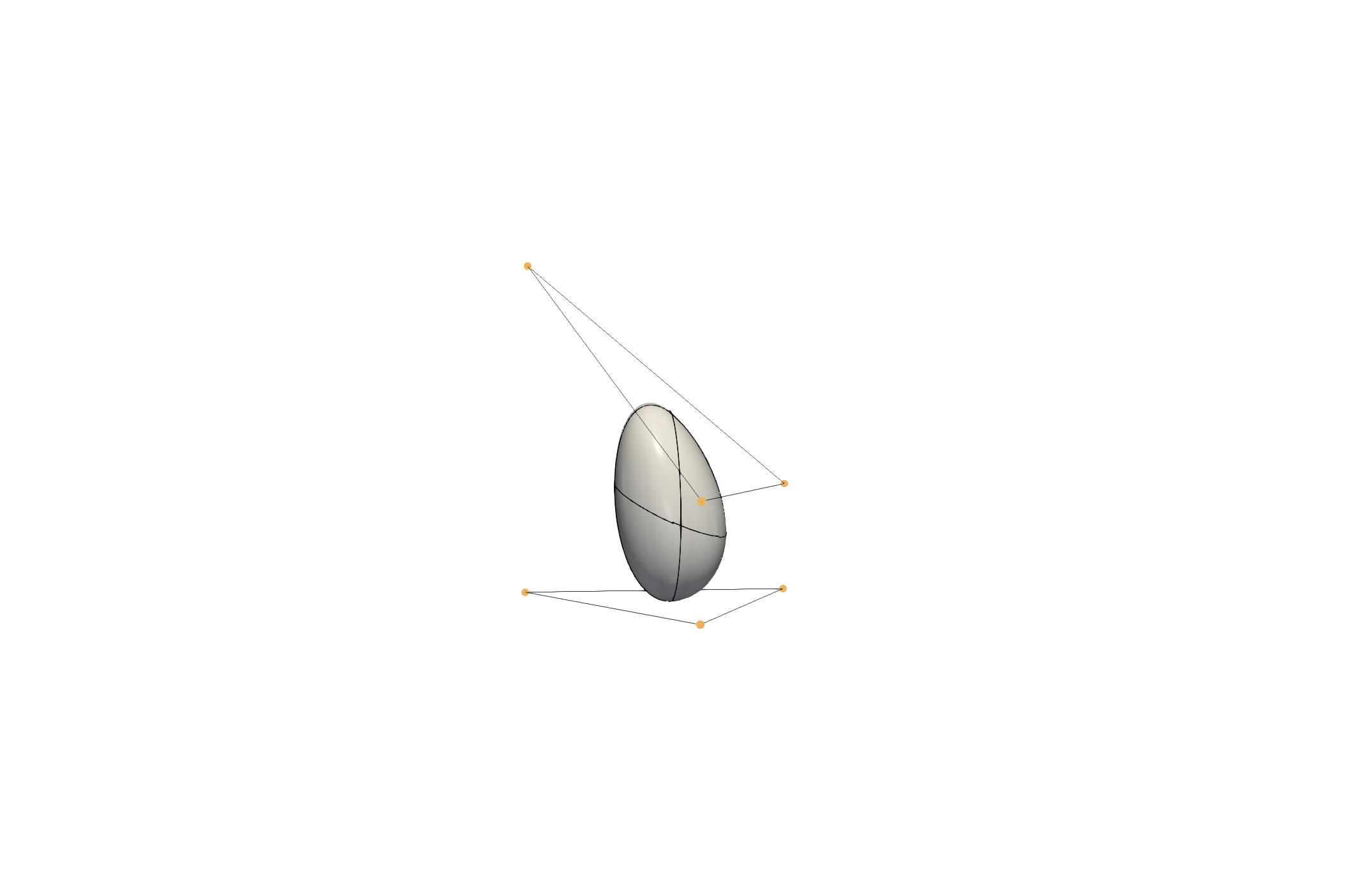}
      \end{mdframed}\vspace*{-0.4cm}
      \begin{mdframed}
         \hspace*{1cm}\includegraphics[trim=600 370 130 320, clip, width=0.9\linewidth]{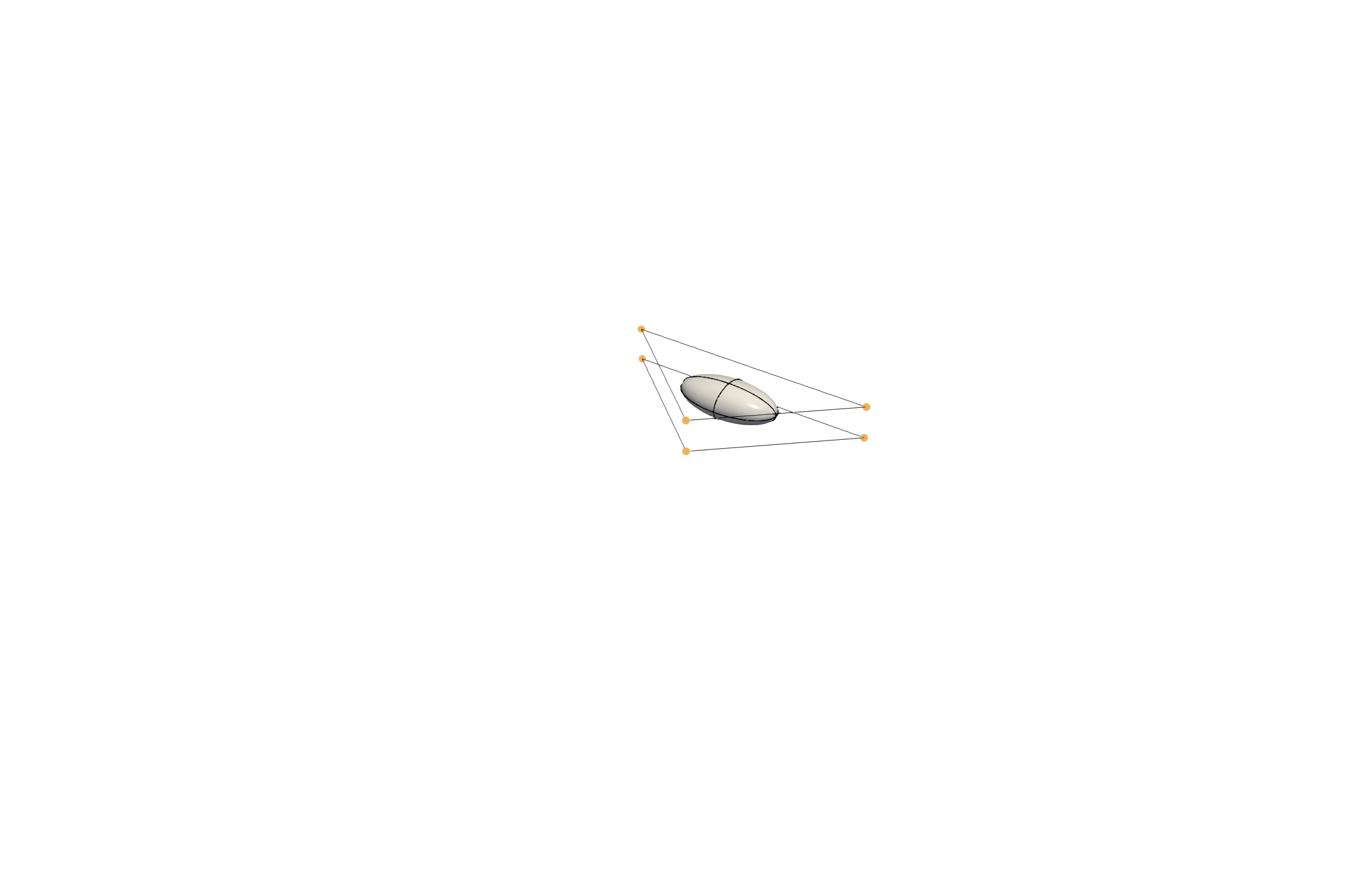}\\
         \includegraphics[trim=140 280 350 420, clip, width=0.9\linewidth]{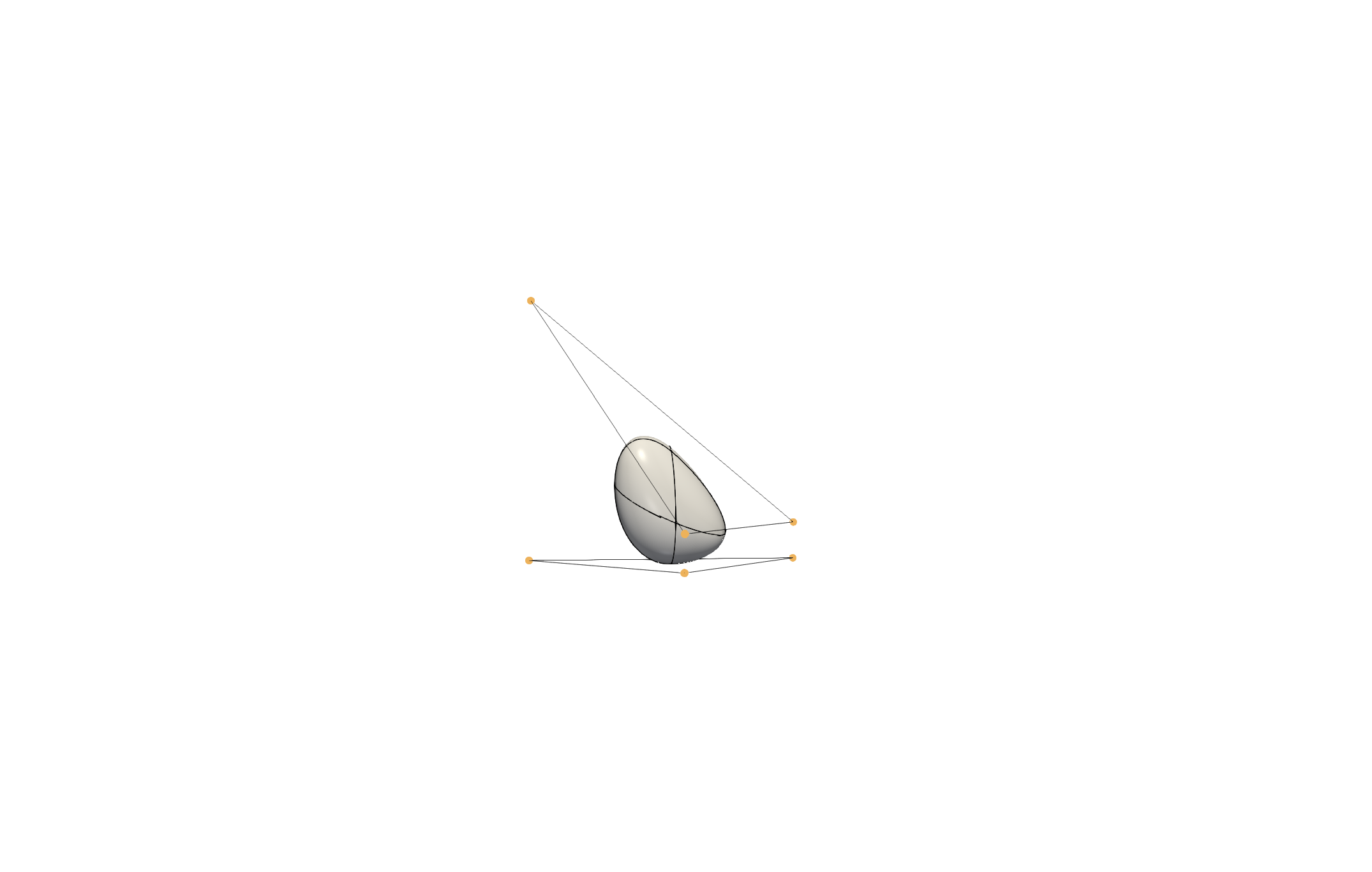}
      \end{mdframed}\vspace*{-0.4cm}
      \caption{degree $(2,2)$}
   \end{subfigure}\hspace*{0.05cm}
   \begin{subfigure}{0.32\linewidth}
      \begin{mdframed}
         \hspace*{0.4cm}\includegraphics[trim=630 370 70 240, clip, width=0.9\linewidth]{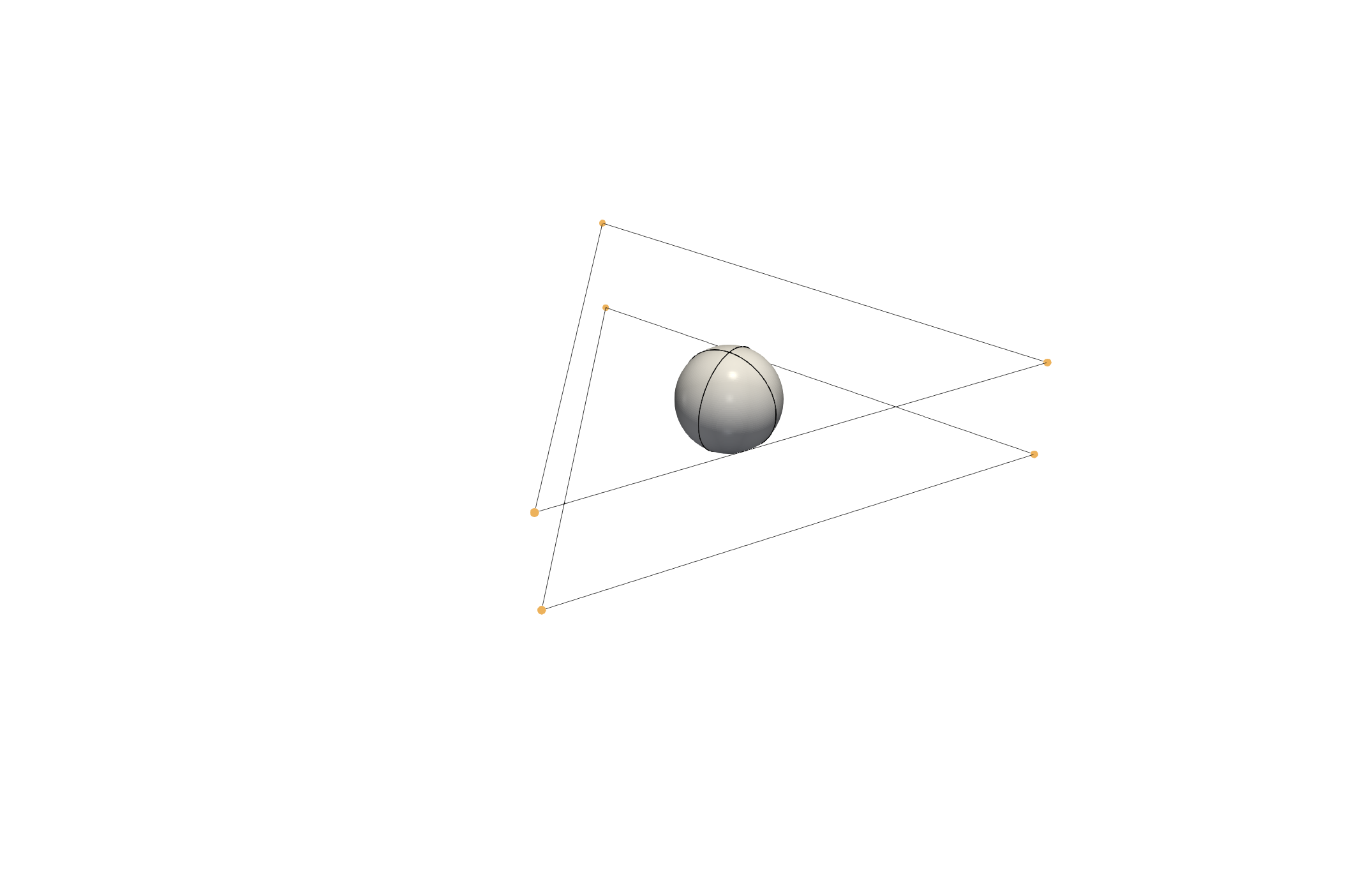}\\
         \includegraphics[trim=170 280 340 400, clip, width=0.9\linewidth]{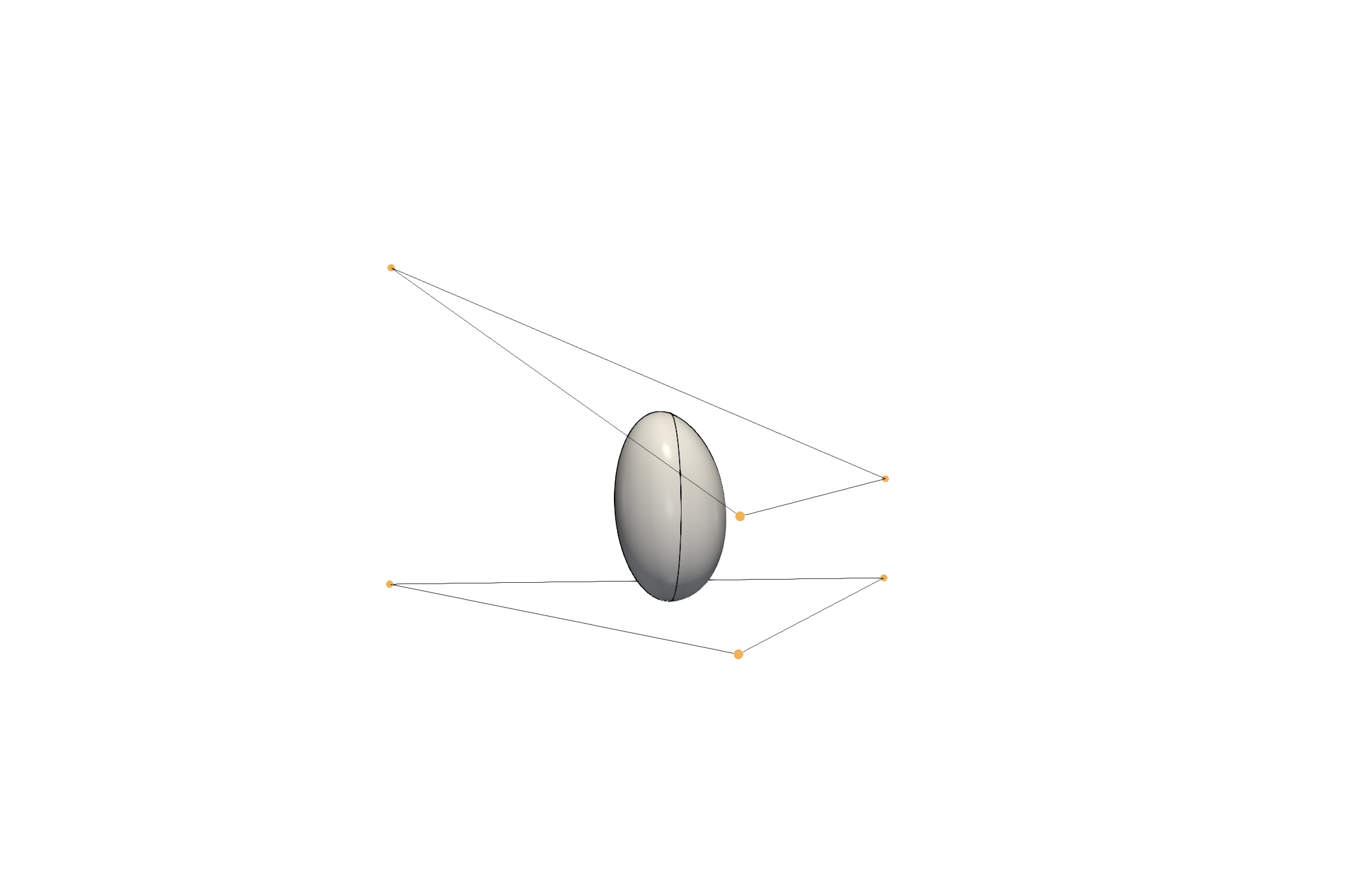}
      \end{mdframed}\vspace*{-0.4cm}
      \begin{mdframed}
         \hspace*{0.3cm}\includegraphics[trim=600 370 130 320, clip, width=0.9\linewidth]{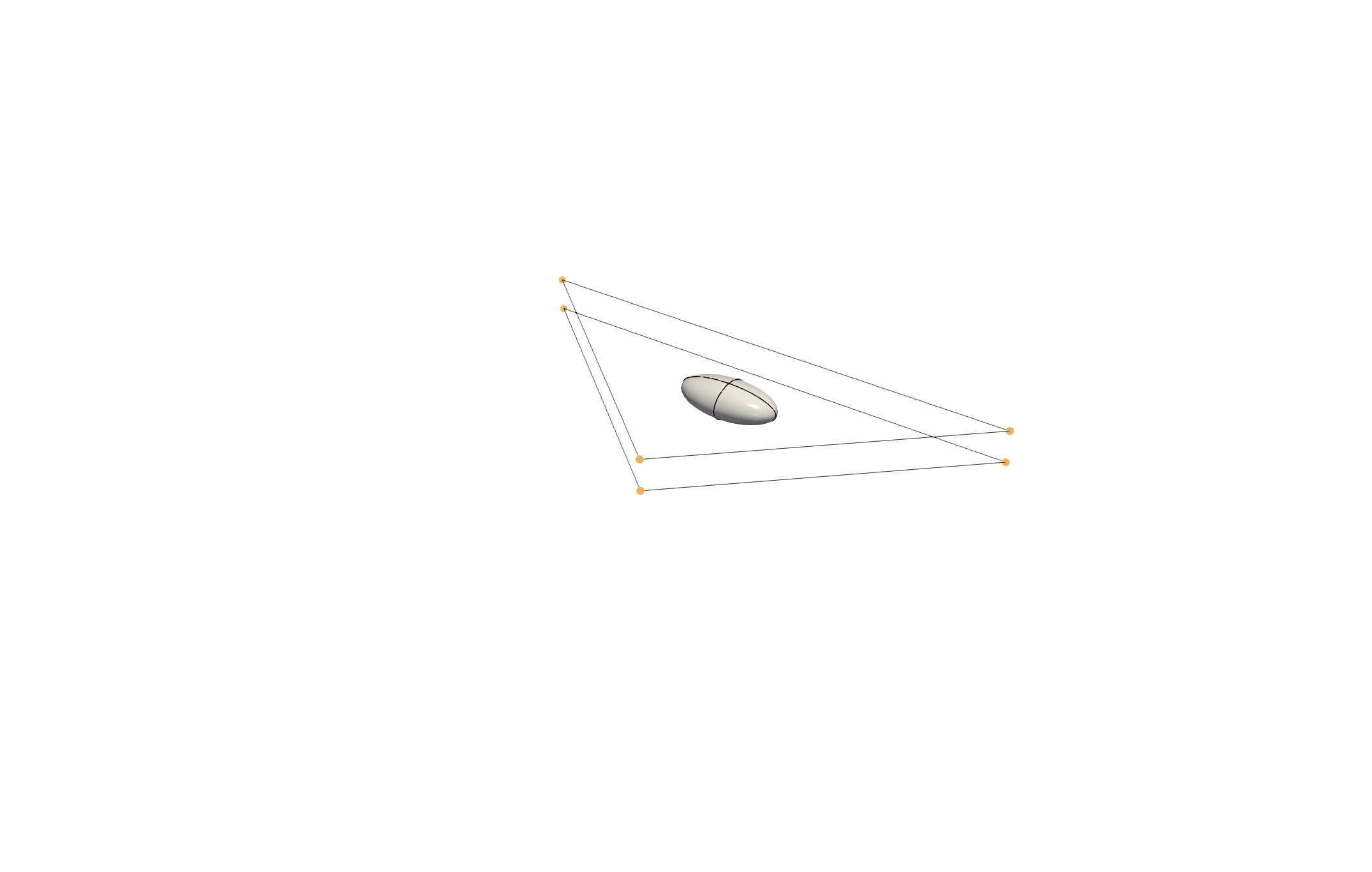}\\
         \hspace*{-0.3cm}\includegraphics[trim=140 280 350 420, clip, width=0.9\linewidth]{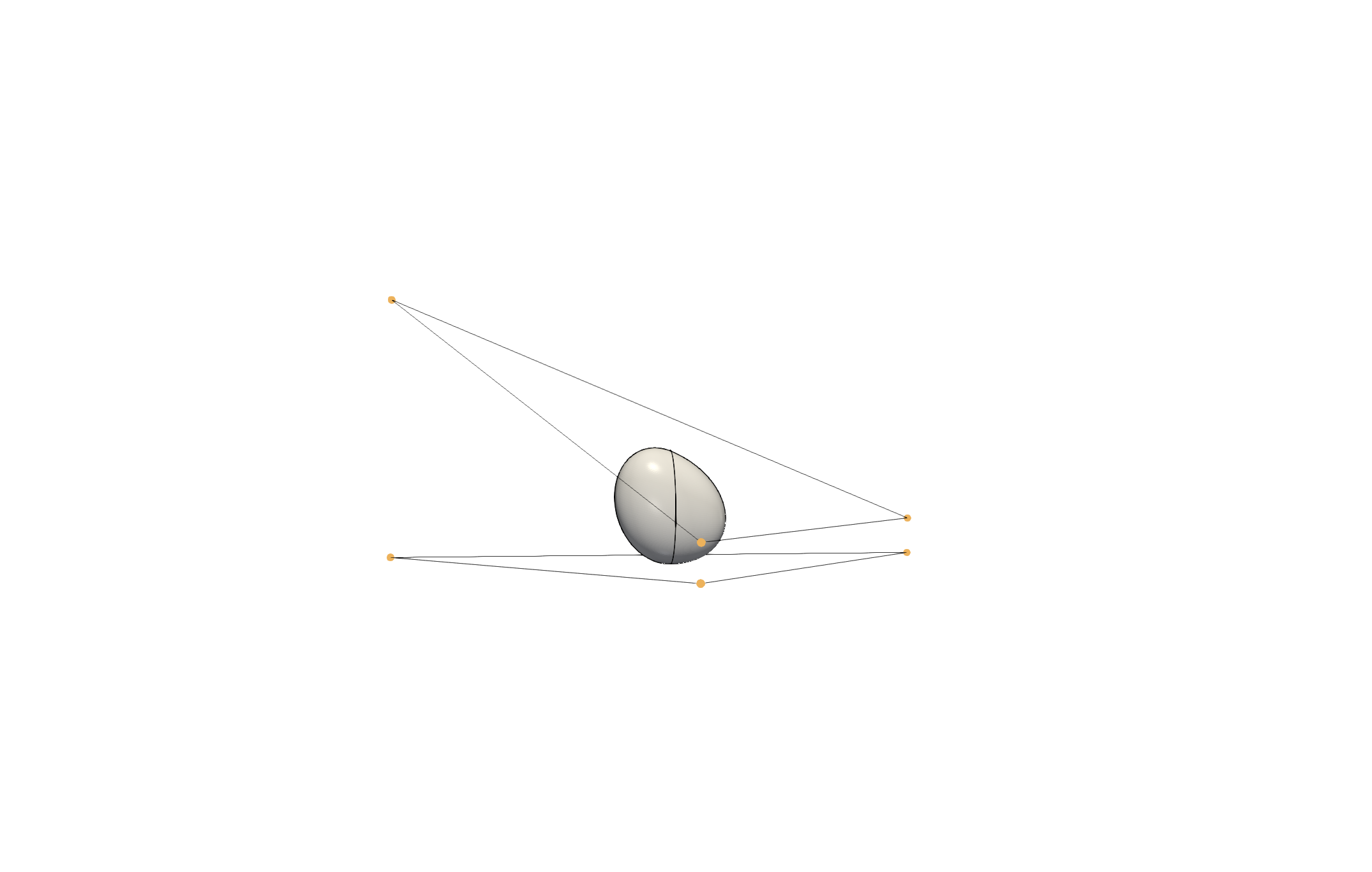}
      \end{mdframed}\vspace*{-0.4cm}
      \caption{degree $(2,3)$}
   \end{subfigure}\hspace*{0.05cm}
   \begin{subfigure}{0.32\linewidth}
      \begin{mdframed}
         \centering
         \includegraphics[trim=630 370 70 240, clip, width=0.9\linewidth]{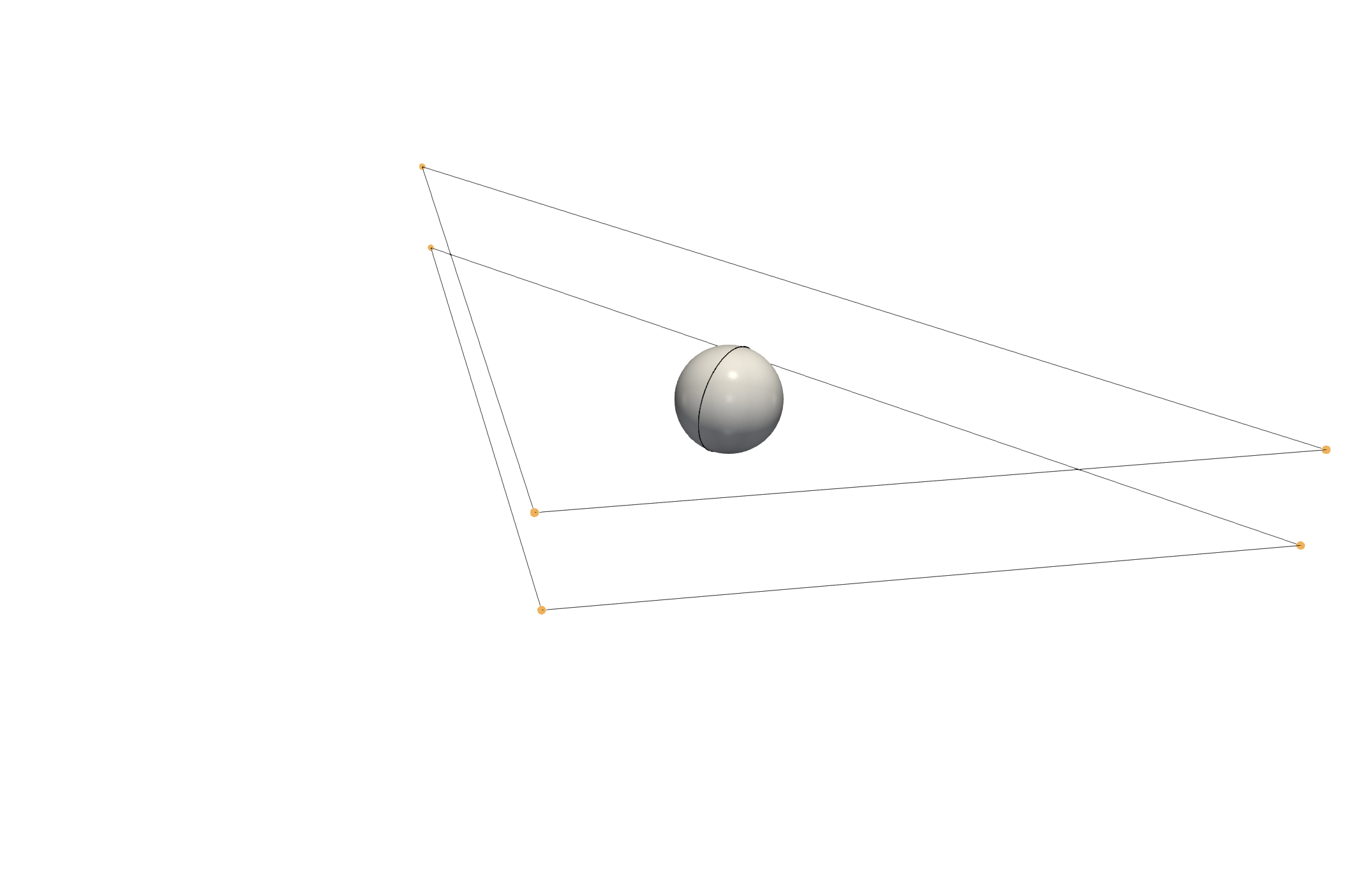}\\
         \includegraphics[trim=170 280 340 400, clip, width=0.9\linewidth]{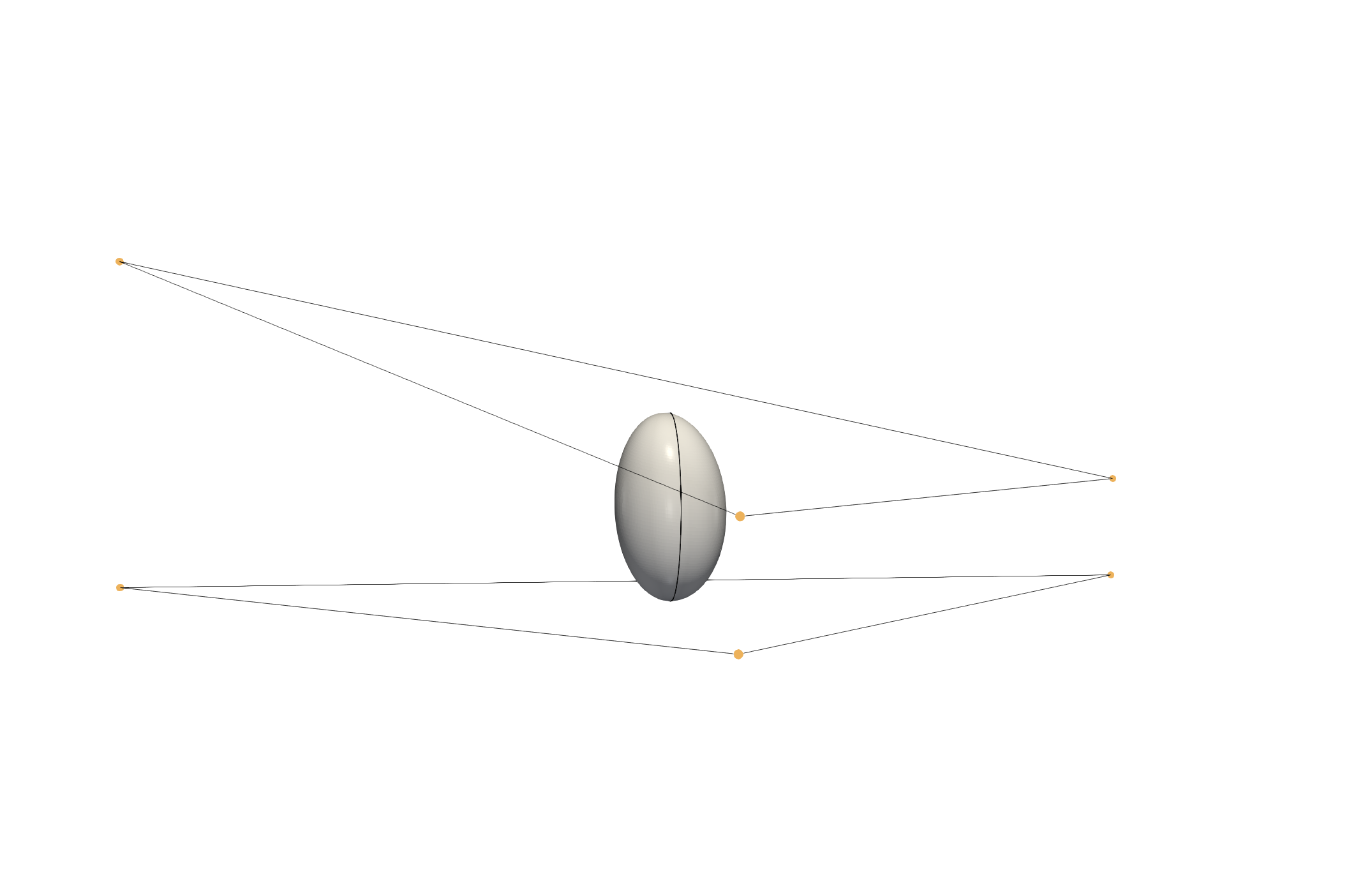}
      \end{mdframed}\vspace*{-0.4cm}
      \begin{mdframed}
         \centering
         \includegraphics[trim=600 370 130 320, clip, width=0.9\linewidth]{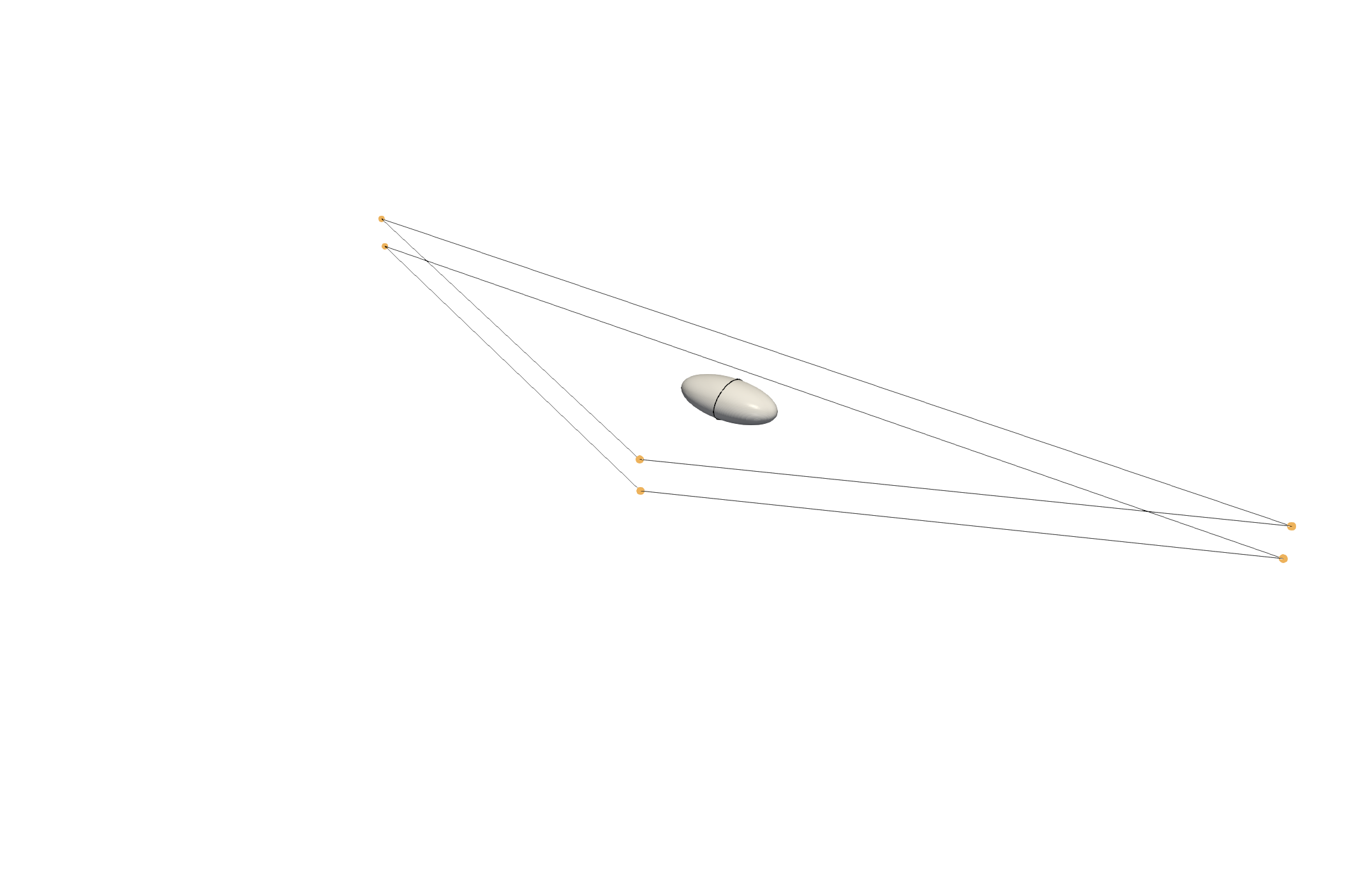}\\
         \includegraphics[trim=140 280 350 420, clip, width=0.9\linewidth]{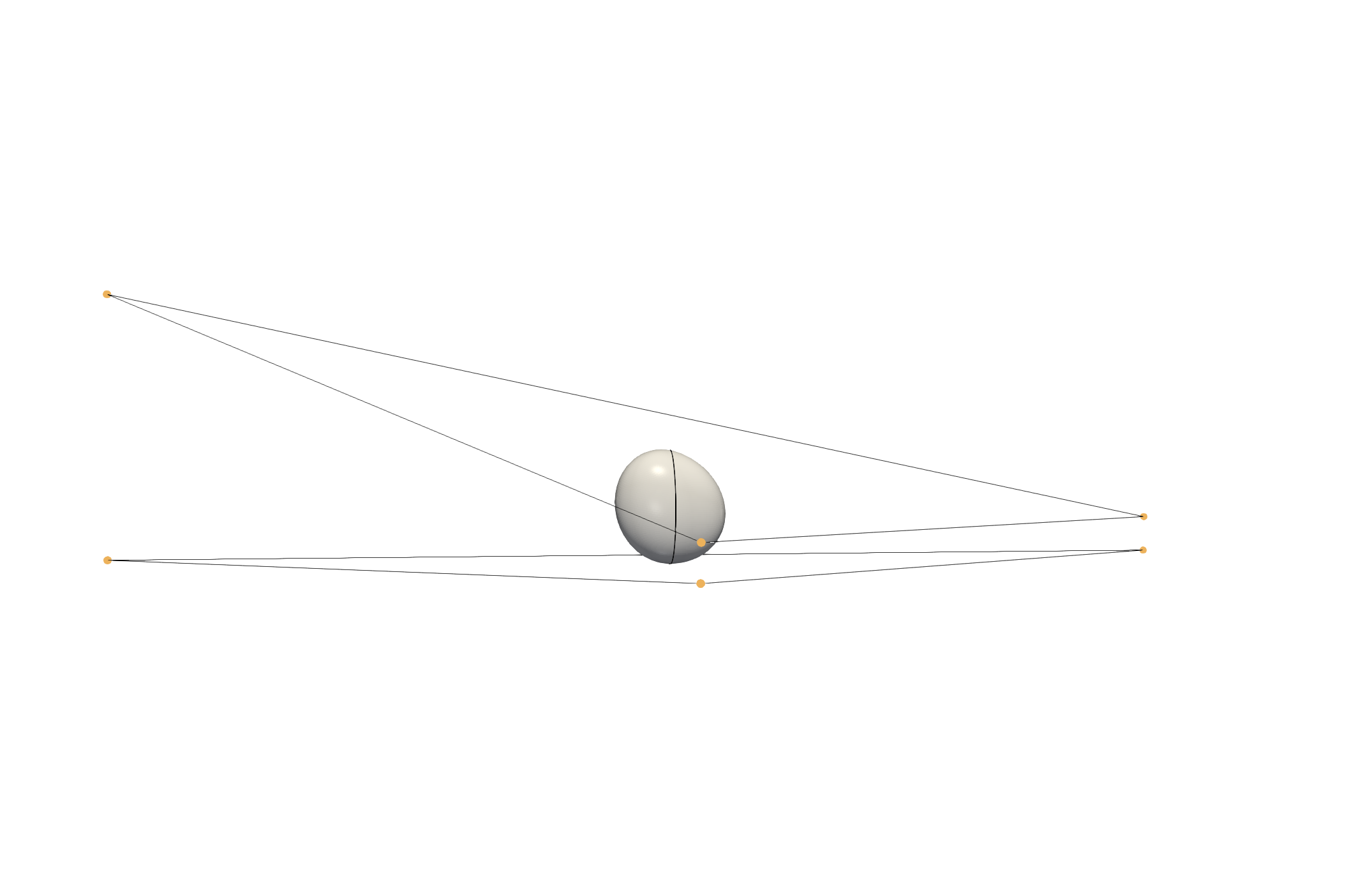}
      \end{mdframed}\vspace*{-0.4cm}
      \caption{degree $(3,3)$}
   \end{subfigure}
   \caption{Different $C^1$ smooth descriptions of unit spheres (top row boxes) and ellipsoids with axis lengths $(1, \frac{1}{2}, \frac{1}{3})$ (bottom row boxes).
   The quadrics in figures (a)--(c) are made up of $8$, $4$ and $2$ rational pieces of bi-degrees $(2,2)$, $(2,3)$ and $(3,3)$, respectively.
   In particular, the top figure in each box shows the exact quadric, while the bottom figure shows a deformed quadric obtained by raising one of the control points of the exact quadric; see Examples~\ref{ex:ellipsoid_22},~\ref{ex:ellipsoid_23} and~\ref{ex:ellipsoid_33} for details.}
   \label{fig:ellipsoid}
\end{sidewaysfigure*}

\subsubsection{$C^1$ description of degree $(2,2)$}\label{ssec:ellipsoid22}
For the first approach, we choose $\domain^\pcoorda = [0, 4]$ and $\domain^\pcoordb = [0, 2]$, and build the univariate rational spline spaces $\splSpace^\pcoorda$ (periodic) and $\splSpace^\pcoordb$ on them using the following sets of parameters:
\begin{alignat*}{6}
\splSpace^\pcoorda &:&
&\left.
\begin{array}{l}
	\seg{i}{\bknot} = [0,0,0,1,1,1]\\
	\seg{i}{\brwts} = \Bigl[1, \tfrac{\sqrt{2}}{2},1\Bigl]
\end{array}
\right\}\;
i = 1, \dots, 4;\\
\splSpace^\pcoordb &:&
&\left.
\begin{array}{l}
\seg{i}{\bknot} = [0,0,0,1,1,1]\\
\seg{i}{\brwts} = \Bigl[1, \tfrac{\sqrt{2}}{2},1\Bigr]
\end{array}
\right\}\;
i = 1, 2.
\end{alignat*}
The corresponding piecewise-NURBS extraction operators are $\extMat^\pcoorda=\periodic{\extMat}$ defined in Equation \eqref{eq:matH-ellipse_2} and $\extMat^\pcoordb$ given by
\begin{equation*}
\extMat^\pcoordb = \begin{bmatrix*}[c]
1 & 0 & 0 & 0 & 0 & 0\\
0 & 1 & \frac{1}{2} & \frac{1}{2} & 0 & 0\\
0 & 0 & \frac{1}{2} & \frac{1}{2} & 1 & 0\\
0 & 0 & 0 & 0 & 0 & 1
\end{bmatrix*},
\end{equation*}
respectively.
The full tensor-product extraction matrix is obtained as
\begin{equation} \label{eq:matH-ellipsoid}
\extMat = \extMat^\pcoorda \otimes \extMat^\pcoordb.
\end{equation}
The next ingredient in the $C^1$ smooth polar construction is the polar extraction operator, and it can be computed, using above formulas, to be
\begin{equation} \label{eq:matE-ellipsoid}
\extMatP = 
\begin{bmatrix*}[c]
\reduced{\extMatP} & \mbf{0}\\
\mbf{0} &  \exchangeMat_3\reduced{\extMatP}\exchangeMat_8
\end{bmatrix*},
\end{equation}
where the matrix $\reduced{\extMatP}$ is given by
\begin{equation*}
	\reduced{\extMatP}
	=
	\frac{1}{3} + 
	\frac{1}{\sqrt{2}}
	\begin{bmatrix*}[c]
	0 & 0 & 0 & 0 & \frac{1}{3} & -\frac{1}{3} & -\frac{1}{3} & \frac{1}{3} \\
	0 & 0 & 0 & 0 & -e_+  & -e_- & e_+ & e_-\\
	0 & 0 & 0 & 0 & e_-  & e_+ & -e_- & -e_+
	\end{bmatrix*},
\end{equation*}
and
\begin{equation*}
	e_+ = \frac{\sqrt{3}+1}{6},\quad e_- = \frac{\sqrt{3}-1}{6};
\end{equation*}
see Equations \eqref{eq:polarC1Extraction} and \eqref{eq:polarC1Extraction-double}--\eqref{eq:polarC1Extraction-double-sub}. \revision{Here the compact notation of adding a scalar to a matrix means the entrywise operation of adding the scalar to each entry of the matrix.}
The matrix $\extMatP$ maps a total of $16$ tensor-product $C^{-1}$ piecewise-NURBS $\bsp_j$ of degree $(2,2)$ to a total of $6$ $C^1$ polar splines $N_\polari$.
Equivalently, $\extMatP\extMat$ maps a total of $72$ NURBS $\bez_j$ of degree $(2,2)$ to a total of $6$ $C^1$ polar splines $N_\polari$.
These relations are encapsulated in the following equation,
\begin{equation*}
N_\polari(\pcoorda,\pcoordb) = \sum_{j = 1}^{16}\extMatPel_{\polari j}\bsp_{j}(\pcoorda,\pcoordb) = \sum_{j = 1}^{16}\extMatPel_{\polari j}\sum_{k = 1}^{72}\extMatel_{jk}\bez_{k}(\pcoorda,\pcoordb).
\end{equation*}
Finally, the required ellipsoid in Equation \eqref{eq:ellipsoid} is obtained by defining the associated $6$ control points as
\begingroup\makeatletter\def\f@size{9}\check@mathfonts
\begin{equation*}
\begin{gathered}
   \mbf{f}_1 = \left(0,2\sqrt{2}a_y,a_z\right),\ 
   \mbf{f}_2 = \left(-\sqrt{6}a_x,-\sqrt{2}a_y,a_z\right),\\
   \mbf{f}_3 = \left(\sqrt{6}a_x,-\sqrt{2}a_y,a_z\right),\ 
   \mbf{f}_4 = \left(-\sqrt{6}a_x,-\sqrt{2}a_y,-a_z\right),\\
   \mbf{f}_5 = \left(\sqrt{6}a_x,-\sqrt{2}a_y,-a_z\right),\ 
   \mbf{f}_6 = \left(0,2\sqrt{2}a_y,-a_z\right).
\end{gathered}
\end{equation*}\endgroup

\begin{example}\label{ex:ellipsoid_22}
   The box at the top in Figure~\ref{fig:ellipsoid} (a) shows a bi-degree $(2,2)$ unit sphere built by choosing $a_x = a_y = a_z = 1$, while the box at the bottom shows a bi-degree $(2,2)$ ellipse with axis lengths $(1, \frac{1}{2}, \frac{1}{3})$ built by choosing $a_x = 2a_y = 3a_z = 1$.
   These descriptions use only $8$ rational pieces.
   In each box, the figure at the top shows the exact quadric, while the figure at the bottom shows the deformed quadric obtained by perturbing the control points as per Equation \eqref{eq:conic3-perturb}.
   The exact and deformed surfaces are all $C^1$ smooth at the poles.
\end{example}

\subsubsection{$C^1$ description of degree $(2,3)$}\label{ssec:ellipsoid23}
For the second approach, we choose $\domain^\pcoorda = [0, 4]$ and $\domain^\pcoordb = [0, 1]$, and build the univariate rational spline spaces $\splSpace^\pcoorda$ (periodic) and $\splSpace^\pcoordb$ on them using the following sets of parameters:
\begin{alignat*}{6}
\splSpace^\pcoorda &:&
&\left.
\begin{array}{l}
\seg{i}{\bknot} = [0,0,0,1,1,1]\\
\seg{i}{\brwts} = \Bigl[1, \tfrac{\sqrt{2}}{2},1\Bigr]
\end{array}
\right\}\;
i = 1, \dots, 4;\\
\splSpace^\pcoordb &:&
&\begin{array}{l}
\seg{1}{\bknot} = [0,0,0,0,1,1,1,1]\\
\seg{1}{\brwts} = \Bigl[1, \tfrac{1}{3}, \tfrac{1}{3},1\Bigr]
\end{array}.
\end{alignat*}
The corresponding piecewise-NURBS extraction operators are $\extMat^\pcoorda=\periodic{\extMat}$ defined in Equation \eqref{eq:matH-ellipse_2} and $\extMat^\pcoordb=\identMat_4$, respectively.
The full tensor-product extraction matrix $\extMat$ is obtained as in Equation \eqref{eq:matH-ellipsoid}.
Moreover, the polar extraction operator $\extMatP$ is equal to the matrix in Equation \eqref{eq:matE-ellipsoid}.
The latter matrix maps a total of $16$ tensor-product $C^{-1}$ piecewise-NURBS $\bsp_j$ of degree $(2,3)$ to a total of $6$ $C^1$ polar splines $N_\polari$.
Equivalently, $\extMatP\extMat$ maps a total of $48$ tensor-product NURBS $\bez_j$ of degree $(2,3)$ to a total of $6$ $C^1$ polar splines $N_\polari$.
These relations are encapsulated in the following equation,
\begin{equation*}
N_\polari(\pcoorda,\pcoordb) = \sum_{j = 1}^{16}\extMatPel_{\polari j}\bsp_{j}(\pcoorda,\pcoordb) = \sum_{j = 1}^{16}\extMatPel_{\polari j}\sum_{k = 1}^{48}\extMatel_{jk}\bez_{k}(\pcoorda,\pcoordb).
\end{equation*}
Finally, the required ellipsoid in Equation \eqref{eq:ellipsoid} is obtained by defining the associated $6$ control points as
\begingroup\makeatletter\def\f@size{9}\check@mathfonts
\begin{equation*}
\begin{gathered}
   \mbf{f}_1 = \left(0,4\sqrt{2}a_y,a_z\right),\ 
   \mbf{f}_2 = \left(-2\sqrt{6}a_x,-2\sqrt{2}a_y,a_z\right),\\
   \mbf{f}_3 = \left(2\sqrt{6}a_x,-2\sqrt{2}a_y,a_z\right),\ 
   \mbf{f}_4 = \left(-2\sqrt{6}a_x,-2\sqrt{2}a_y,-a_z\right),\\
   \mbf{f}_5 = \left(2\sqrt{6}a_x,-2\sqrt{2}a_y,-a_z\right),\ 
   \mbf{f}_6 = \left(0,4\sqrt{2}a_y,-a_z\right).
\end{gathered}
\end{equation*}\endgroup
Note that, since we are using a higher-degree representation compared to Section~\ref{ssec:ellipsoid22}, the control points move farther away from the spline surface, mimicking the behavior of classical NURBS. 

\begin{example}\label{ex:ellipsoid_23}
   The box at the top in Figure~\ref{fig:ellipsoid} (b) shows a bi-degree $(2,3)$ unit sphere built by choosing $a_x = a_y = a_z = 1$, while the box at the bottom shows a bi-degree $(2,3)$ ellipse with axis lengths $(1, \frac{1}{2}, \frac{1}{3})$ built by choosing $a_x = 2a_y = 3a_z = 1$.
   These descriptions use $4$ rational pieces.
   In each box, the figure at the top shows the exact quadric, while the figure at the bottom shows the deformed quadric obtained by perturbing the control points as per Equation \eqref{eq:conic3-perturb}.
   The exact and deformed surfaces are all $C^1$ smooth at the poles.
\end{example}

\subsubsection{$C^1$ description of degree $(3,3)$}\label{ssec:ellipsoid33}
Finally, for the third approach, we choose $\domain^\pcoorda = [0, 2]$ and $\domain^\pcoordb = [0, 1]$, and build the univariate rational spline spaces $\splSpace^\pcoorda$ (periodic) and $\splSpace^\pcoordb$ on them using the following sets of parameters:
\begin{alignat*}{6}
\splSpace^\pcoorda &:&
&\left.
\begin{array}{l}
\seg{i}{\bknot} = [0,0,0,0,1,1,1,1]\\
\seg{i}{\brwts} = \Bigl[1, \tfrac{1}{3},\tfrac{1}{3},1\Bigr]
\end{array}
\right\}\;
i = 1, 2;\\
\splSpace^\pcoordb &:&
&\begin{array}{l}
\seg{1}{\bknot} = [0,0,0,0,1,1,1,1]\\
\seg{1}{\brwts} = \Bigl[1,\tfrac{1}{3},\tfrac{1}{3},1\Bigr]
\end{array}.
\end{alignat*}
The corresponding piecewise-NURBS extraction operators are $\extMat^\pcoorda=\periodic{\extMat}$ defined in Equation \eqref{eq:matH-ellipse_3} and $\extMat^\pcoordb=\identMat_4$, respectively.
The full tensor-product extraction matrix $\extMat$ is obtained as in Equation \eqref{eq:matH-ellipsoid}.
Moreover, the polar extraction operator $\extMatP$ is equal to the matrix in Equation \eqref{eq:matE-ellipsoid}.
The latter matrix maps a total of $16$ tensor-product $C^{-1}$ piecewise-NURBS $\bsp_j$ of degree $(3,3)$ to a total of $6$ $C^1$ polar splines $N_\polari$.
Equivalently, $\extMatP\extMat$ maps a total of $32$ tensor-product NURBS $\bez_j$ of degree $(3,3)$ to a total of $6$ $C^1$ polar splines $N_\polari$.
These relations are encapsulated in the following equation,
\begin{equation*}
N_\polari(\pcoorda,\pcoordb) = \sum_{j = 1}^{16}\extMatPel_{\polari j}\bsp_{j}(\pcoorda,\pcoordb) = \sum_{j = 1}^{16}\extMatPel_{\polari j}\sum_{k = 1}^{32}\extMatel_{jk}\bez_{k}(\pcoorda,\pcoordb).
\end{equation*}
Finally, the required ellipsoid in Equation \eqref{eq:ellipsoid} is obtained by defining the associated $6$ control points as
\begingroup\makeatletter\def\f@size{9}\check@mathfonts
\begin{equation*}
\begin{gathered}
   \mbf{f}_1 = \left(0,4\sqrt{2}a_y,a_z\right),\ 
   \mbf{f}_2 = \left(-4\sqrt{6}a_x,-2\sqrt{2}a_y,a_z\right),\\
   \mbf{f}_3 = \left(4\sqrt{6}a_x,-2\sqrt{2}a_y,a_z\right),\ 
   \mbf{f}_4 = \left(-4\sqrt{6}a_x,-2\sqrt{2}a_y,-a_z\right),\\
   \mbf{f}_5 = \left(4\sqrt{6}a_x,-2\sqrt{2}a_y,-a_z\right),\ 
   \mbf{f}_6 = \left(0,4\sqrt{2}a_y,-a_z\right).
\end{gathered}
\end{equation*}\endgroup
Observe again that, since we are using a higher-degree representation compared to Sections~\ref{ssec:ellipsoid22} and~\ref{ssec:ellipsoid23}, the control points move even farther away from the spline surface, mimicking the behavior of classical NURBS. 

\begin{example}\label{ex:ellipsoid_33}
   The box at the top in Figure~\ref{fig:ellipsoid} (c) shows a bi-degree $(3,3)$ unit sphere built by choosing $a_x = a_y = a_z = 1$, while the box at the bottom shows a bi-degree $(3,3)$ ellipse with axis lengths $(1, \frac{1}{2}, \frac{1}{3})$ built by choosing $a_x = 2a_y = 3a_z = 1$.
   These descriptions use only $2$ rational pieces.
   In each box, the figure at the top shows the exact quadric, while the figure at the bottom shows the deformed quadric obtained by perturbing the control points as per Equation \eqref{eq:conic3-perturb}.
   The exact and deformed surfaces are all $C^1$ smooth at the poles.
\end{example}

\begin{figure}
  \centering
  \includegraphics[trim=500 100 630 100, clip, width=0.4\textwidth]{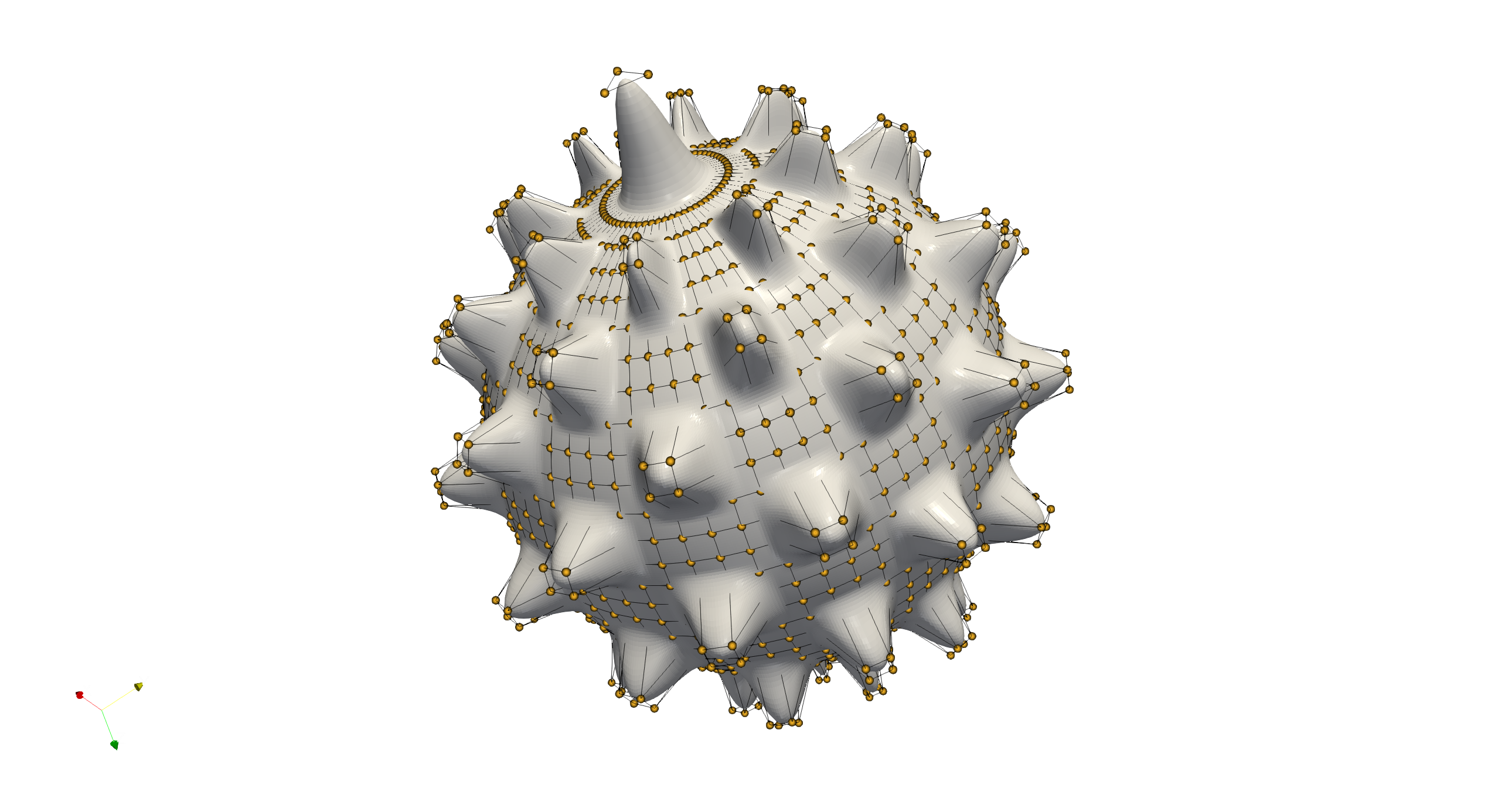} 
  \caption{A smoothly deformed sphere built from Figure~\ref{fig:ellipsoid} (a) by refining the surface and then modifying the control points.}\label{fig:sphere-deformed}
\end{figure}

\begin{remark}
   The examples presented here have focused on the simplest possible $C^1$ descriptions of quadrics, namely descriptions that either use lowest-degree splines --- bi-degree $(2,2)$ --- or the smallest number of polynomial pieces --- two.
   Unsurprisingly, these simplest descriptions can lead to large control triangles since each control point influences a large portion of the spline surface.
   Nevertheless, localized control is easily attained upon refinement (see Section \ref{sec:refinement-surfaces}) and, in particular, refinement also leads to much smaller control triangles that offer much finer geometric control.
   The surface shown in Figure~\ref{fig:sphere-deformed} illustrates this point.
   This surface has been obtained by refining and modifying the control points of the bi-degree $(2,2)$ sphere from Figure~\ref{fig:ellipsoid} (a).
   The smaller control triangle is visible near the top of the figure.
\end{remark}

	\section{Conclusions} \label{sec:conclusions}

We have presented a general class of $C^1$ smooth rational splines that allow for the construction and refinement of $C^1$ smooth curves and (polar) surfaces. They are built by gluing together multiple sets of NURBS basis functions with $C^1$ smoothness using a \DTA{} extraction matrix. The main features of the \HS{splines we have built} are the following:
\begin{itemize}
   \item all standard properties of NURBS, including support for intuitive control-point-based design,
   \item (local) degree elevation and knot insertion based on classical NURBS refinement,
   \item low-degree $C^1$ descriptions of exact ellipses and ellipsoids, and
   \item compatibility with CAD or CAE software through the explicit representation in terms of NURBS.
\end{itemize}
In particular, with regard to the last two bullets above, \revision{we believe} that the explicit, NURBS-compatible $C^1$ descriptions of ellipses and ellipsoids provided herein will be of \revision{use to geometric modellers \cite{karvciauskas2020polar} and computational scientists \cite{toshniwal2020discretedifferential} alike. For instance, the exact $C^1$ (re)parameterizations at polar points may make the design of algorithms more stable and efficient; it may also avoid the need for special treatment of polar points.}

	\section*{Acknowledgements}
	H. Speleers was partially supported by the Beyond Borders Program of the University of Rome Tor Vergata through the project ASTRID (CUP E84I19002250005) and by the MIUR Excellence Department Project awarded to the Department of Mathematics, University of Rome Tor Vergata (CUP E83C18000100006). He is a member of Gruppo Nazionale per il Calcolo Scientifico, Istituto Nazionale di Alta Matematica. 
	
	\bibliographystyle{plain}
	\bibliography{./sections_2col_rev/bibliography}

\end{document}